\documentclass[11pt]{article}
\usepackage{fancyhdr}
\usepackage{isomath}
\usepackage{mathtools} % asmmath included here
\usepackage{amsbsy}
\usepackage{amssymb}
\usepackage{amscd,amsfonts}
\usepackage{bigints}
\usepackage{graphicx}
\usepackage{verbatim}
\usepackage{euscript}
\usepackage{alltt}
\usepackage{stmaryrd}
\usepackage{relsize}
\usepackage{enumerate}
\usepackage{url}
\usepackage[stable]{footmisc}
\usepackage{breakurl}
\usepackage{hyperref}
\usepackage{comment}
\usepackage[numbers]{natbib}
\usepackage[font=small,labelfont=bf]{caption}
\usepackage[font=small,labelfont=bf]{subcaption}
\usepackage{float}
\usepackage{mwe}
\usepackage{algorithm}
\usepackage{algpseudocode}
\usepackage{xcolor}
\usepackage[super]{nth}
\usepackage{soul}
\usepackage{centernot}
\usepackage{bm}

\DeclareGraphicsExtensions{.eps,.pdf}
\usepackage{amsmath}
\usepackage{amsbsy}
\usepackage{amssymb}
\usepackage{amscd}
\usepackage{amsfonts}

\newcommand{\bfa}{{\mathbold a}}

\newcommand{\bfe}{{\mathbold e}}
\newcommand{\bff}{{\mathbold f}}
\newcommand{\bfg}{{\mathbold g}}

\newcommand{\bfv}{{\mathbold v}}

\newcommand{\bfI}{{\mathbold I}}
\newcommand{\bfJ}{{\mathbold J}}

\newcommand{\bfL}{{\mathbold L}}

\newcommand{\bfN}{{\mathbold N}}

\newcommand{\bfR}{{\mathbold R}}

\newcommand{\bfW}{{\mathbold W}}

\newcommand{\beq}{\begin{equation}}
\newcommand{\eeq}{\end{equation}}
\newcommand{\beqs}{\begin{eqnarray}}
\newcommand{\eeqs}{\end{eqnarray}}
\newcommand{\beql}{\begin{equation} \label}

\newcommand{\bfomega}{\mathbold{\omega}}

\newcommand{\bflambda}{\mathbold{\lambda}}

\newcommand{\bfOmega}{\mathbold{\Omega}}

\newcommand{\bfzero}{\mathbf{0}}

\newcommand{\parderiv}[2]{\frac{\partial #1}{\partial #2}}
\newcommand{\deriv}[2]{\frac{d #1}{d #2}}

\newcommand{\veps}{\varepsilon}

\newcommand{\p}{\partial}

\usepackage[margin=1in]{geometry}

\newcommand{\dee}{\mathcal{D}}
\newcommand{\scl}{\mathcal{L}}

\date{}
\begin{document}
\title{Hidden convexity in the heat, linear transport, and Euler's rigid body equations: A computational approach}

\author{Uditnarayan Kouskiya\thanks{Department of Civil \& Environmental Engineering, Carnegie Mellon University, Pittsburgh, PA 15213, email: udk@andrew.cmu.edu.} $\qquad$ Amit Acharya\thanks{Department of Civil \& Environmental Engineering, and Center for Nonlinear Analysis, Carnegie Mellon University, Pittsburgh, PA 15213, email: acharyaamit@cmu.edu.}}

\maketitle
\begin{abstract}
\noindent 
A finite element based computational scheme is developed and employed to assess  a duality based variational approach to the solution of the linear heat and transport PDE in one space dimension and time, and the nonlinear system of ODEs of Euler for the rotation of a rigid body about a fixed point. The formulation turns initial-(boundary) value problems into degenerate elliptic boundary value problems in (space)-time domains representing the Euler-Lagrange equations of suitably designed dual functionals in each of the above problems. We demonstrate reasonable success in approximating solutions of this range of parabolic, hyperbolic, and ODE primal problems, which includes energy dissipation as well as conservation, by a unified dual strategy lending itself to a variational formulation. The scheme naturally associates a family of dual solutions to a unique primal solution; such `gauge invariance' is demonstrated in our computed solutions of the heat and transport equations, including the case of a transient dual solution corresponding to a steady primal solution of the heat equation. Primal evolution problems with causality are shown to be correctly approximated by non-causal dual problems.
\end{abstract}

\section{Introduction}
The goal of this paper is to test a duality based technique proposed in \cite{action_2, action_3} for solving differential equations. We implement the proposed scheme on the linear heat and first-order transport equations posed in bounded domains in one space dimension, and the nonlinear system of ODEs of Euler for the motion of a rigid body with a fixed point (which applies no torque to the body). Euler's system comprises three coupled ODEs for the components of the angular velocity vector of the rigid body on the rotating principal axis frame of its moment of inertia tensor. In all cases, the initial value problems are converted to a boundary value problem in time by the scheme. We compute and assess approximate solutions to the dual equations by the finite element method and verify our procedures with exact or approximate solutions of the primal equations, computed directly from the primal descriptions. Obviously, our effort is \textit{not} meant to be interpreted as the development of competing algorithms for more standard procedures of approximation for the test problems; instead, it is meant to explore the duality scheme in simple settings where the outcomes are clearly understood. The potential utility of the duality scheme is in generating dual variational principles for fairly general classes of nonlinear systems of algebraic, ordinary, and partial differential equations (the primal equations), see \cite{action_2, action_3, brenier_book, brenier2018initial}. Solutions to the Euler-Lagrange equations of the dual variational principle correspond to those of the primal equations in a well-defined, albeit formal, sense. Specifically, the scheme generates an adapted change of variables which constitutes a mapping of a dual to a primal solution. Such a scheme is most useful for primal equations for which standard definitions of a solution and methods for their approximation are not available, e.g. systems of second-order Hamilton-Jacobi equations. In the context of partial differential equations, our scheme also seems to have the unifying property of producing at worst degenerate elliptic dual problems, with the unusual feature of involving \textit{oblique} natural boundary conditions \cite{oblique}.

The outline of the paper is as follows: The dual formulation of the heat equation is presented in Sec.~2.  In Sec.~3, a dual formulation of the linear transport equation is developed. Sec.~4 describes the dual formulation of Euler's system of ODE for the motion of a rigid body about a fixed point. Sec.~5 describes the finite element implementations of the above problems and presents computational results. Sec.~6 contains concluding remarks.

\section{Dual formulation of the heat equation}
Following the ideas in \cite{action_2, action_3} reviewed in Appendix \ref{app:rev_duality} for the convenience of the reader, we write the heat equation in first order form:
\begin{equation}\label{eq:primal_ht}
    \p_x \theta - \pi = 0; \qquad \p_t \theta - \p_x (k \pi) = 0; \qquad (x,t)  \in \Omega = (0,L) \times (0,T).
\end{equation}
Using the Lagrange multiplier fields $p,l$ for the two equations, define 
\begin{equation*}
    \widehat{S}[ \theta, l, \pi, p] = \int_\Omega dtdx \,   - \theta \p_x p - p \pi - \theta \p_t l + k \pi \p_x l + H(\theta, \pi),
\end{equation*}
where the requirement on the function $H$ is that, defining $\dee: = (\p_x p , p, \p_t l, \p_x l)$, the following equations are solvable for $\theta, \pi$ in terms of $\dee$:
\begin{equation}\label{eq:L_trans}
    \begin{aligned}
          - \p_x p - \p_t l + \p_\theta H & = 0\\
          - p + k \p_x l + \p_\pi H & = 0.
    \end{aligned}
\end{equation}
Choosing 
\begin{equation}\label{eq:heat_H}
   H(\theta, \pi) = \frac{1}{2} (\theta^2 + \pi^2) 
\end{equation} 
we obtain the following \textit{dual-to-primal mapping} (DtP)
\begin{equation}\label{eq:DtP}
    \theta = \p_x p + \p_t l; \qquad \pi = p - k \p_x l.
\end{equation}
Now eliminate the fields $\theta, \pi$ in $\widehat{S}$ by substituting these functions for $\theta$ and $\pi$ to define the \textit{dual} functional
\begin{equation}\label{eq:pre_S}
   - \frac{1}{2} \int_\Omega  dtdx \, (\p_x p + \p_t l)^2 + (p - k \p_x l)^2.
\end{equation}
The requirement \eqref{eq:L_trans} ensures that the dual E-L system of \eqref{eq:pre_S} is simply the primal system \eqref{eq:primal_ht} (cf. Appendix \ref{app:rev_duality} for more details), in which the primal fields $(\theta, \pi)$ are expressed in terms of the dual fields through the DtP mapping \eqref{eq:DtP}.

Denoting $(p,l)$ as the dual fields corresponding to the heat equation \eqref{eq:primal_ht}, the `dual heat PDE' system for the quadratic $H$ chosen in \eqref{eq:heat_H} is the following constant coefficient second order system of PDE
\begin{equation}\label{eq:dual_heat_pde}
\begin{aligned}
      \p_x^2 p  + \p_{xt} l  + k \p_x l  - p & = 0,  \qquad k > 0 \ \mbox{a constant} \\
      k^2 \p_x^2 l   + \p_t^2 l     +  \p_{tx} p  - k \p_x p   & = 0.
\end{aligned}
\end{equation}
The quadratic form arising from its principal part is given by
\begin{equation*}
    Q_{dual heat} := (\p_x p + \p_t l)^2 + k^2 (\p_x l)^2
\end{equation*}
which shows that the system is \textit{degenerate elliptic} since $Q_{dualheat}$ is non-negative for all `gradient' matrices $F$ with $F_{11} = \p_x p, F_{12} = \p_t p, F_{21} = \p_x l, F_{22} = \p_t l$, and ellipticity fails on the rank-one direction $(a, 0) \otimes (0, 1)$, for any value of $a$.

\subsection{Formal uniqueness of solutions of the dual heat equation}\label{heat_unique}
We utilize uniqueness arguments to infer boundary conditions on space-time domains to define well-set dual problems.

The first variation of the dual functional \eqref{eq:pre_S} implies that the first dual equation arises from variations $\delta p$ and the second one from variations $\delta l$. Without imposing any boundary conditions, the following statement holds, assuming the dual E-L equations are satisfied:
\begin{equation}\label{eq:weak_dual_heat}
    \begin{aligned}
      - &  \int_0^T dt \int_0^L dx \, \p_x \delta p \,(\p_x p + \p_t l) - \int_0^T dt \int_0^L dx \, \delta p \,(p - k \p_x l) \\
      - & \int_0^T dt \int_0^L dx \, \p_t \delta l\, ( \p_x p + \p_t l) + \int_0^T dt \int_0^L dx \, \p_x \delta l \, (k (p - k \p_xl)) \\
      +  & \int_0^T dt \, \delta p (\p_x p + \p_t l)\big|^L_0 + \int^L_0 dx \, \delta l \, (\p_x p + \p_t l) \big|^T_0 - \int_0^T dt \, \delta l (k (p - k \p_x l)) \big|^L_0 = 0. 
    \end{aligned}
\end{equation}

Let $(p,l)$ be the difference of two (sufficiently smooth, as required) solutions to the system \eqref{eq:dual_heat_pde} and let each of the solutions $(p^i, l^i)$, $i = 1,2$ satisfy 
\begin{enumerate}
\item either $(\p_x p^i + \p_t l^i)$ or $p^i$ equal specified values on the left and right boundaries of the space-time domain;
\item $(\p_x p^i + \p_t l^i)$ equals a specified initial condition on $\theta$ on the bottom boundary, $t = 0$, of the domain;
\item $l^i$ equals an arbitrarily specified function of $x$ (w.l.o.g. chosen to  be 0) on the top boundary of the domain, $t = T$;
\item either $(p^i - k \p_x l^i)$ or $l^i$ equal specified values on the left and right boundaries of the space-time domain.
\item on any point of the left and right boundaries, either $(\p_x p^i + \p_t l^i)$ or $(p^i - k \p_x l^i)$ are specified.
\end{enumerate}
Then, choosing $\delta p = p$ and $\delta l = l$ in \eqref{eq:weak_dual_heat} results in all the boundary terms vanishing and
\begin{equation*}
- \int_0^T dt \int_0^L dx \, (\p_x p + \p_t l)^2 + (p - k \p_xl)^2 = 0
\end{equation*}
which implies that
\begin{equation}\label{eq:heat_unique}
    \p_x p + \p_t l = 0 \ \mbox{and} \  p - k \p_xl = 0  \Longrightarrow \p_x p = k \p_{xx} l \Longrightarrow  \p_tl + k \p_{xx} l = 0   \qquad \mbox{in} \ (0,L) \times (0,T),
\end{equation}
which is the backward heat equation in $l$ for $k > 0$. Since at any left or right boundary point either $l^i$ or $(p^i - k \p_x l^i)$ is specified, the latter implying that at such points $p^i$ is specified (by b.c.s 4.,5.,1. above), \eqref{eq:heat_unique} (and 4.) imply that either $l = 0$ or $\p_x l = 0$ at such points. Then, defining $\Lambda(t) := \frac{1}{2} \int_0^L dx \,  l^2 (x,t) \geq 0$, \eqref{eq:heat_unique} gives $\p_t \Lambda(t) \geq 0 \Longrightarrow \Lambda(T) - \Lambda(t) \geq 0$ for all $t \in [0,T]$. But by the `top' boundary condition $\Lambda(T) = 0$, so that $0 \leq \Lambda(t) \leq 0$, and we have $l = 0$ almost everywhere in the space-time domain, and since $p - k \p_xl = 0$ in the same domain, we have uniqueness.

\subsection{Weak formulation for the dual heat equation}
For the sake of definiteness, let us assume that the primal heat equation is posed with the following initial and boundary conditions:
\begin{equation}\label{eq:primal_bc}
    \theta(x,0) = \theta_0(x); \qquad \theta(0,t) = \theta_l(t); \qquad \pi(L,t) = \pi_r(t).
\end{equation}
Then the following weak form for the problem suffices to compute unique (approximate) solutions to the dual problem, from which the unique primal solution to \eqref{eq:primal_ht} and \eqref{eq:primal_bc} can be approximated through the DtP mapping (for arbitrarily specified $T$):
\begin{equation}\label{eq:weak_heat}
     \begin{aligned}
     -  &  \int_0^T dt \int_0^L dx \, \p_x \delta p \,(\p_x p + \p_t l) - \int_0^T dt \int_0^L dx \, \delta p \,(p - k \p_x l) \\
      - & \int_0^T dt \int_0^L dx \, \p_t \delta l\, ( \p_x p + \p_t l) + \int_0^T dt \int_0^L dx \, \p_x \delta l \, (k (p - k \p_xl)) \\
      - & \int_0^T dt \, \delta p(0,t) \theta_l(t) - \int^L_0 dx \, \delta l(x,0) \, \theta_0(x) - \int_0^T dt \, \delta l(L,t) (k \pi_r(t)) = 0; \\
      & \delta l (x,T) = 0; \qquad \delta p (L,t) = 0; \qquad \delta l (0,t) = 0;\\
      & l(x,T) = l_T(x); \qquad p(L,t) = p_r(t); \qquad l(0,t) = l_l(t)
    \end{aligned}
\end{equation}
 with  $l_T(x), p_r(t), l_l(t)$ being arbitrarily specified continuous functions subject to $l_T(0) = l_l(T)$. This weak statement corresponds to the first variation of the dual functional
 \begin{equation*}
 \begin{aligned}
     S[l,p] & = - \frac{1}{2} \int_\Omega dtdx  \ (\p_x p + \p_t l)^2 + (p - k \p_x l)^2 \\
     & \quad - \int_0^T dt \ p(0,t) \theta_l(t)  - \int^L_0 dx \ l(x,0) \, \theta_0(x)  - \int_0^T dt \ l(L,t) (k \pi_r(t)).
     \end{aligned}
 \end{equation*}
If a Dirichlet boundary condition of the form $\theta(L,t)=\theta_r(t)$ is applied to the right boundary instead of $\pi(L,t)=\pi_r(t)$, the formulation and weak form of the problem remain the same. The only difference is that instead of $p(L,t)$, the dual field $l(L,t)$ is specified arbitrarily as $l_r(t)$ with $\delta l (L,t) = 0$ and the boundary term in \eqref{eq:weak_heat} for $\pi_r$ on the right boundary is replaced by the boundary term produced based on $\theta_r(t)$, following what is done for the Dirichlet condition for the left boundary.
\section{Dual formulation of the linear transport equation}
An exactly similar procedure as for the heat equation can be followed for the linear, constant coefficient wave equation 
\begin{equation}\label{eq:primal_wave}
    \p_t u + c \, \p_x u = 0 \  \mbox{in} \  (0,L) \times (0,T), \quad c > 0; \quad u(0,t) = u_l(t); \quad u(x, 0) = u_0(x)
\end{equation}
(with the functions $u_l, u_0$ being specified) to obtain its dual functional with natural boundary conditions
\begin{equation}\label{eq:wave_dual}
  S[\lambda] = - \int_0^T dt \int_0^L dx \, \frac{1}{2} (\p_t \lambda + c \p_x \lambda)^2 - \int_0^L dx \, \lambda(x, 0) u_0(x) - \int_0^T dt \, \lambda(0,t) u_l(t) c
\end{equation}
and dual PDE
\begin{equation}\label{eq:dual_wave_EL}
    \p_t (\p_t \lambda + c \p_x \lambda) + c \, \p_x (\p_t \lambda + c \p_x \lambda) = 0 = \p_t^2 \lambda + 2c \, \p_{tx} \lambda + c^2 \, \p_x^2 \lambda.
\end{equation}
Here, 
\[H(u) = \frac{1}{2} u^2
\]
and the DtP mapping emerges as
\begin{equation}\label{eq:DtP_wave}
    u = \p_t \lambda + c \p_x \lambda.
\end{equation}
We supplement the dual PDE with the essential boundary conditions 
\begin{equation*}
\lambda(x, T) = 0; \qquad \lambda(L, t) = 0,
\end{equation*}
where the specifications of these functions are `arbitrary,' as concerns recovery of the (formally) unique primal solution by the scheme.

As before, the quadratic form in gradients is
\begin{equation*}
    Q_{dual1-wave} = (\p_t \lambda + c \p_x \lambda)^2,
\end{equation*}
and this is again positive semi-definite for all gradient vectors $(\p_x \lambda, \p_t \lambda)$ with ellipticity failing for the vector direction $a \, (1, -c)$, for any value of $a$.

Equivalently,
\[ 
grad ({A \, grad \lambda}) =
\begin{bmatrix} \p_x
\left( \begin{bmatrix}
    c^2 & c \\
    c & 1
\end{bmatrix}
\begin{bmatrix}
    \p_x \lambda \\
    \p_t \lambda
\end{bmatrix} \right)
& & 
\p_t
\left( \begin{bmatrix}
    c^2 & c \\
    c & 1
\end{bmatrix}
\begin{bmatrix}
    \p_x \lambda \\
    \p_t \lambda
\end{bmatrix} \right)
\end{bmatrix}
\]
so that $div (A \, grad \lambda)  = 0$ is \eqref{eq:dual_wave_EL}, and ellipticity is equivalent to $n_i A_{ij} n_j \neq 0$ for all $n \neq 0$, which does not hold for the $n$ given.

Weak solutions of this dual system \eqref{eq:dual_wave_EL} obviously exist, generated explicitly by first integrating the primal equation and boundary and initial conditions by the method of characteristics to obtain, say $u(x,t)$, and then again using the method of characteristics utilizing the DtP mapping, $u(x,t) = \p_t \lambda + c \p_x \lambda$ as the governing equation, and using the `top' and `right' boundary conditions on $\lambda$. Clearly, the smoothness of the solution depends on the smoothness of the specified data $u_l, u_0$ (with the arbitrary data on $\lambda$ always chosen to be smooth, w.l.o.g.).

The degenerate elliptic dual problem allows gradient-discontinuities which translate, through the DtP mapping, to strong discontinuities in the primal $u$. Curiously, the argument in the previous paragraph indicates that  functions $\lambda$ that are discontinuous across characteristics can be constructed as solutions of the method of characteristics, for appropriately chosen data on $u$. Whether such functions can qualify as solutions to the E-L equations of the functional in \eqref{eq:wave_dual} is a question that remains to be answered.
\subsection{Formal uniqueness of solutions of the dual linear transport equation}
Let $\lambda$ be the difference of two solutions $\lambda^i$, $i = 1,2$ of \eqref{eq:dual_wave_EL}, each satisfying the following boundary conditions:
\begin{itemize}
\item $(\p_t \lambda^i + c \, \p_x \lambda^i)$ equals the specified functions $u_l$ on $x = 0$ and $u_0$ on $t = 0$;
\item $\lambda^i(x,T) = \lambda_T(x)$, $\lambda^i(L,t) = \lambda_r(t)$, where $\lambda_T, \lambda_r$ are specified functions.
\end{itemize}
Then, a direct consequence of the `energy method' (i.e., multiply by $\lambda$, integrate by parts, apply boundary conditions on the space-time domain) is that
\[
\p_t \lambda + c \, \p_x \lambda = 0 \qquad \mbox{in}  \ (0,L) \times (0,T).
\]
Applying the energy method again to the statement above,  and defining $\Lambda(t) := \frac{1}{2} \int_0^L dx \,  \lambda^2 (x,t) \geq 0$, one obtains $\p_t \Lambda (t) + \frac{1}{2} c (\lambda^2(L,t) - \lambda^2(0,t)) = 0$ which, on applying the b.c. at $x = L$, integrating in time from $t$ to $T$, and applying the b.c. at $t = T$ gives $  0 \geq - \Lambda(t) = \frac{1}{2} \int^T_t dt \, c \, \lambda^2(0,t) \geq 0$ since $c > 0$, and we have $\lambda(x,t) = 0$ almost everywhere in $(0,L) \times (0,T)$ and uniqueness.

\subsection{Weak formulation for the dual linear transport equation}
For the primal first-order wave equation along with the applied boundary conditions stated in \eqref{eq:primal_wave}, we obtain the weak form described as follows. For any field $\delta \lambda$ satisfying the conditions stated below, we intend to find the dual field $\lambda$ which satisfies the following equations:
\begin{equation} \label{eq:weak_dual_wave}
     \begin{aligned}
     -  &  \int_0^T dt \int_0^L dx \, \p_t \delta \lambda\, (\lambda_t + c \lambda_x) - \int_0^T dt \int_0^L dx \,
     c\,\p_x \delta\lambda\, (\lambda_t + c \lambda_x)\\
      & - \int_0^T dt \, \lambda(0,t) u_l(t) - \int_0^L dx \, \lambda(x,0) u_0(x) = 0; \\
      & \delta \lambda(x,T) = 0; \qquad \delta \lambda(L,t) = 0; \\
      & \lambda(x,T) = \lambda_T(T); \qquad
      \lambda(L,t) = \lambda_r(t),
    \end{aligned}
\end{equation}
where $\lambda_T(\cdot)$ and $\lambda_r(\cdot)$ are arbitrarily specified functions satisfying $\lambda_T(L) = \lambda_r(T)$.
The dual scheme guarantees that the solution to \eqref{eq:weak_dual_wave} implies the solution to the set of equations \eqref{eq:primal_wave}. We make use of the above weak form to compute an approximate solution for the dual field and utilize the DtP mapping \eqref{eq:DtP_wave} to obtain the corresponding field for the primal problem i.e.~the first-order wave equation. The weak form corresponds to the first variation of the dual functional \eqref{eq:wave_dual}.

\section{Dual variational principle for Euler's system for motion of a rigid body, with and without viscosity}\label{sec:euler_formulation}
Following the derivation in Appendix \ref{app} of Euler's system for the components of the angular velocity of a rigid body rotating about a fixed point under the possible action of a viscous damping force proportional to its angular momentum, we consider the following primal system of ODEs
\begin{equation}\label{eq:Euler}
    I_i \, \dot{\omega}_i + c_i \, \omega_{i+1} \omega_{i+2} + \nu I_i \,\omega_i = 0, \qquad i = 1,2,3, \ \ c_i = I_{i+2} - I_{i+1}, \ \ \nu \geq 0,
\end{equation}
where a superscript dot represents a derivative w.r.t time and \textit{all subscript indices in this Section are modulo 3}. The system is subject to the initial condition
\[
\omega_i(0) = \omega_i^0.
\]
Using the auxiliary potential
\[
H(\omega, t) = \frac{1}{2} \sum_{i=1}^3 a (\omega_i - \tilde{\omega}_i(t))^2,
\]
where $\tilde{\omega_i}$ is an arbitrarily chosen `base state' subject only to facilitating a dual solution \cite[Sec.~5]{action_3}, denoting the dual functions by $\lambda_i$, and using our standard protocol for generating a dual variational principle, one obtains the DtP mapping for this problem as
\begin{equation}\label{eq:euler_DtP}
 \omega_i - \tilde{\omega}_i =  \sum_j \mathbb{K}^{-1}_{ij} \left(I_j \dot{\lambda}_j - \nu I_j \lambda_j \right),
\end{equation}
with
\begin{equation*}
\mathbb{K} = 
\begin{bmatrix}
a & c_3 \lambda_3 & c_2 \lambda_2 \\
c_3 \lambda_3 & a & c_1 \lambda_1 \\
c_2 \lambda_2 & c_1 \lambda_1 & a
\end{bmatrix},
\end{equation*}
where $a$ can a constant of arbitrary magnitude and we choose $a = 1$.

The dual problem is obtained by substituting the change of variables \eqref{eq:euler_DtP} into the primal system \eqref{eq:Euler} appended with the boundary conditions
\begin{equation*}
    \label{eq:dual_euler_bc}
    \omega^0_i(0) = \tilde{\omega}_i(0) + \sum_j \mathbb{K}^{-1}_{ij}\left (I_j \dot{\lambda}_j (0) - \nu I_j \lambda_j (0) \right); \qquad \lambda_i(T) = 0.
\end{equation*}

The dual problem so obtained is, by design, the E-L equations of the dual functional
\begin{equation}\label{eq:euler_dual_functional}
\begin{aligned}
     S[\lambda] &:= \int^T_0  dt \, \sum_i \left( - I_i \, \omega_i \, \dot{\lambda}_i + \lambda_i (c_i \, \omega_{i+1} \omega_{i+2} + \nu I_i \omega_i ) \right) \\
    & \qquad - \sum_i I_i \, \omega^0_i \lambda_i (0)  + H(\omega, t),
\end{aligned}
\end{equation}
\textit{with all occurrences of $\omega_i$ replaced by its DtP representation \eqref{eq:euler_DtP} in terms of $(\tilde{\omega}_i, \dot{\lambda}, \lambda)$}.

It is important to note that the dual variational problem defined by the choice of base state $\tilde{\omega}_i$ given by a solution of the primal initial value problem (ivp) guarantees the existence of a solution to the dual problem given by $\lambda_i = 0$. Thus, it may be expected that for base states that are `close' to actual solutions of the primal problem, the dual problem has a solution.
\subsection{Weak formulation of the dual Euler rigid body system}
For the coupled system of equations \eqref{eq:Euler}, we employ the following weak form to obtain the approximate solution for the dual fields and subsequently use the DtP mapping \eqref{eq:euler_DtP} to evaluate the primal field variables:

\begin{equation}\label{eq:weak_euler}
\begin{gathered}
        \int_0^T dt \, \sum_i (-I_i \, \omega_i \, \dot{ \delta\lambda_i} + c_i \, \omega_{i+1} \, \omega_{i+2} \, \delta\lambda_i + \nu I_i \, \omega_i \, \delta\lambda_i) - \sum_i I_i\,\omega_i^0\,\delta\lambda_i(0) = 0; \\
        \delta\lambda_i(T) = 0; \qquad \lambda_i(T) = \lambda_i^T,
\end{gathered}
\end{equation}
where $\lambda_i^T$ are specified arbitrarily and each of the $\omega_i$ are understood to be expressed in terms of the dual fields and their derivatives using \eqref{eq:euler_DtP}. As usual, the solution to the above weak form defines an (approximate) solution to the set of equations \eqref{eq:Euler} through the DtP mapping. This weak form corresponds to the first variation of the functional \eqref{eq:euler_dual_functional}.
\section{Results}\label{result}
The following sections describe examples for each of the problems developed in the previous sections. For the heat and transport equations, we solve the dual problems as space-time boundary value problems (bvp), while for Euler's system of ODE (with and without damping), we solve the dual problem as a two-point bvp in time. We use the Finite Element (FE) method to discretize all problems. In the case of heat/transport equations, we define the space-time domain $\Omega=\{(x,t):x \in (0,L), t \in (0,T)\}$, while for Euler's system of ODE, the time domain is defined as $\Omega=\{t:t \in (0,T)\}$. A linear span of globally continuous, piecewise smooth finite element shape functions corresponding to a FE mesh for $\Omega$ is used to achieve this discretization. These shape functions are represented by $N^{(\cdot)}$, where $(\cdot)$ denotes the index of the node under consideration. For the following examples, the \textit{error} for any field with respect to a reference field within the domain  is given by:
\[
\% \mbox{ error} = \left|\frac{\mbox{dual solution}-{\mbox{reference solution}}}{\mbox{reference solution}}\right| \times 100,
\]
 where the `dual solution' above refers to the primal field obtained by the DtP mapping from the corresponding dual approximation. In most cases, this measure is evaluated pointwise. The measure is not implementable at points where the reference field approaches zero. In such cases, alternative measures are used, as described in the relevant sections.

In this work, we use $C^0$ FE shape functions to approximate the dual fields. Hence, the direct evaluation of primal fields through the DtP mapping exhibit discontinuities in the domain in general, even when approximating continuous primal solutions. This is because the DtP mappings involve derivatives of the dual fields. To deal with this feature, our algorithm introduces the following device - essentially, a commonly used $L^2$ projection - to posit the corresponding primal solution, which we then compare with the exact primal solutions to test the accuracy of the overall scheme (including the projection). Let $u$ represent such a discontinuous primal approximation, determined from the dual solution. To achieve a $C^0$ continuous approximation in the domain $\Omega\subset\mathbb{R}^2 \mbox{ or } \mathbb{R}$ (depending on the problem), we employ the following method: Let $u_h$ represent the projection of $u$ onto a space $V_h$ formed by the linear span of globally continuous, piecewise smooth finite element shape functions corresponding to a FE mesh for $\Omega$. We enforce the following condition upon $u_h$:
\begin{equation*}
     u_h = \operatorname*{arg\,min}_{v\in V_h}  \int_\Omega \frac{1}{2} |u-v|^2 \, d\Omega.\label{eq:weak_L2}
\end{equation*} 
Let $N^A$ represent the basis functions associated with any node with an index $A$. The discrete version of the optimality condition of the above statement is:
\begin{equation}\label{eq:L2_project}
    \sum_{A=1}^N\sum_{B=1}^N\delta u^A_h\biggl(\int_\Omega N^A\, N^B \,d\Omega \biggl) u_h^B=  \sum_{A=1}^N\delta u^A_h \int_\Omega N^A\, u\,d\Omega
\end{equation}
where $u_h^B$ denotes the nodal value at node $B$ of the sought continuous projection of $u$, and  $\delta u_h := \delta u^A_h N^A$ is a test function. In solving \eqref{eq:L2_project} we impose any known data, e.g.~initial and boundary conditions of the primal problem, as known function values $u^A_h$, at corresponding nodes $A$ where the data is known, with $\delta u^A_h = 0$. Using the arbitrariness of the remaining `free' $\delta u^A_h$, we obtain a system of linear equations to solve for the unconstrained nodal values $u_h^B, B \in \{1, \ldots, N\}$.

For several examples presented in the following sections, partitioning the domain into smaller domains (stages) along the time direction simplifies computation and aids in refinement. Each stage, indexed by $s$, uses initial conditions based on the previous stage's results, solving the full problem up to the required final time by concatenating primal fields obtained by DtP mappings for successive stages along increasing time.
We will refer to this process as `time-slicing'. The details of this method are explained in the relevant sections. As also explained later, higher errors in the primal solution occur near the final time of a stage since high gradients in the dual solution have to be resolved, in general, due to the imposition of an arbitrarily fixed final-time boundary condition (this phenomenon also occurs at the outflow boundary of the linear transport equation for the same reason). However,  a portion of the results near the final time of any stage is, and can be, discarded without loss of generality. Thus, the domain is extended beyond the required final time of any specific problem and the obtained results are retained only up to the specified final time.

In the following, we use summation over repeated vector or tensor indices, as well as FE formulation related summations, with summation performed over the appropriate index ranges, unless explicitly mentioned otherwise.

\subsection{Heat equation}
In order to solve the dual heat equation, the weak form \eqref{eq:weak_heat} is employed. Let $N^A$ represent the shape function corresponding to any node $A$. Using the discrete approximations 
\begin{equation*}
    \begin{gathered}
      p(x,t) = d_1^A N^A(x,t);  \qquad  l(x,t)=d_2^AN^A(x,t), \\  
    \end{gathered}
\end{equation*}
where $N$ is the total number of nodes in the FE mesh,
we obtain the following set of equations:
\begin{equation*}\label{eq:heat_discrete}
        K_{ij}^{AB} d_j^B = R_i^A;
\end{equation*}
where
\begin{equation*}
    \begin{gathered}
        \qquad K_{11}^{AB} = 
            \int_0^T dt \int_0^L dx \Bigl(\,- \p_x N^A(x,t) \,\p_xN^B(x,t) - \,  N^A(x,t)\, N^B(x,t)\Bigl);\\
                    \qquad K_{12}^{AB} = \int_0^T dt \int_0^L dx \Bigl(\,- \p_x N^A(x,t) \,\p_tN^B(x,t) + \, k N^A(x,t) \,\p_x N^B(x,t)\Bigl);  \\
                   \qquad K_{21}^{AB} =  \int_0^T dt \int_0^L dx \Bigl(\,- \p_t N^A(x,t) \,\p_xN^B(x,t) + \, k \, \p_x N^A(x,t)\, N^B(x,t)\Bigl);  \\
                   \qquad K_{22}^{AB} =  \int_0^T dt \int_0^L dx \Bigl(\,- \p_t N^A(x,t) \,\p_tN^B(x,t) - \, k^2 \p_x N^A(x,t) \,\p_x N^B(x,t)\Bigl) ;  \\
R_1^A = \int_0^T dt \, N^A(0,t) \, \theta_l(t); \qquad 
R_2^A = \int_0^T dt \, N^A(L,t) \, k\pi_r(t) 
+ \int_0^L dx \, N^A(x,0) \, \theta_0(x).
    \end{gathered}
\end{equation*}
Additionally, we apply the Dirichlet Boundary conditions as
\begin{equation*}
\begin{gathered}
l(0,t) = l_l(t);\qquad p(L,t) = p_r(t); \qquad l(x,T) = l_T(x).
\end{gathered}
\end{equation*}
If the boundary condition at $x=1$ is altered from $\pi(L,t)=\pi_r(t)$ to $\theta(L,t)=\theta_r(t)$ in \eqref{eq:primal_bc}, $R_i^A$ gets modified in the obvious manner following the left boundary condition:
\begin{equation*}
R_1^A = \int_0^T dt \, \bigl(N^A(0,t) \, \theta_l(t) 
-  \, N^A(L,t) \, \theta_r(t)\bigl); \qquad 
R_2^A = 
\int_0^L dx \, N^A(x,0) \, \theta_0(x),
\end{equation*} 
and the Dirichlet boundary condition $p(L,t)=p_r(t)$ changes to  $l(L,t)=l_r(t)$ on the right boundary. Each of the integrals in the above expressions are evaluated using a two-point Gauss quadrature scheme in each direction. Once the discretized dual fields $p(x,t)$ and $l(x,t)$ are obtained, the primal field $\theta(x,t)$ is evaluated at the Gauss points and the scheme explained in the preamble of Sec.~\ref{result} is used to generate its continuous projection.
\subsubsection{Steady solution of the Heat equation}\label{sec:steady_heat}
\begin{figure}[ht]
\captionsetup[sub]{font=small,justification=centering}
\centering
\begin{subfigure}{.32\textwidth}
  \centering
  % include first image
  \includegraphics[width=0.9\linewidth]{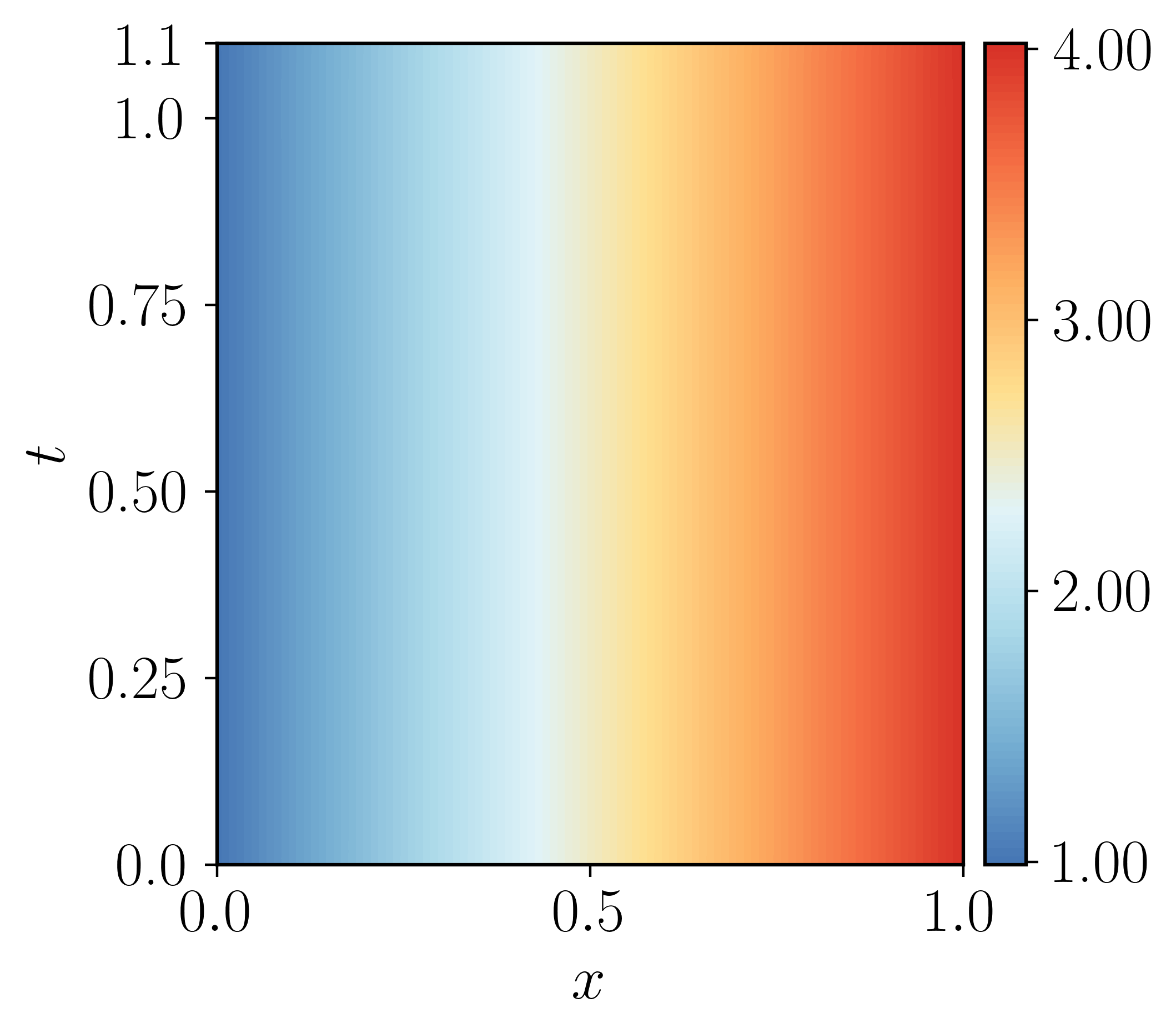}  
  \caption{$\theta(x,t)$}
  \label{fig:sub-heat_p1_theta_1.1}
\end{subfigure}
\begin{subfigure}{.32\textwidth}
  \centering
  % include second image
  \includegraphics[width=0.9\linewidth]{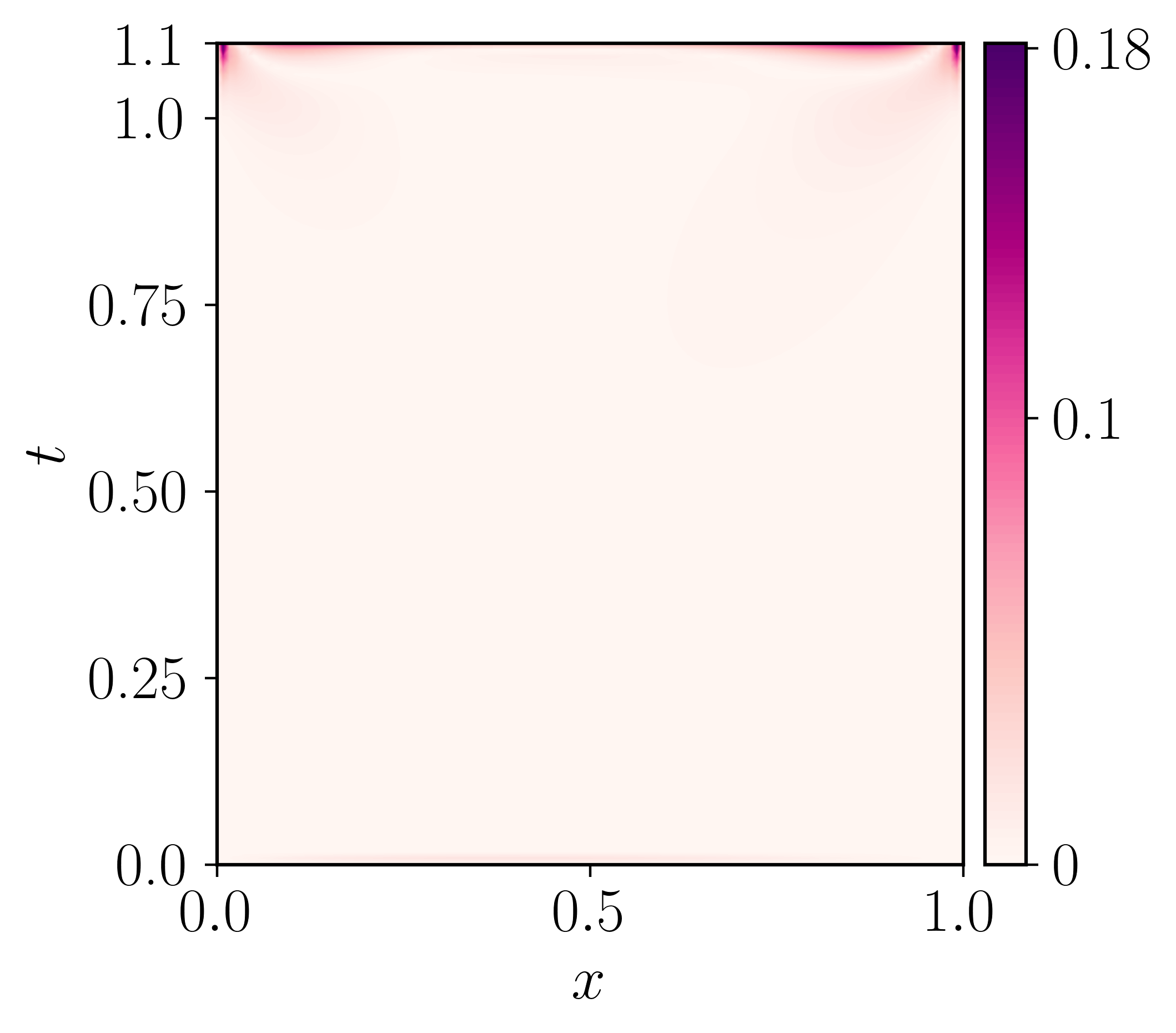}  
  \caption{ $\%$ error in $\theta(x,t)$}
  \label{fig:sub-heat_p1_error_1.1}
\end{subfigure}
\begin{subfigure}{.32\textwidth}
  \centering
  % include first image
  \includegraphics[width=0.9\linewidth]{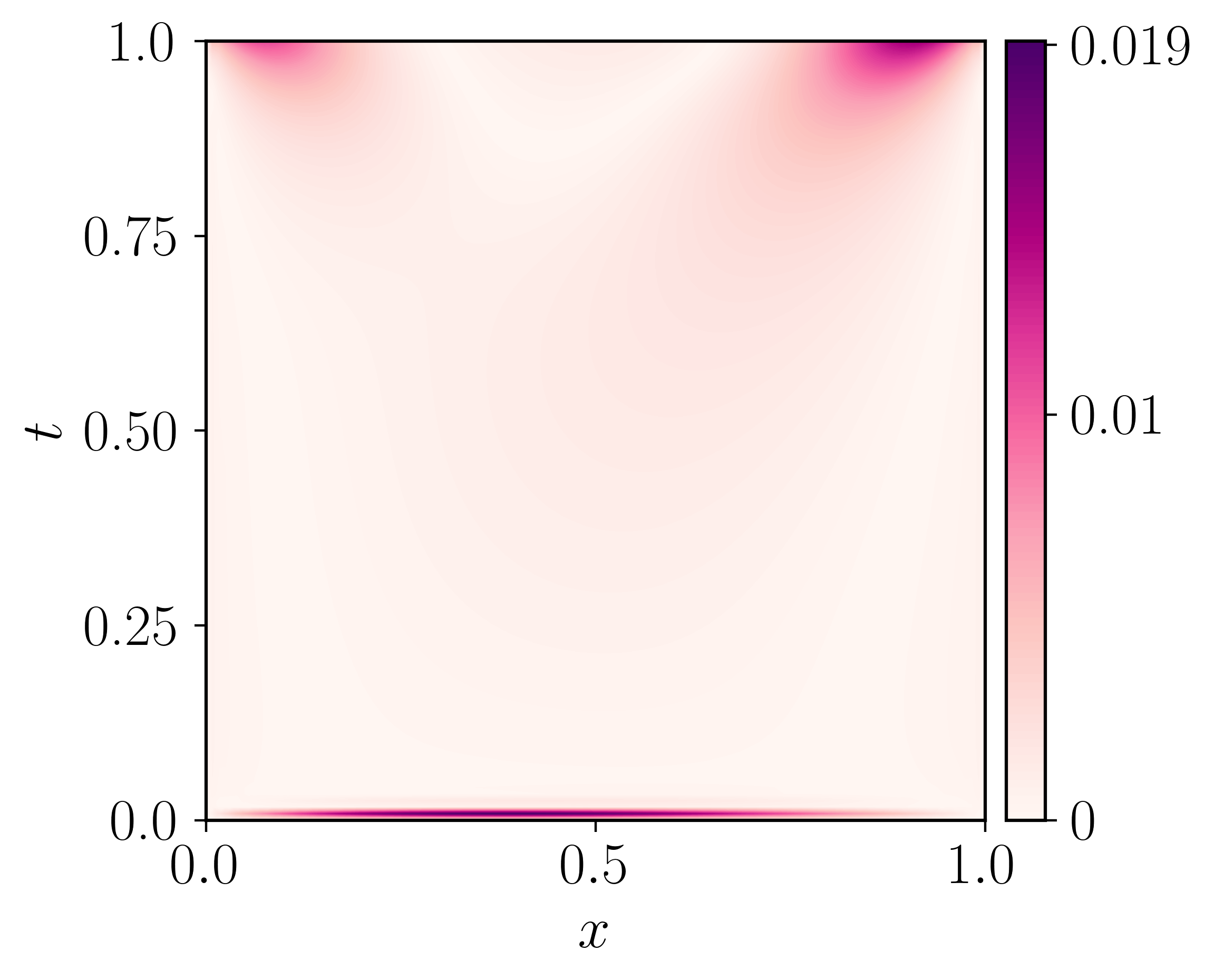}  
  \caption{$\%$ error in $\theta(x,t)$ \\  up to $T=1$}
  \label{fig:sub-heat_p1_error}
\end{subfigure}
\caption{DtP mapping generated primal field $\theta$ for the steady state heat equation. The mesh is $100\times110$, $T=1.1$ and $k=1$.}
\label{fig:heat_p1_theta_error}
\end{figure}

A steady solution of the heat equation \eqref{eq:primal_ht} for $k = 1$ given by 
\begin{equation*}\label{eq:heat_p1_exact_sol}
  \theta(x,t) = 3x+1 %%\qquad \forall t\in[0,\infty)  
\end{equation*}
is considered on a domain $\Omega$ with $L=1$, satisfying the following initial and boundary conditions:
\begin{equation}\label{eq:heat_p1_primal_bcs}
    \theta_0(x) = 3x+1; \qquad \theta_l(t) = 1; \qquad \theta_r(t) = 4.
\end{equation}
In order to evaluate our dual scheme, we approximate the solution to the current problem by applying the boundary conditions 
\begin{equation}\label{eq:dual_heat_hom_l_p}
l_l(t)=l_r(t)=0; \qquad l_T(x)=0
\end{equation}
in the weak form \eqref{eq:weak_heat}. The Dirichlet boundary conditions at $x=0$ and $x=L$ in the primal problem are transformed into oblique \textit{natural} boundary conditions in the dual formulation. For a final time of $T=1.1$ and a mesh of $100\times110$ elements in the space-time domain, the field $\theta(x,t)$ and the error obtained with respect to the exact solution are shown in Fig.~\ref{fig:sub-heat_p1_theta_1.1} and Fig.~\ref{fig:sub-heat_p1_error_1.1},  respectively.

As evident from \ref{fig:sub-heat_p1_error_1.1}, while the errors are minimal, they are higher at at $T=1.1$. Given the arbitrariness in the choice of imposition of b.c.s for the dual problem at the final time and the latter's value, given a final time of interest for the primal problem, we can always choose $T$ appropriately and discard a reasonable portion of data from the dual solution near the final time and focus on results within the remaining domain. For illustration, Fig.~\ref{fig:sub-heat_p1_error} displays the results until $T=1$, which has a significantly lower maximum error compared to until $T=1.1$.

    \begin{figure}[ht]
    \captionsetup[subfigure]{justification=centering}
        \centering
        \begin{subfigure}[b]{0.475\textwidth}
            \centering
            \includegraphics[width=0.9\textwidth]{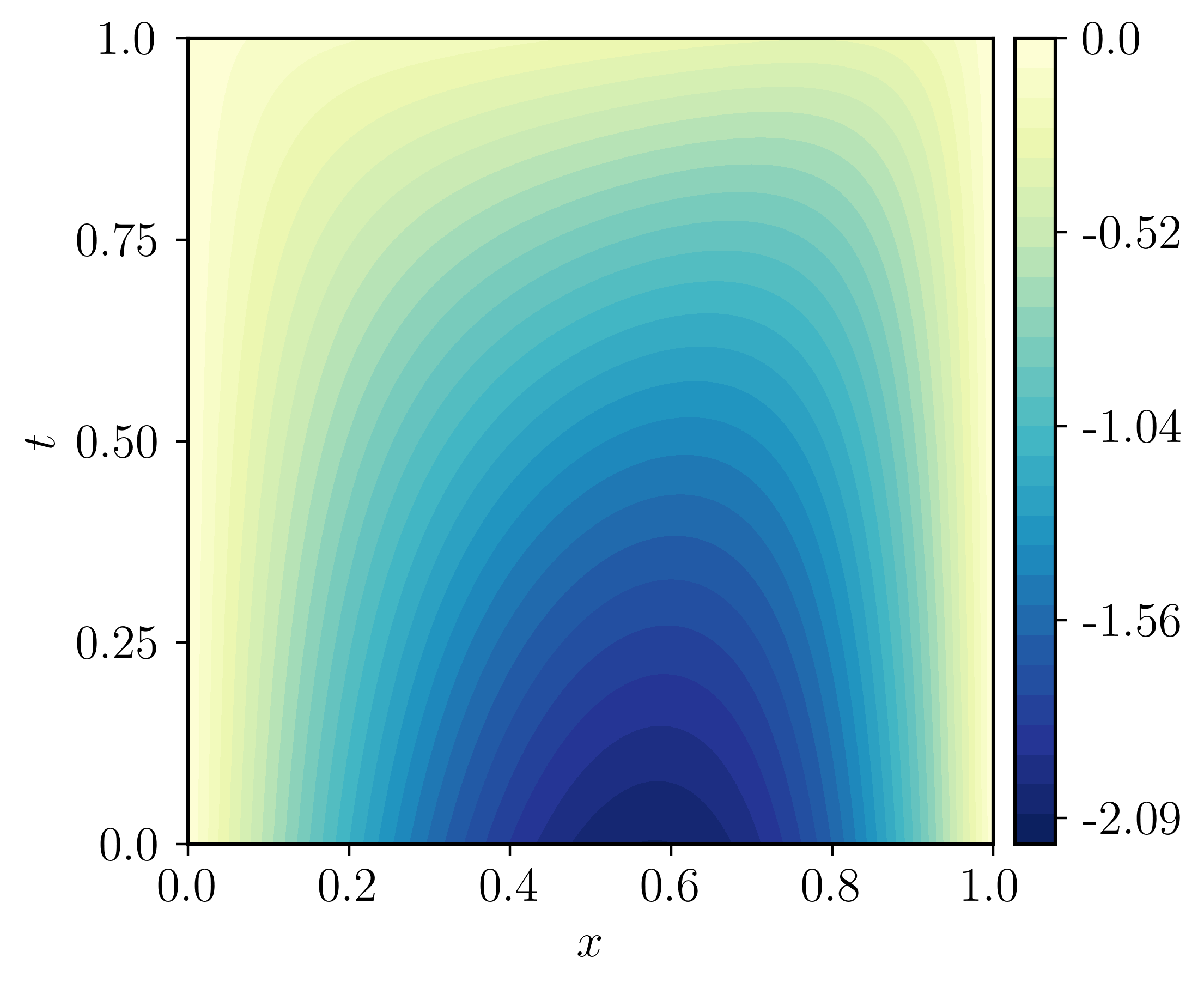}
            \caption[l]%
            {Transient Dual field $l(x,t)$}    
            \label{fig:sub-heat_p1_l1}
        \end{subfigure}
        \hfill
        \begin{subfigure}[b]{0.475\textwidth}  
            \centering 
            \includegraphics[width=0.9\textwidth]{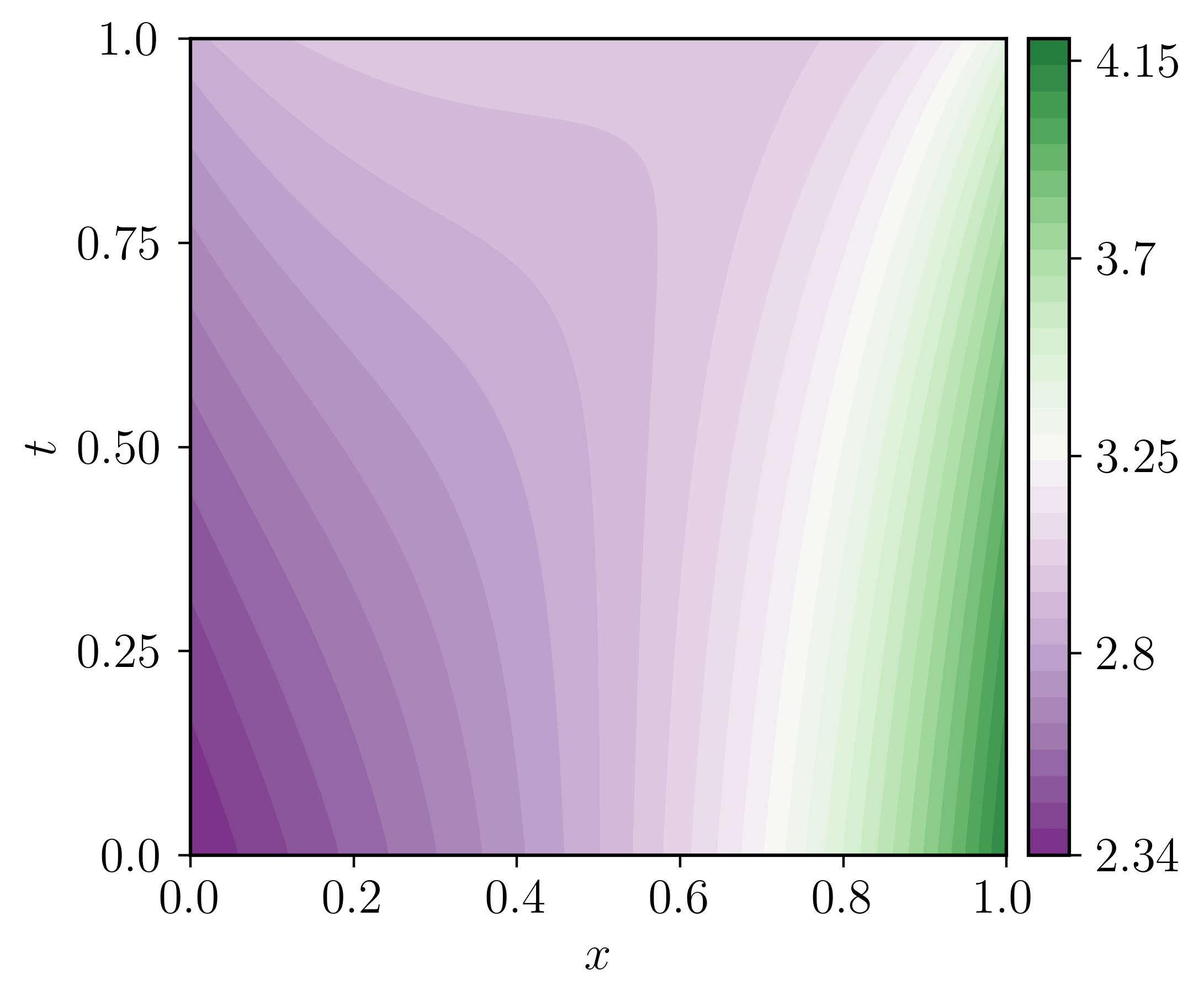}
            \caption[]%
            {Transient Dual field $p(x,t)$}    
            \label{fig:sub-heat_p1_p1}
        \end{subfigure}
        \vskip\baselineskip
        \begin{subfigure}[b]{0.475\textwidth}   
            \centering 
            \includegraphics[width=0.9\textwidth]{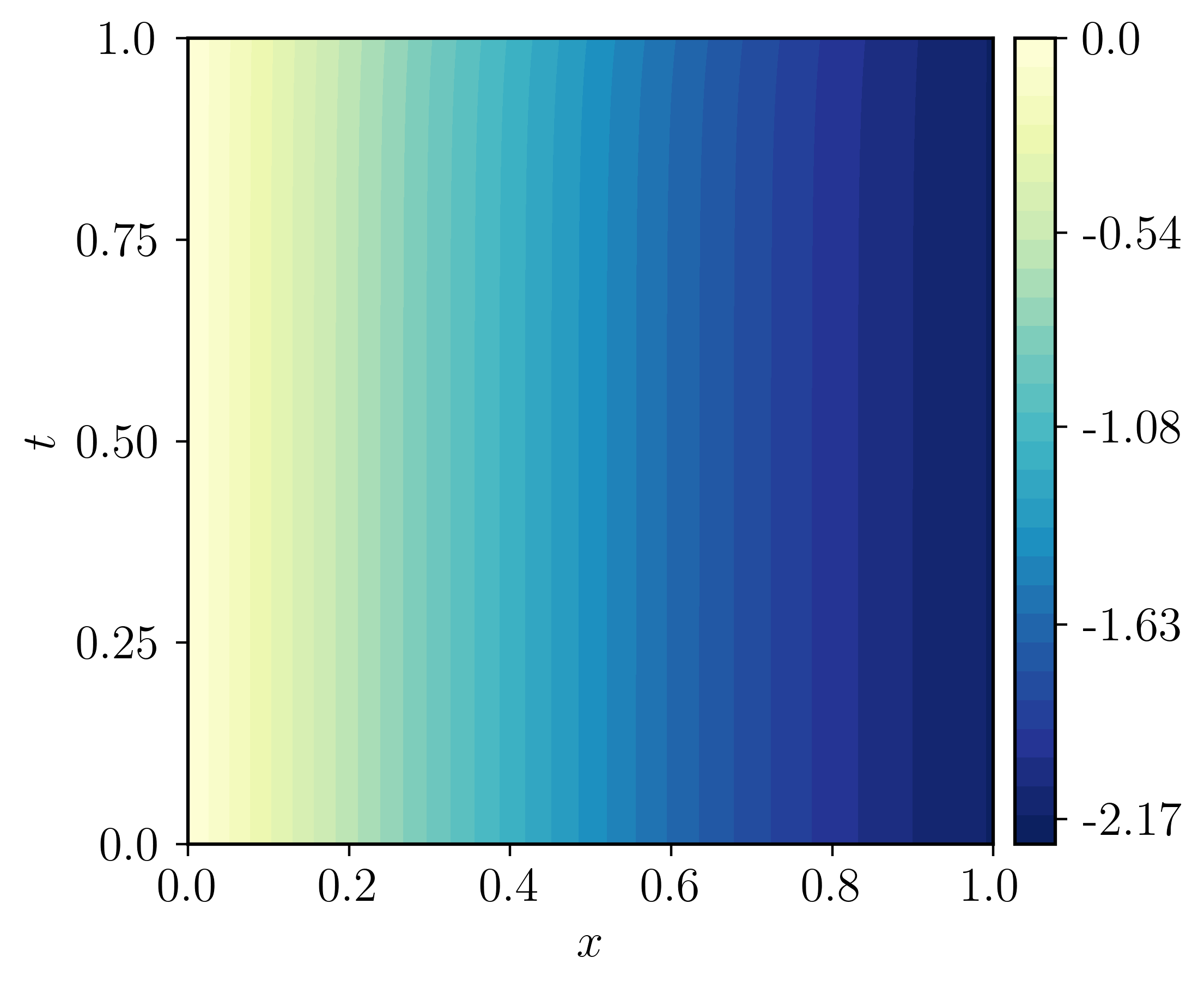}
            \caption[]%
            {Steady Dual field $l(x,t)$}   
            \label{fig:sub-heat_p1_l2}
        \end{subfigure}
        \hfill
        \begin{subfigure}[b]{0.475\textwidth}   
            \centering 
            \includegraphics[width=0.9\textwidth]{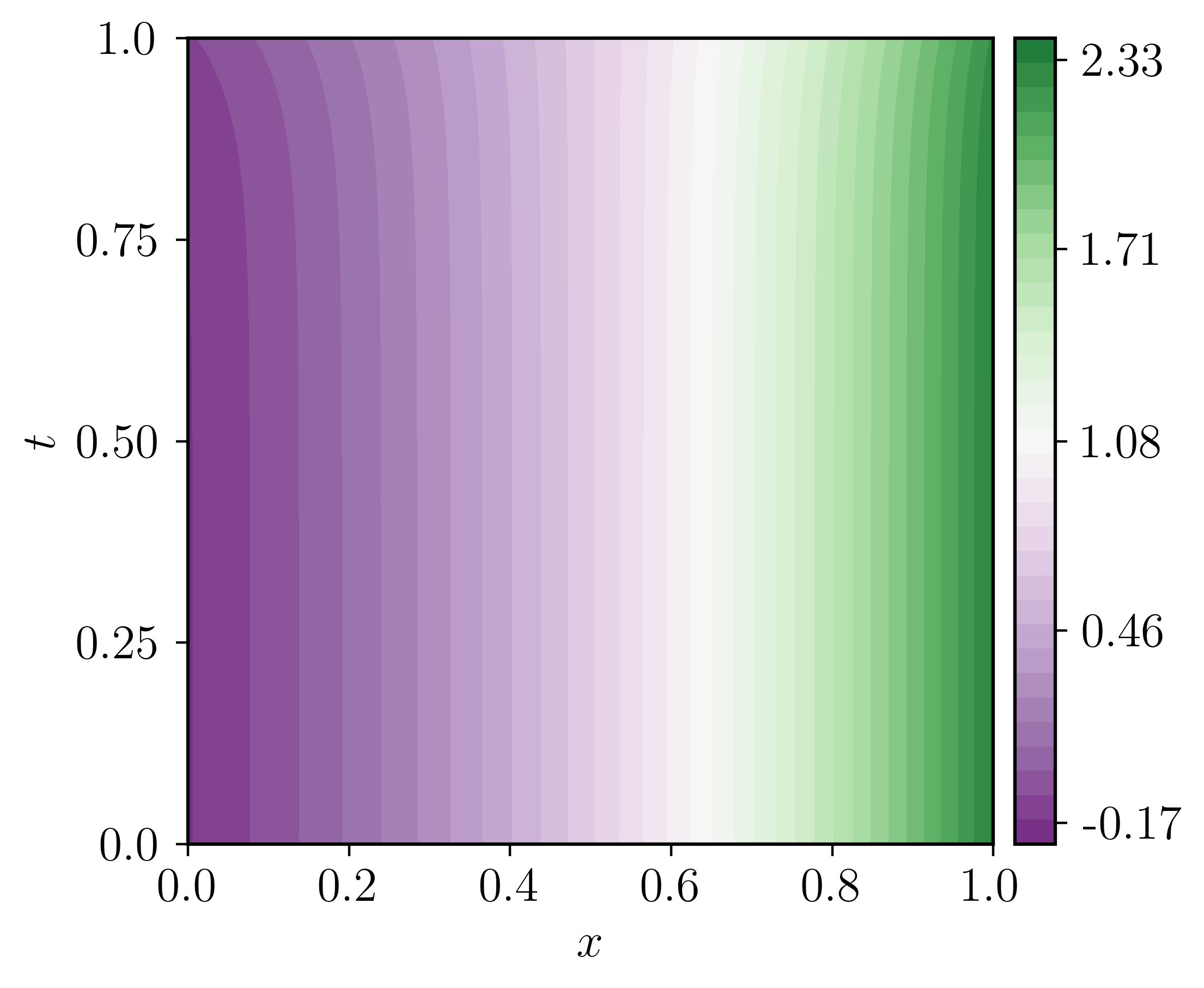}
            \caption[]%
            {Steady Dual field $p(x,t)$}    
            \label{fig:sub-heat_p1_p2}
        \end{subfigure}
        \caption{Dual fields for steady state heat equation. (a) and (b) represent dual fields when $l_l(t)=l_r(t)=l_T(x)=0$ are imposed. (c) and (d) represent the dual fields when the boundary condition set \eqref{eq:heat_p1_steady_bc} is imposed.}
        \label{fig:heat_p1_duals}
    \end{figure}

The dual fields for this problem are presented in Fig.~\ref{fig:sub-heat_p1_l1}-\ref{fig:sub-heat_p1_p1}. Although these fields are transient in nature, the corresponding primal solution is steady in time.
As explained in \cite{action_3}, when the primal problem has a unique solution, \textit{any} solution to the dual problem, regardless of the choice of $H$ and the space-time b.cs imposed on the dual problem beyond the imposition of the primal set, must recover the unique primal solution through the DtP mapping.  As an illustration of this invariance of the DtP mapping generated primal solution from two different dual solutions, we set the following Dirichlet Boundary condition on the dual fields:
\begin{equation}\label{eq:heat_p1_steady_bc}
l_l(t)=0;\qquad l_r(t)=-\frac{13}{6};\qquad l_T(x)= \frac{x^3}{3} + \frac{x^2}{2}-3x.
\end{equation}
The computed dual solution with these b.c.s are shown in Fig.~\ref{fig:sub-heat_p1_l2}-\ref{fig:sub-heat_p1_p2}. It can be easily verified that
\begin{equation}
\label{eq:heat_p1_steady_duals}
p(x,t) = \frac{3x^2}{2}+2; \qquad l(x,t) =\frac{1}{k}\left(\frac{x^3}{3} + \frac{x^2}{2}-3x\right)
\end{equation}
satisfies the boundary conditions \eqref{eq:heat_p1_steady_bc} and the dual heat equation \eqref{eq:dual_heat_pde}. Furthermore, on application of the DtP mapping \eqref{eq:DtP}, it results in the expression $\theta(x,t) = 3x+1$ which is the unique solution to the primal problem, as must be. We note that dual solution \eqref{eq:heat_p1_steady_duals}  corresponding to b.cs \eqref{eq:heat_p1_steady_bc} is steady, unlike the dual solution corresponding to the b.cs \eqref{eq:dual_heat_hom_l_p}. This fact is also evident in the simulation results Fig.~\ref{fig:sub-heat_p1_l2}-\ref{fig:sub-heat_p1_p2} (except close to $T = 1$, without any essential loss of generality, as explained earlier).
\subsubsection{Transient solution of the Heat equation by duality}
The purpose of the following example is to demonstrate an approximation of a time-dependent solution of the heat equation for the following combination of Dirichlet and Neumann boundary conditions:
\begin{equation*}
    \theta_l(t) = 1; \qquad \pi_r(t) = 0,  
    \end{equation*}
and the initial condition 
    $$\theta_0(x) = \sin{\frac{\pi x}{2}} + 1.$$

\begin{figure}[h]
\centering
\begin{subfigure}{.32\textwidth}
  \centering
  % include the first image
  \includegraphics[width=.9\linewidth]{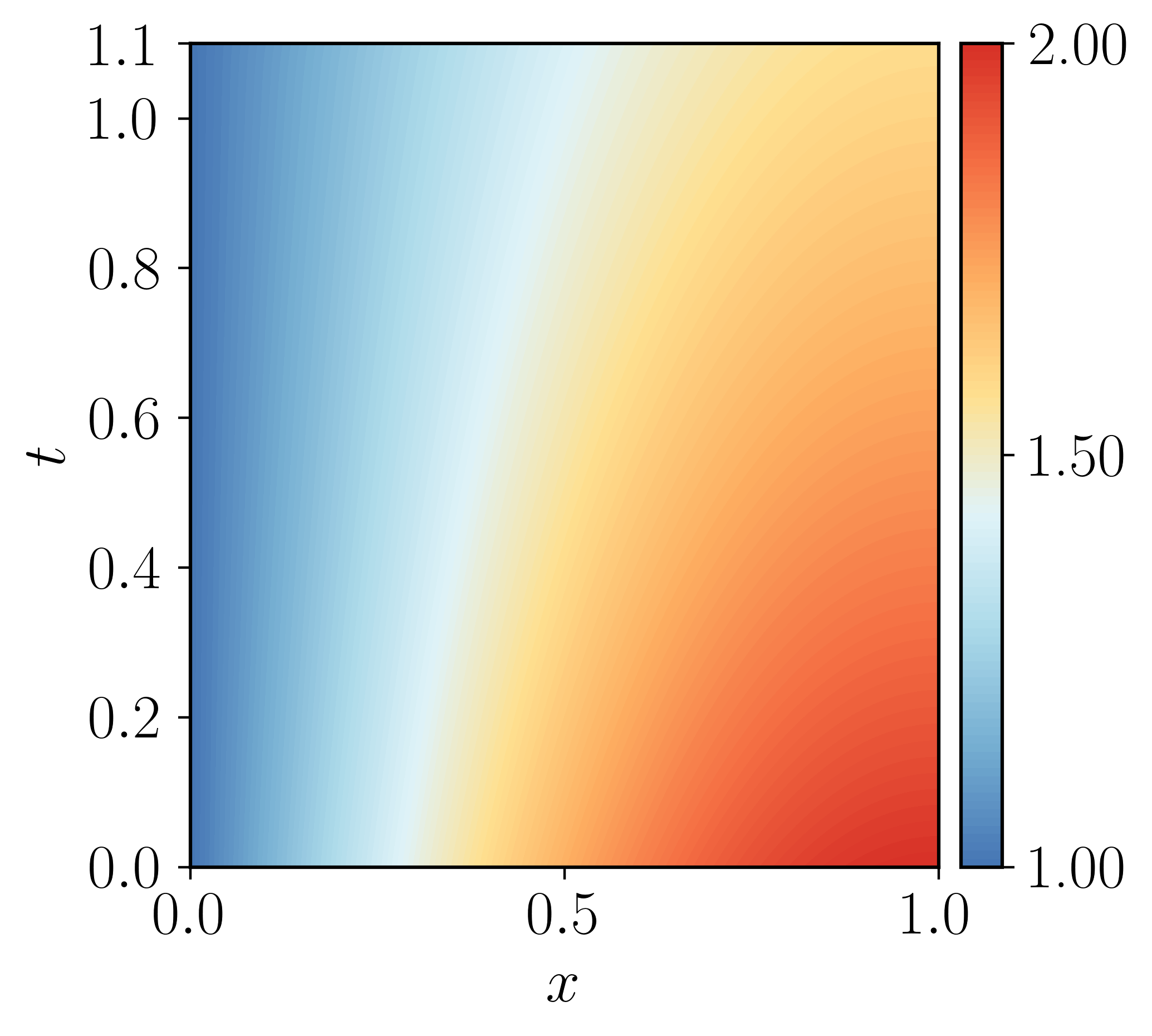}  
  \caption{$\theta(x,t)$}
  \label{fig:sub-heat_p2_theta_1.1}
\end{subfigure}
\begin{subfigure}{.32\textwidth}
  \centering
  % include second image
  \includegraphics[width=.9\linewidth]{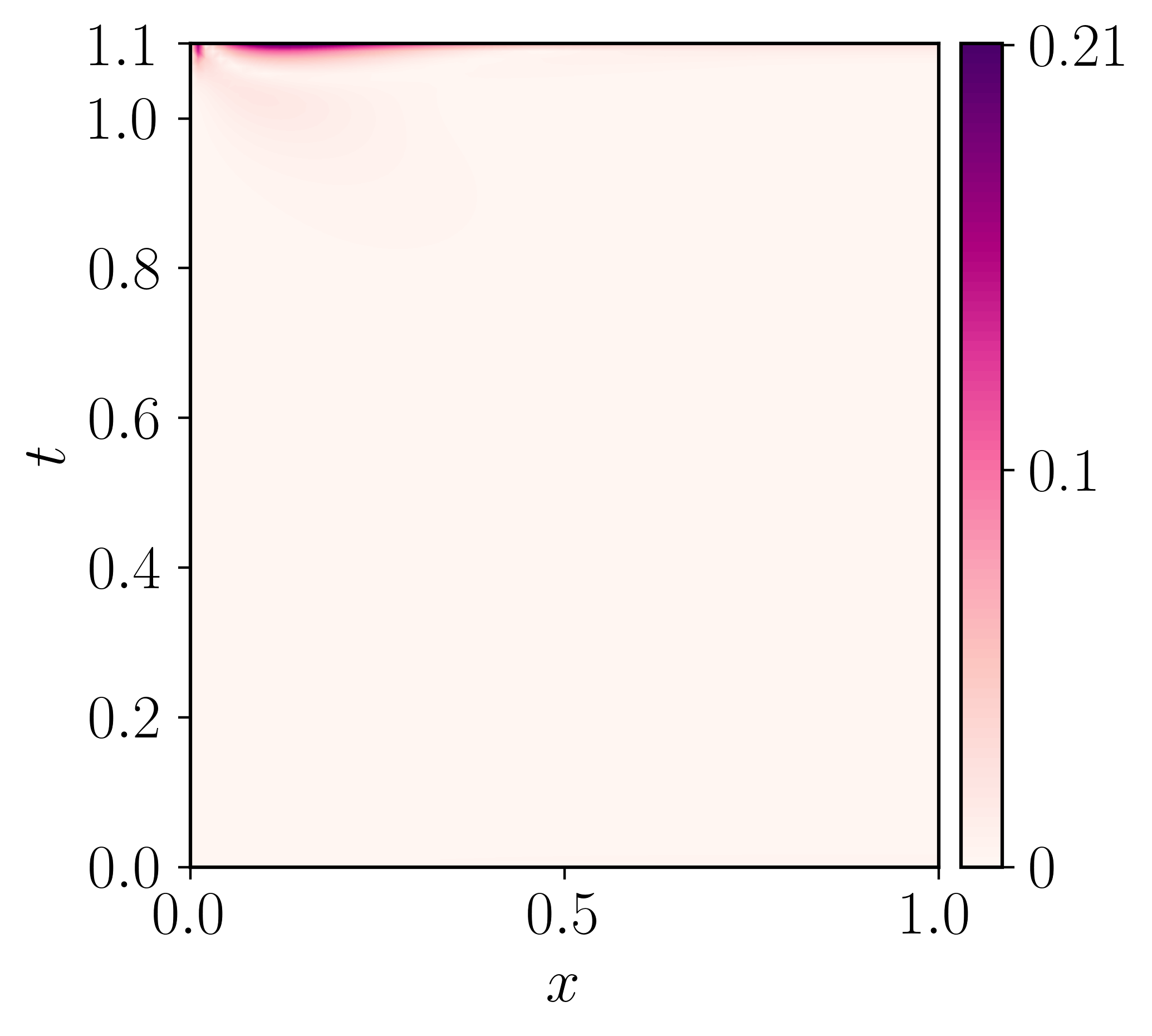}  
  \caption{$\%$ error in $\theta(x,t)$}
  \label{fig:sub-heat_p2_error_1.1}
\end{subfigure}
\begin{subfigure}{.32\textwidth}
  \centering
  % include the first image
  \includegraphics[width=.9\linewidth]{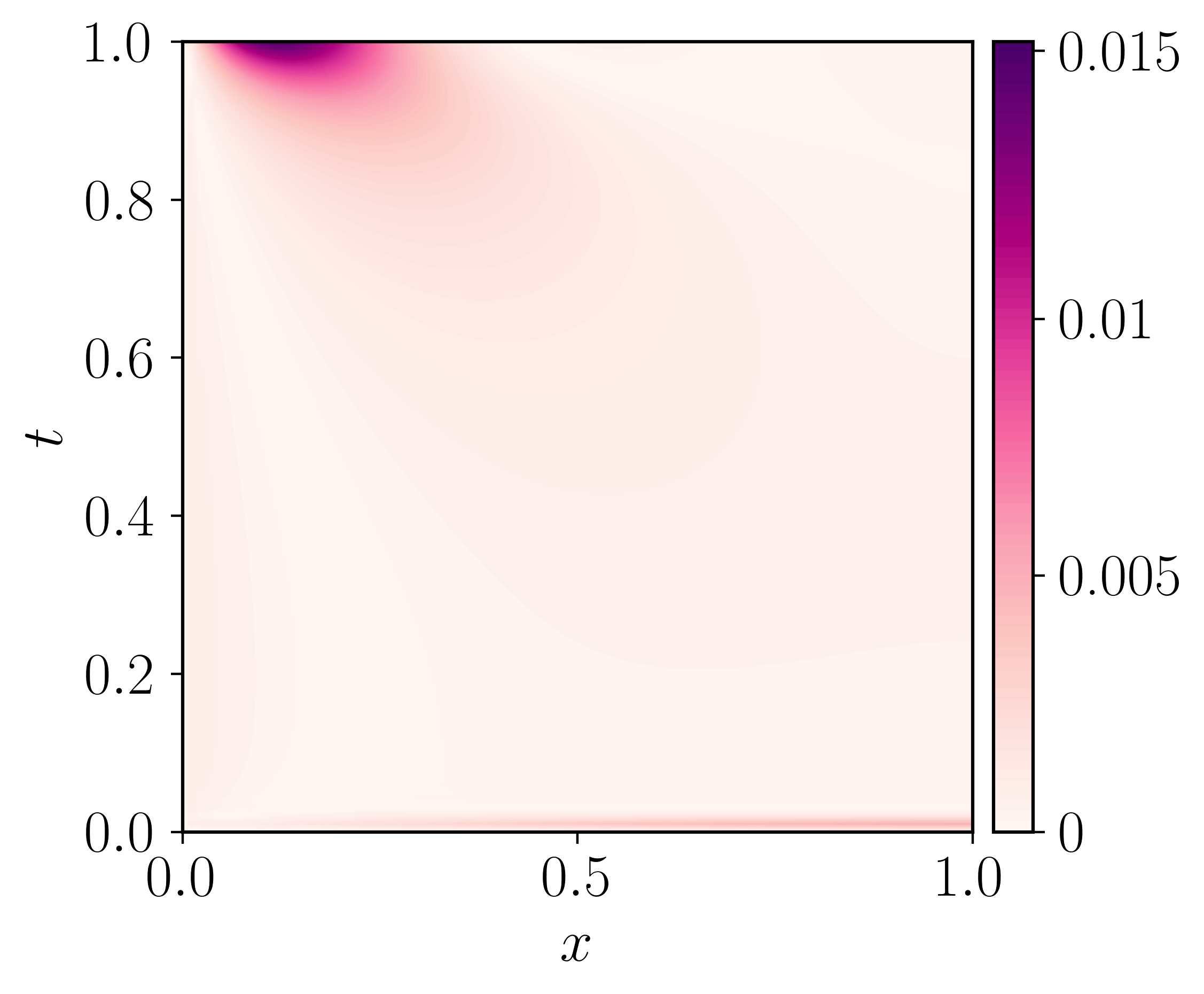}  
  \caption{$\%$ error in $\theta(x,t)$ \\ upto $T=1$}
  \label{fig:sub-heat_p2_error}
\end{subfigure}
\caption{DtP mapping generated primal field $\theta$ for the transient heat equation. Mesh is $100\times110$, $T=1.1$ and $k=0.2$.}
\label{fig:heat_p2_theta_error}
\end{figure}

For $L=1$, the exact solution in this case is given by:
$$\theta(x,t) = \sin\left(\frac{\pi x}{2}\right) e^{-\frac{\pi^2 kt}{4}}+1. $$
Applying $l_T(x)=p_r(t)=l_l(t) =0$ as the boundary conditions along with $k=0.2$, we utilize \eqref{eq:weak_heat} and simulate the problem up to $T=1.1$. The results for a mesh of $100\times110$ are presented in Fig.~\ref{fig:sub-heat_p2_theta_1.1} and Fig.~\ref{fig:sub-heat_p2_error_1.1}. The error plot up to a time $T=1$ is also shown in Fig.~\ref{fig:sub-heat_p2_error} which illustrates that the error in the interior of the domain is smaller relative to that near the top boundary. 
\subsubsection{Discontinuous initial data for the Heat equation}
\begin{figure}[!ht]
\centering
\begin{subfigure}{.32\textwidth}
  \centering
  % include first image
  \includegraphics[width=0.95\linewidth]{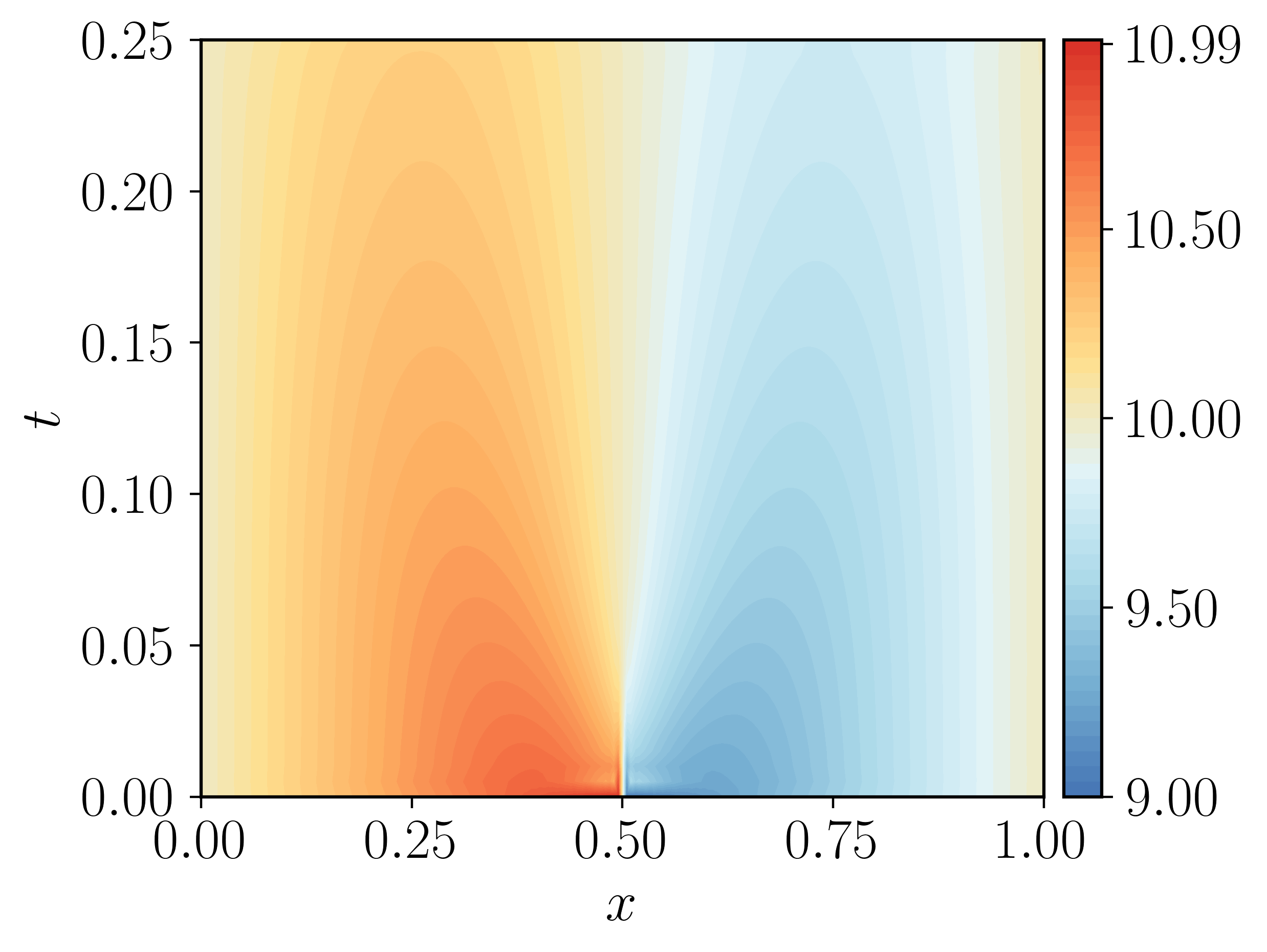}  
  \caption{$\theta(x,t)$}
  \label{fig:heat_p3_theta}
\end{subfigure}
\begin{subfigure}{.32\textwidth}
  \centering
  % include first image
  \includegraphics[width=0.95\linewidth]{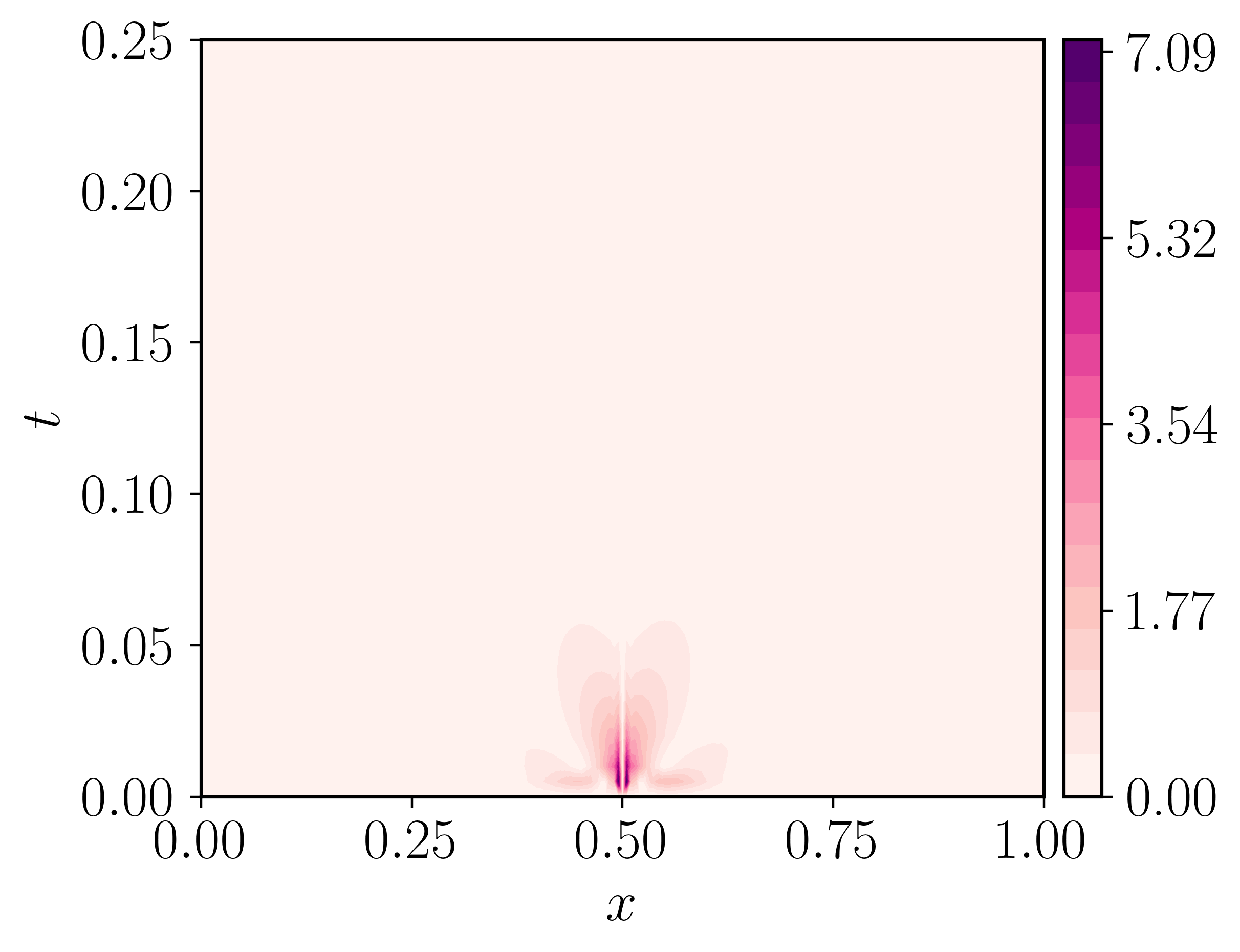}  
  \caption{$\%$ Error for $\theta(x,t)$}
  \label{fig:sub-heat_p3_error_ana}
\end{subfigure}
\begin{subfigure}{.32\textwidth}
  \centering
  % include second image
  \includegraphics[width=0.95\linewidth]{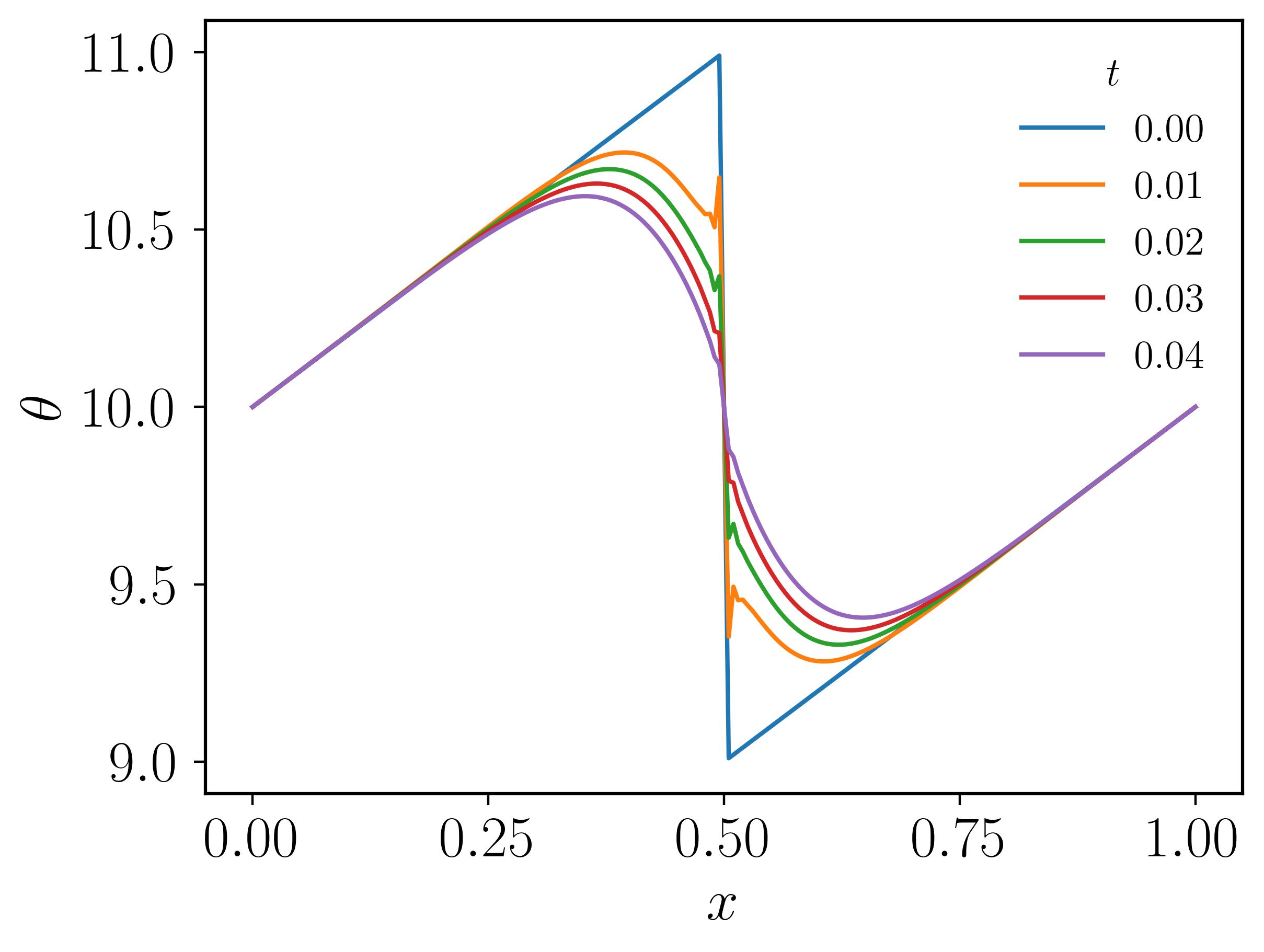}  
  \caption{Line plots for $\theta(x,t)$}
  \label{fig:heat_p3_line_plot}
\end{subfigure}
\caption{DtP mapping generated primal field $\theta$ for a discontinuous initial condition. Mesh is $200\times50$, $T=0.25$ and $k=0.1$. The line plots for $\theta(x,t)$ close to the initial time are illustrated in (c).}
\label{fig:heat_p3_disconitnuous_ic}
\end{figure}

As a final example involving the heat equation, we investigate how the dual methodology deals with any problem which involves a discontinuity in initial conditions, especially since primal Dirichlet b.cs and initial conditions get imposed as natural boundary conditions in the dual problem. Consider the following set of initial and boundary conditions for $L=1$:
\begin{equation}\label{eq:heat_p3_initial_condition}
\theta_0(x) = \begin{cases}
10+2x &  \mbox{for } 0 \leq x< 0.5 \\ 
8+2x &  \mbox{for }  0.5 < x \leq 1
\end{cases}; 
\qquad \theta_l(t) = 10;
\qquad \theta_r(t) = 10.
\end{equation}
In this case, separation of variables along with the application of Fourier series yields the solution
\begin{equation}\label{heat_p3_exact}
    \theta(x,t) = 10 + 2\sum_{m=1}^\infty \frac{(-1)^{m+1}}{\pi m} \sin{(2\pi mx)}e^{-(2m)^2\pi^2kt}\qquad \text{ for }0\leq x\leq1.
\end{equation}
It has been verified that the series converges in the $L^2$ sense to a value within machine precision on summing the first $10^5$ terms and beyond, and we retain $10^5$ terms for generating the exact solution. Additionally, at the point of discontinuity, the function value approximates to  the average of the values obtained from the left and the right of the jump. Consequently, for the initial condition in the simulation, we provide a value of $10$, the average jump, at $x=0.5$. The results are presented in Fig.~\ref{fig:heat_p3_disconitnuous_ic}.
\begin{figure}[!ht]
\centering
\end{figure}

\begin{figure}[ht]
    \captionsetup[subfigure]{justification=centering}
        \centering
        \begin{subfigure}[b]{0.475\textwidth}
            \centering
            \includegraphics[width=0.9\textwidth]{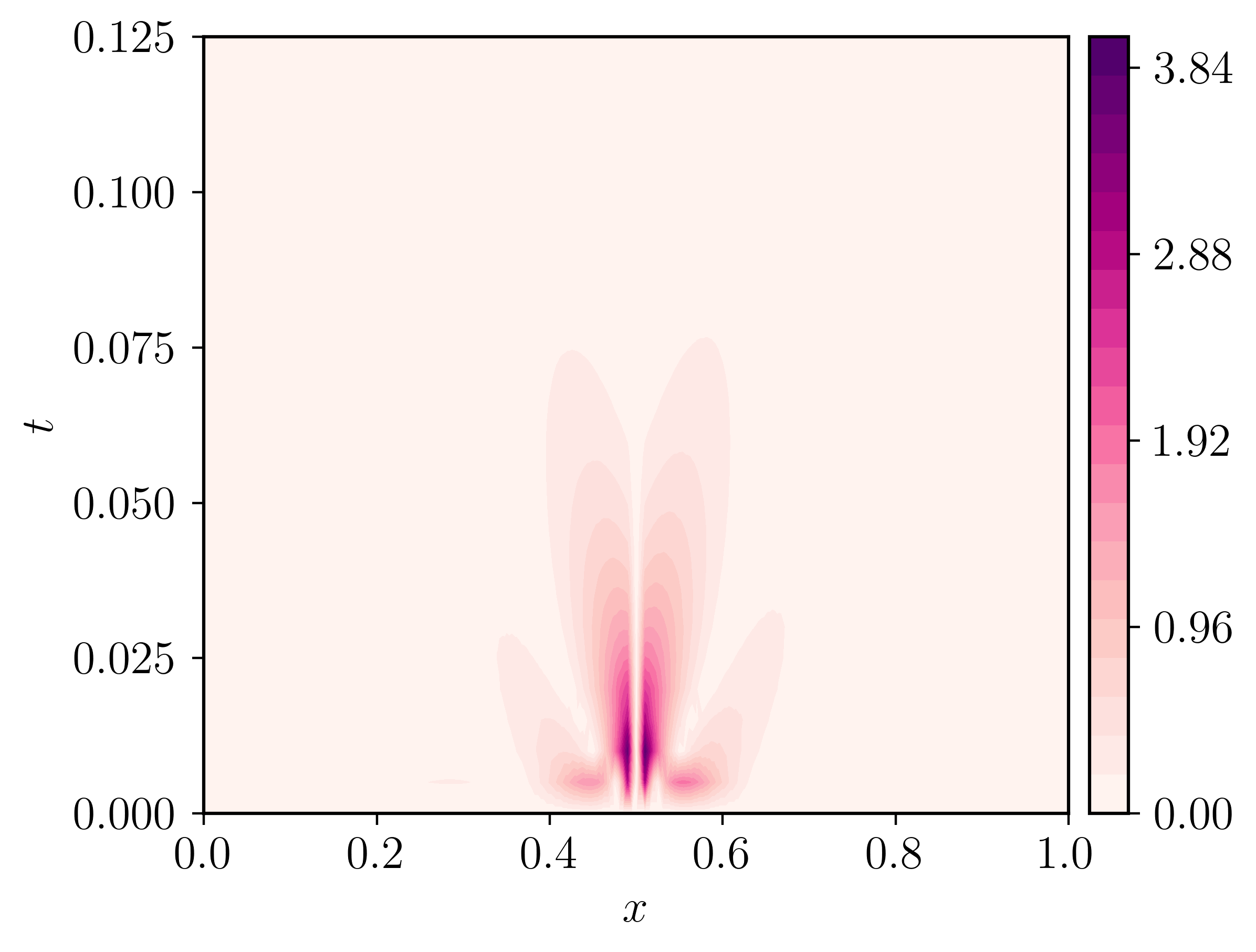}
            \caption[l]%
            {$\%$ Error for $\theta(x,t)$}    
            \label{fig:sub-heat_p3_eps200_error}
        \end{subfigure}
        \hfill
        \begin{subfigure}[b]{0.475\textwidth}  
            \centering 
            \includegraphics[width=0.9\textwidth]{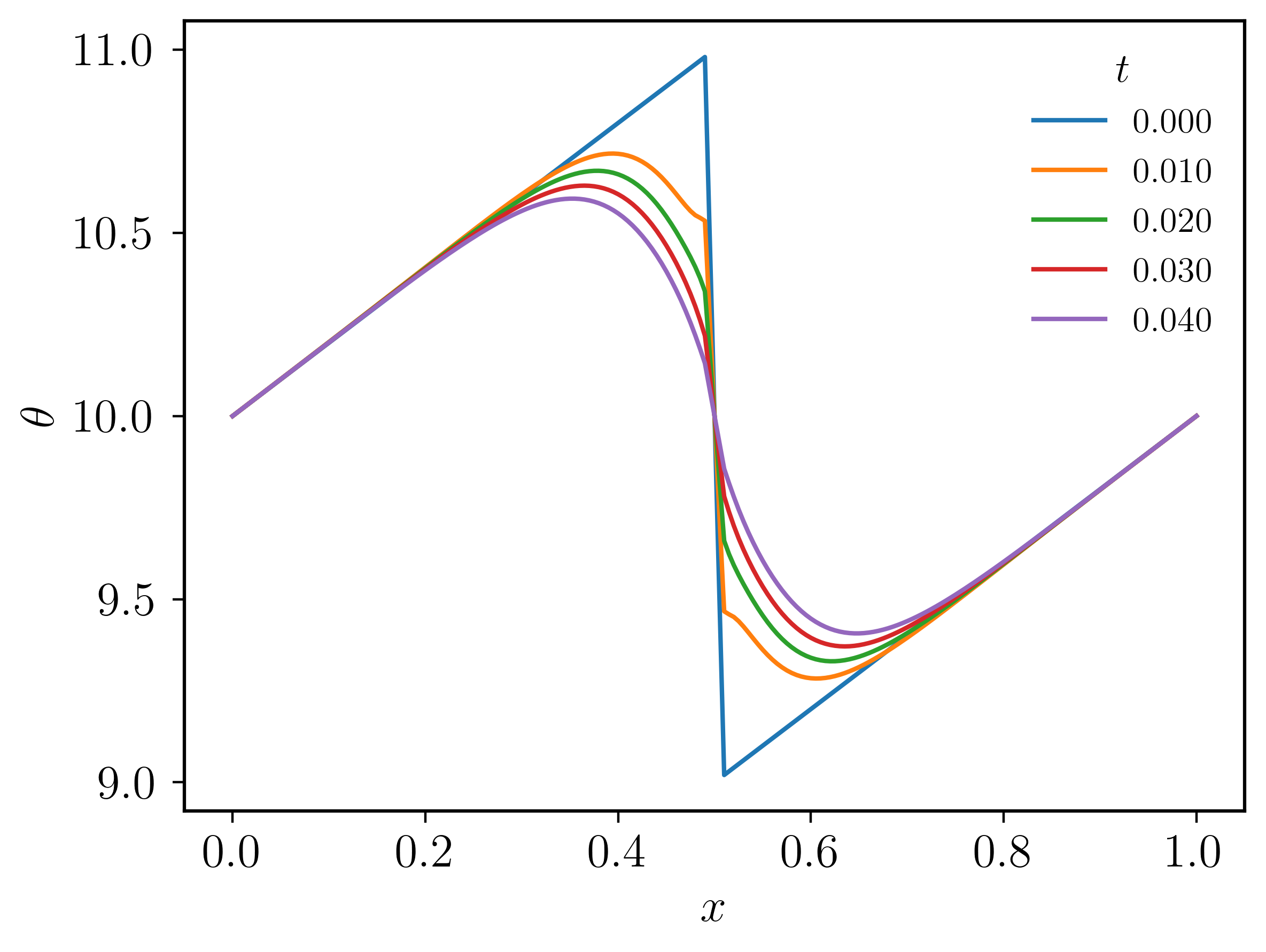}
            \caption[]%
            {Line plot for $\theta(x,t)$}    
            \label{fig:sub-heat_p3_eps200_line_plot}
        \end{subfigure}
        \vskip\baselineskip
        \begin{subfigure}[b]{0.475\textwidth}   
            \centering 
            \includegraphics[width=0.9\textwidth]{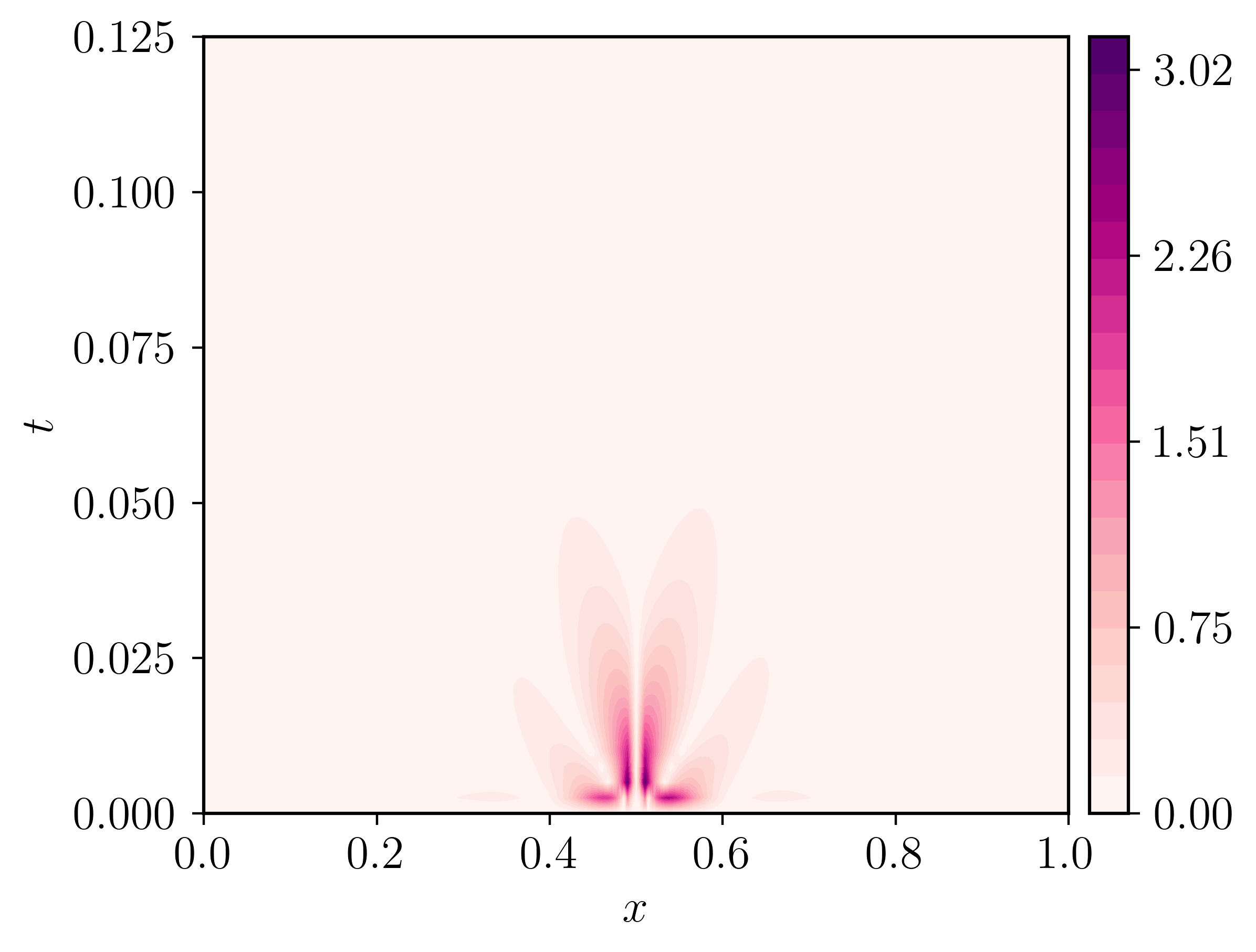}
            \caption[]%
            {$\%$ Error for $\theta(x,t)$}   
            \label{fig:sub-heat_p3_eps400_error}
        \end{subfigure}
        \hfill
        \begin{subfigure}[b]{0.475\textwidth}   
            \centering 
            \includegraphics[width=0.9\textwidth]{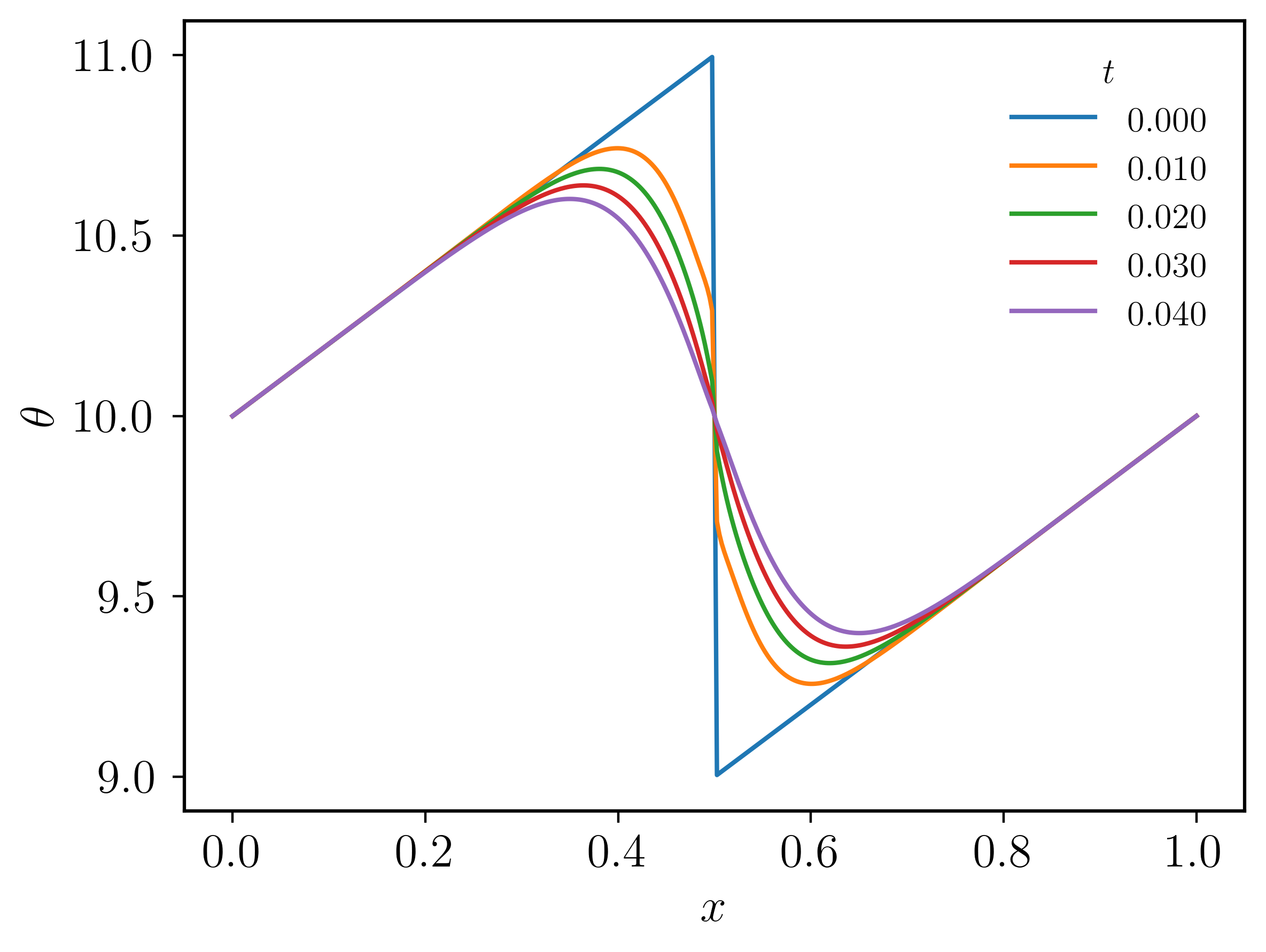}
            \caption[]%
            {Line plots for $\theta(x,t)$}    
            \label{fig:sub-heat_p3_eps400_line_plot}
        \end{subfigure}
        \caption[short 1]
        {Errors and line plots for DtP mapping generated primal field $\theta$ with a smoothed jump in the initial condition.  $\beta=10$. (a) and (b) are obtained on a mesh of $200\times25$. (c) and (d)  are obtained on a mesh of $400\times50$.}
        \label{fig:heat_p3_eps}
    \end{figure}

 Even though Fig.~\ref{fig:heat_p3_disconitnuous_ic} provides a satisfactory result in attempting to capture a discontinuous problem, achieving convergence to \eqref{heat_p3_exact} on a finite mesh through refinement is practically not possible. 
 This is mainly due to the fact that the linear span of continuous FE basis functions can represent a discontinuity only in the limit of infinite mesh refinement. To test our code when convergence on a finite mesh is expected for a sharp transition, we modify the initial condition \ref{eq:heat_p3_initial_condition} to the following: 
 \begin{equation*}
    \theta_0(x) = \begin{cases}
\beta+2x &  \mbox{for } 0 \leq x< 0.5-\veps  \\ 
kx+c    &  \mbox{for } 0.5 - \veps \leq x \leq 0.5+ \veps \\
\beta-2+2x &  \mbox{for }  0.5+\veps < x \leq 1
\end{cases};    \qquad k = \frac{2\veps-1}{\veps}; \qquad c = \beta-\frac{2\veps-1}{2\veps}.
\end{equation*}
With $\beta=10$, the above-stated initial condition tends to \eqref{eq:heat_p3_initial_condition} for $\veps \to 0$.
To obtain the dual results, we set $\veps = 0.01$, $l_l(t)=l_r(t)=0, l_T(x)=0$, $k =0.1$, and perform the simulation on a $200\times25$ and a $400\times50$ mesh up to $T=0.125$. For this case, no data is discarded near the final time. The results for these simulations are shown in Fig.~\ref{fig:heat_p3_eps}. The solution is assessed by comparing it to the analytical result, which is obtained by considering the Fourier series of the modified initial condition, similar to the discontinuous initial condition. This analytical result has been stated in Appendix \ref{app:heat_discontinuous}. Comparing Fig.~\ref{fig:sub-heat_p3_eps200_error} with Fig.~\ref{fig:sub-heat_p3_eps400_error}, it is evident that refinement improves the results.

\begin{figure}[!ht]
\centering
\begin{subfigure}{.32\textwidth}
  \centering
  % include first image
\includegraphics[width=0.95\linewidth]{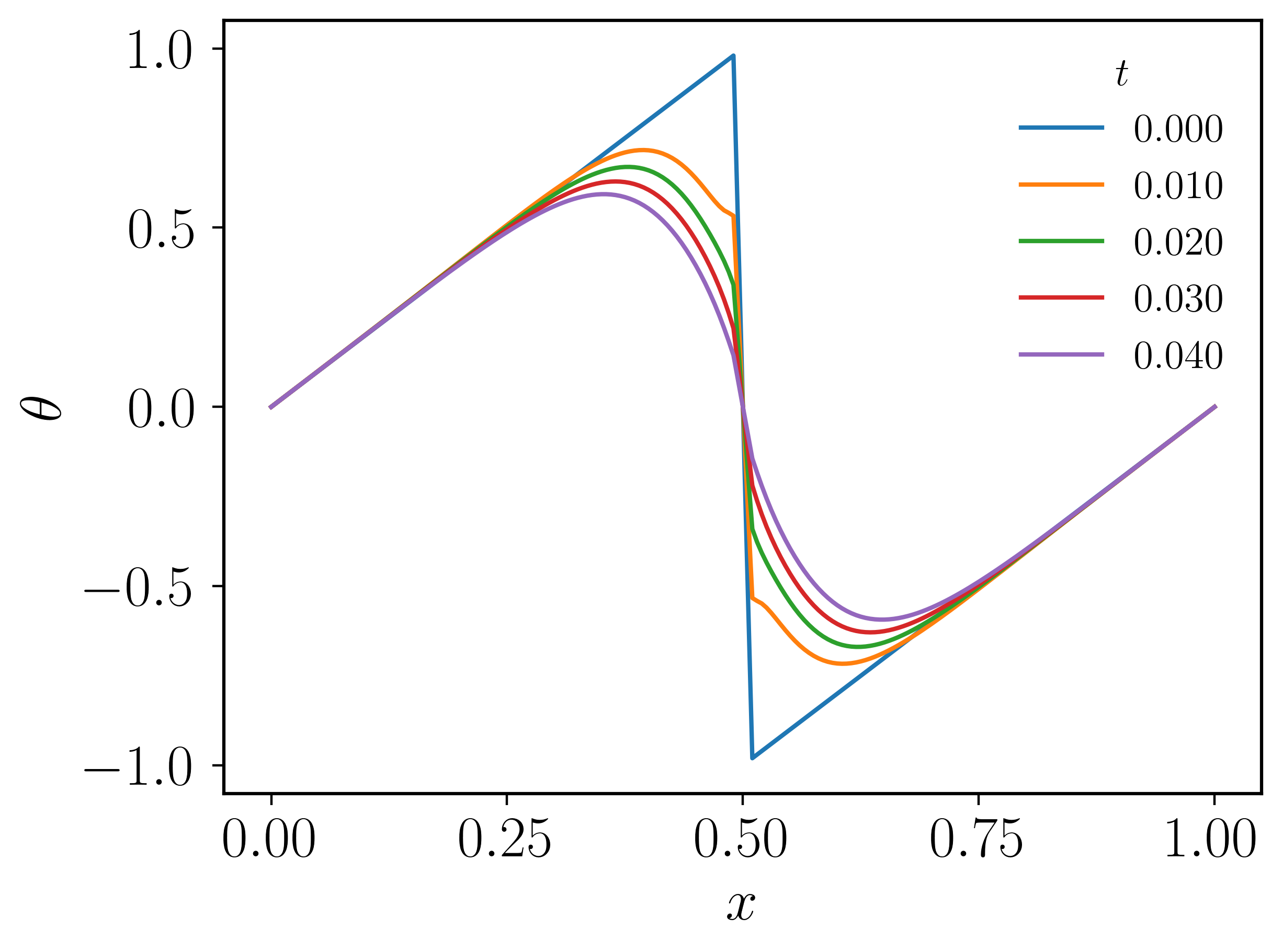}  
  \caption{Line plots for $\theta(x,t)$}
  \label{fig:sub-heat_p3_rms_200_line_plot}
\end{subfigure}
\begin{subfigure}{.32\textwidth}
  \centering
  % include second image
  \includegraphics[width=0.95\linewidth]{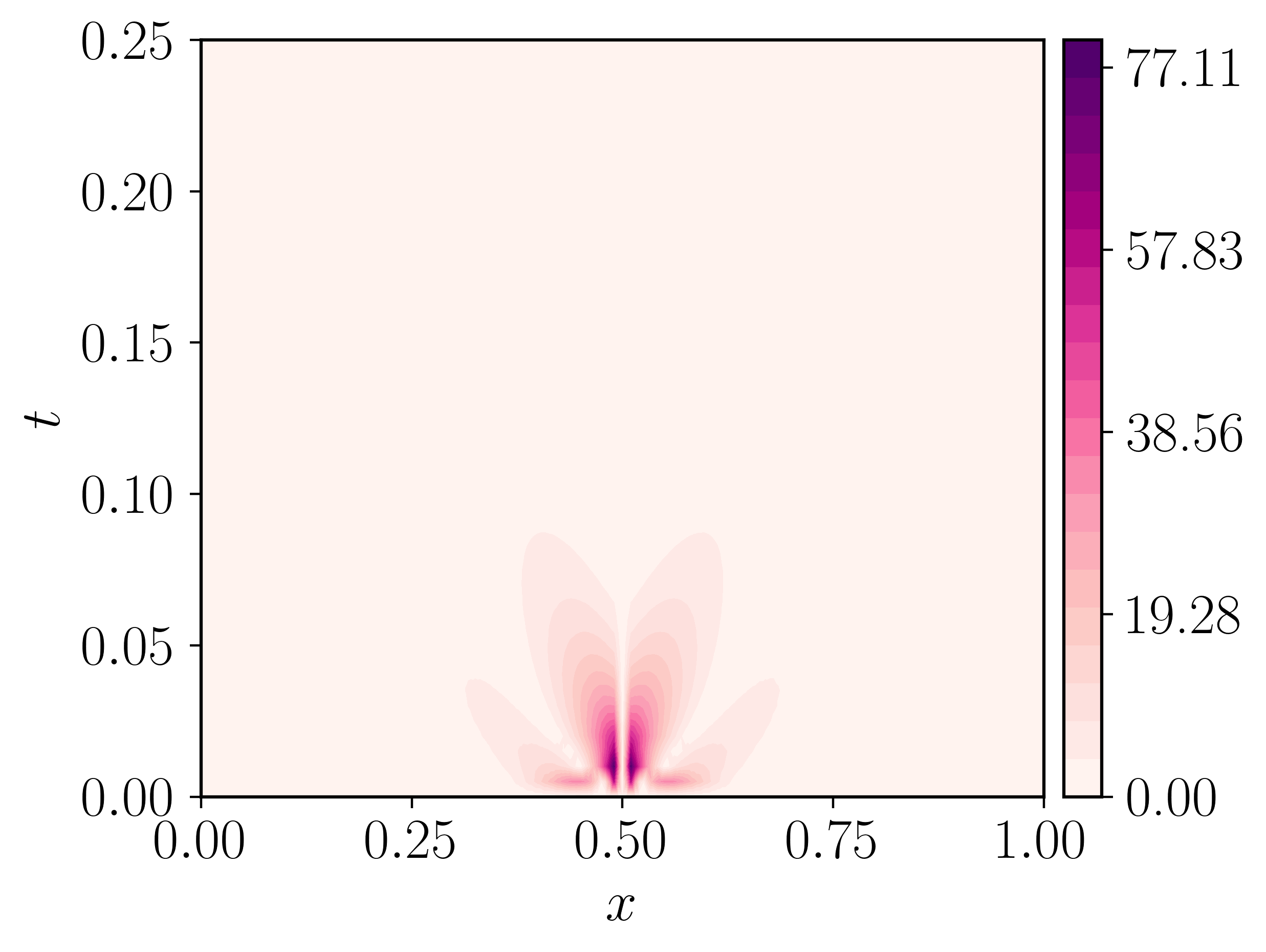}  
  \caption{$\% \,err_1(x,t)$}
  \label{fig:sub-heat_p3_rms_200_error}
\end{subfigure}
\begin{subfigure}{.32\textwidth}
  \centering
  % include second image
  \includegraphics[width=0.95\linewidth]{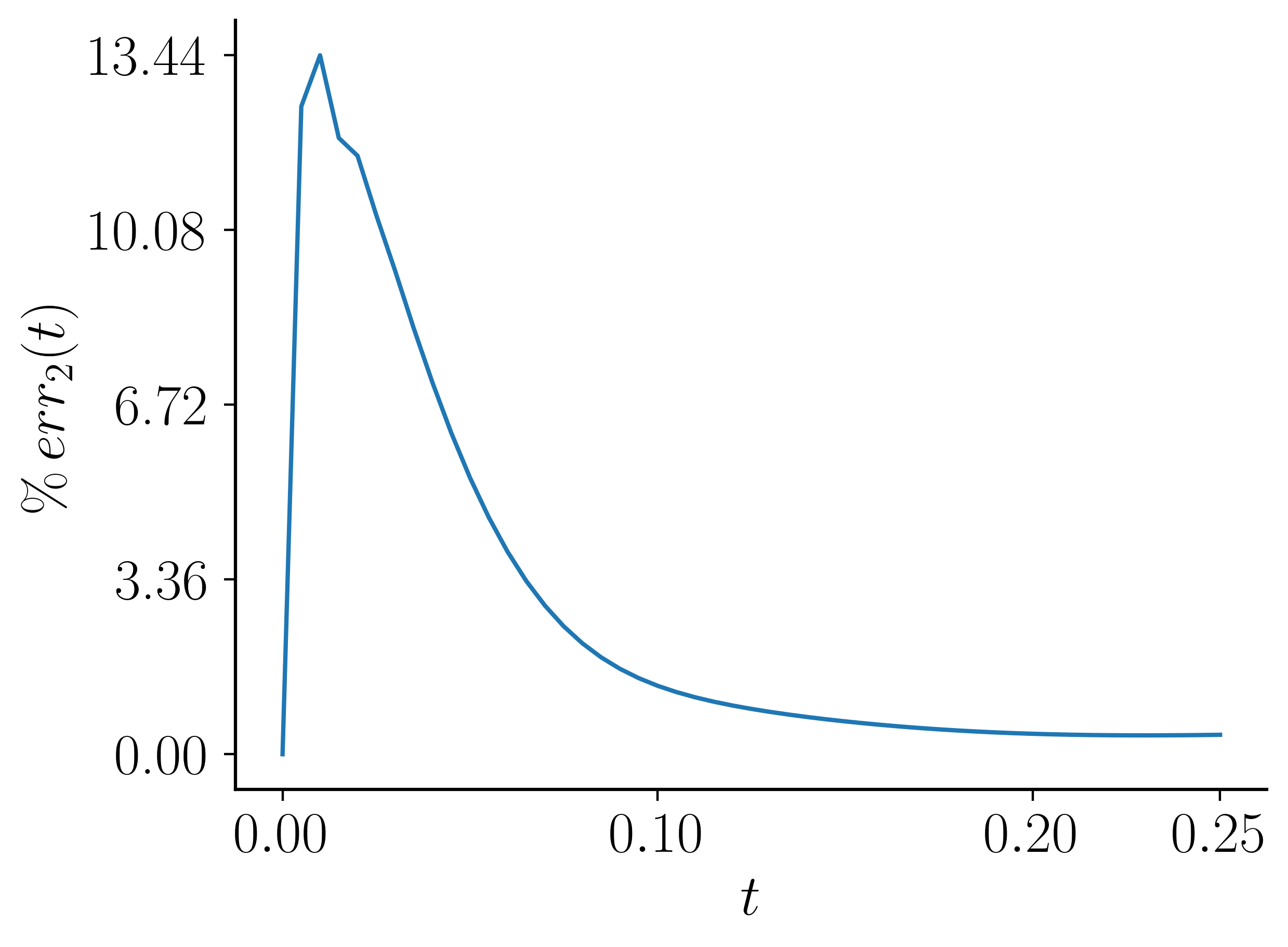}  
  \caption{$\%\,err_2(t)$}
  \label{fig:sub-heat_p3_rms_200_rms}
\end{subfigure}
\caption{Errors and line plots for DtP mapping generated primal field $\theta$ with a smoothed jump in the initial condition. Mesh is  $200\times50$, $\beta=0$, $k=0.1$ and $T=0.25$. }
\label{fig:heat_p3_rms}
\end{figure}

The parameter $\beta$ determines the constant factor by which the function profile is shifted with respect to $x$-axis. By increasing $\beta$, one can artificially create a smaller error field. Furthermore, when $\beta$ is set close to zero, there are regions in the domain where the function value approaches zero, making it difficult to evaluate errors in these regions. To address these issues, we compute the following unbiased error measures ($err_1$ and $err_2$) as well,
\begin{equation*}
    \% \,err_1(x,t) = \frac{|u(x,t) - u^e(x,t)|}{\mbox{rms}(u^e,t)}\times 100; 
    \qquad \% \,err_2(t) = \frac{\mbox{rms}(u-u^e,t)}{\mbox{rms}(u^e,t)}\times 100,
\end{equation*}
where we define rms$(u,t)$ to be the spatial root mean square value of any function given by
\begin{equation*}
    \mbox{rms}(u,t) = \sqrt{\frac{1}{L}\int_0^L\bigl(u(x,t)\bigl)^2\,dx}
\end{equation*}
and $u^e(x,t)$ represents the exact solution which can be obtained using the expression presented in Appendix \ref{app:heat_discontinuous}. $err_1(x,t)$ represents a local measure of error while $err_2(t)$ represents a global measure of error. For $\beta=10$, we also find that the maximum value of $err_2(t)$  decreases from $0.93\%$ to $0.82\%$.
when the mesh is refined from $200\times25$ to $400\times50$. Furthermore, we compute the results for a value of $\beta=0$ for which only $err_1(x,t)$ and $err_2(t)$ are evaluated. The results for this problem along with the error profiles are shown in Fig.~\ref{fig:heat_p3_rms}.

\begin{figure}[!ht]
\centering
\end{figure}

\subsection{Linear transport equation}
% \subsubsection{Example 1: Travelling Sine wave}
% Consider the wave equation \eqref{eq:primal_wave} with the following boundary conditions:
% \begin{equation*}
% u_0(x) = 3 + \sin{5\pi x} \qquad u_l(t) = 3 - \sin{\frac{5\pi t}{2}} \qquad c = 0.5
% \end{equation*}
% The exact solution to this problem is given by
% $$u(x,t) = 3 + \sin{\left(5\pi\left(x-\frac{t}{2}\right)\right)}$$
% Utilizing the weak form \eqref{eq:weak_dual_wave} along with the Dirichlet boundary conditions $\lambda_T(x) = \lambda_r(t) = 0$, we approximate the dual field in $\Omega$. The primal field $u$ along with its error with respect to the exact solution, corresponding to these dual fields, are shown in fig. \ref{fig:wave_p1_primals}.  The time plot for this traveling wave based on the dual formulation is shown in fig. \ref{fig:wave_p1_time_plot}.
% %%
% %%
% \subsubsection{Example: A discontinuous wave}
In order to solve the dual wave equation \eqref{eq:dual_wave_EL}, the weak form \eqref{eq:weak_dual_wave} is employed. Using the discrete approximation
\begin{equation*}
      \lambda(x,t) = \lambda^A N^A(x,t),  
\end{equation*}
where $N$ is the total number of nodes in the FE mesh,
we obtain the following set of equations:
\begin{equation*}\label{eq:weak_dual_wave_discrete}
    K^{AB} \lambda^B = R^A,
\end{equation*}
where
\begin{multline*}
    K^{AB} = - \int_0^T dt \int_0^L dx \, \Bigl( \p_t N^A(x,t)\, (\p_tN^B(x,t) + c \p_xN^B(x,t))  \\ + 
     c\,\p_x N^A(x,t)\, (\p_tN^B(x,t) + c \p_xN^B(x,t)) \Bigl);
\end{multline*}
\begin{equation*} 
     \begin{gathered}
      R^A = 
    \int_0^T dt \, N^A(0,t) u_l(t) +
    \int_0^L dx \, N^A(x,0) u_0(x).
    \end{gathered}
\end{equation*}
Additionally the Dirichlet boundary conditions are given by $$\lambda(L,t) = \lambda_r(t); \qquad \lambda(x,T) = \lambda_T(x).$$
To approximate the integrals in each of the above expressions, we use a two-point Gauss quadrature rule in each direction. Upon obtaining the discretized field $\lambda$, we evaluate the primal field $u$, at the Gauss points of each element and compute its continuous projection as described at the beginning of Sec.~\ref{result}.

% \udk{Additionally to simulate the problem for a long period of time, we utilize the 'time-slicing' strategy explained in the preamble of Sec.~\ref{result}.}

\subsubsection{A Discontinuous Wave}
Consider the wave equation \eqref{eq:primal_wave} with the following boundary conditions  on a domain $\Omega$ with $L=2$:
\begin{equation*}
u_0(x) = \begin{cases}
2 &  \mbox{for } x< 0.2 \\ 
4 &  \mbox{for } x>0.2
\end{cases};  \qquad u_l(t) = 2; \qquad c = 0.25.
\end{equation*}
The exact solution to this problem is given by
\begin{equation} \label{eq:wave_p1_exact}
u(x,t) = \begin{cases}
2 &  \mbox{for } x< 0.2 + 0.25 t \\ 
4 &  \mbox{for } x>0.2+ 0.25 t.
\end{cases}
\end{equation}
We utilize the weak form \eqref{eq:weak_dual_wave} along with the Dirichlet boundary conditions $\lambda_T(x) = \lambda_r(t) = 0$ to compute the dual field in $\Omega$. 

\begin{figure}[h]
\begin{subfigure}{.49\textwidth}
  \centering
  % include first image
  \includegraphics[width=0.9\linewidth]{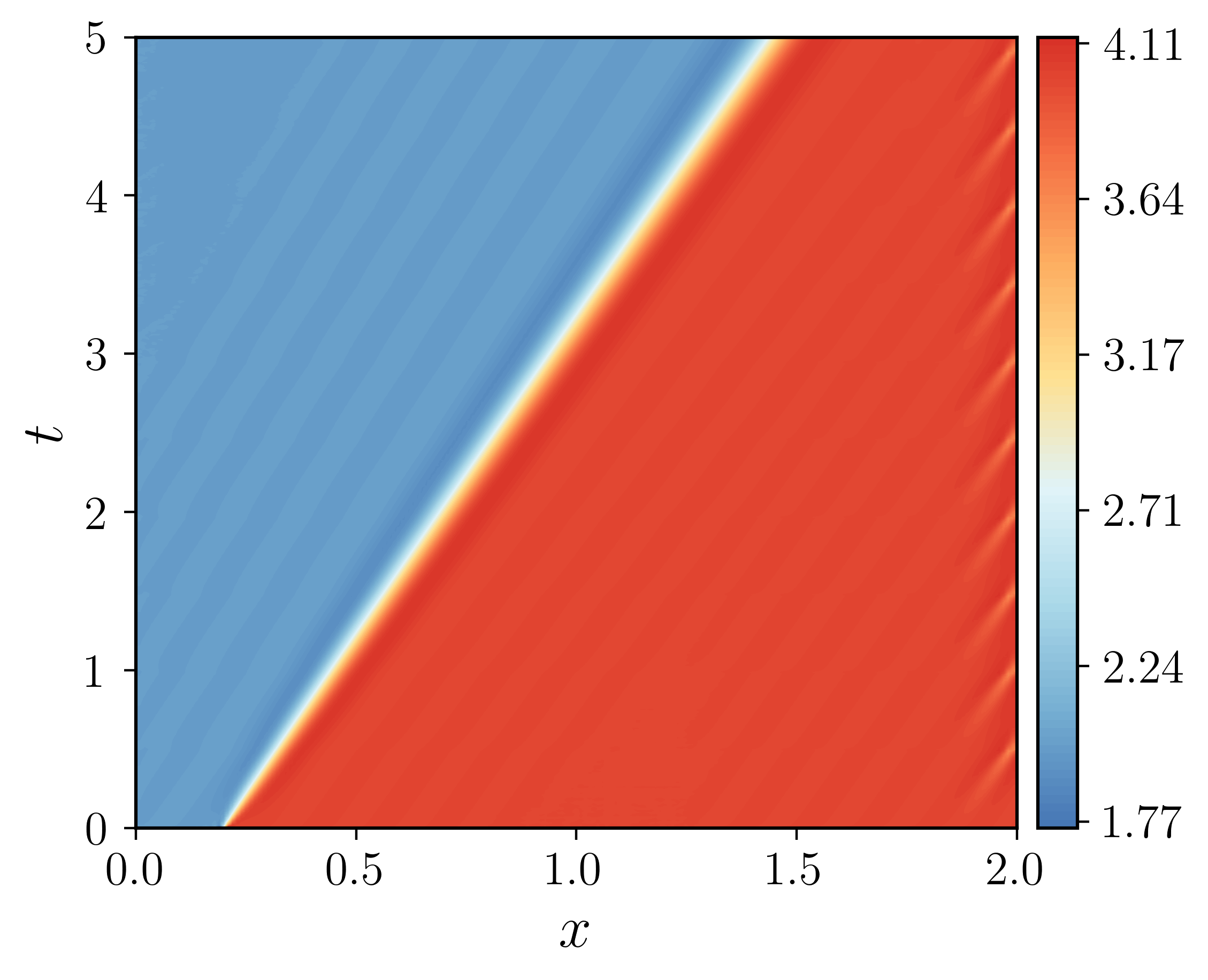}  
  \caption{$u(x,t)$}
  \label{fig:sub-wave_p1_u}
\end{subfigure}
\begin{subfigure}{.49\textwidth}
  \centering
  % include second image
  \includegraphics[width=0.9\linewidth]{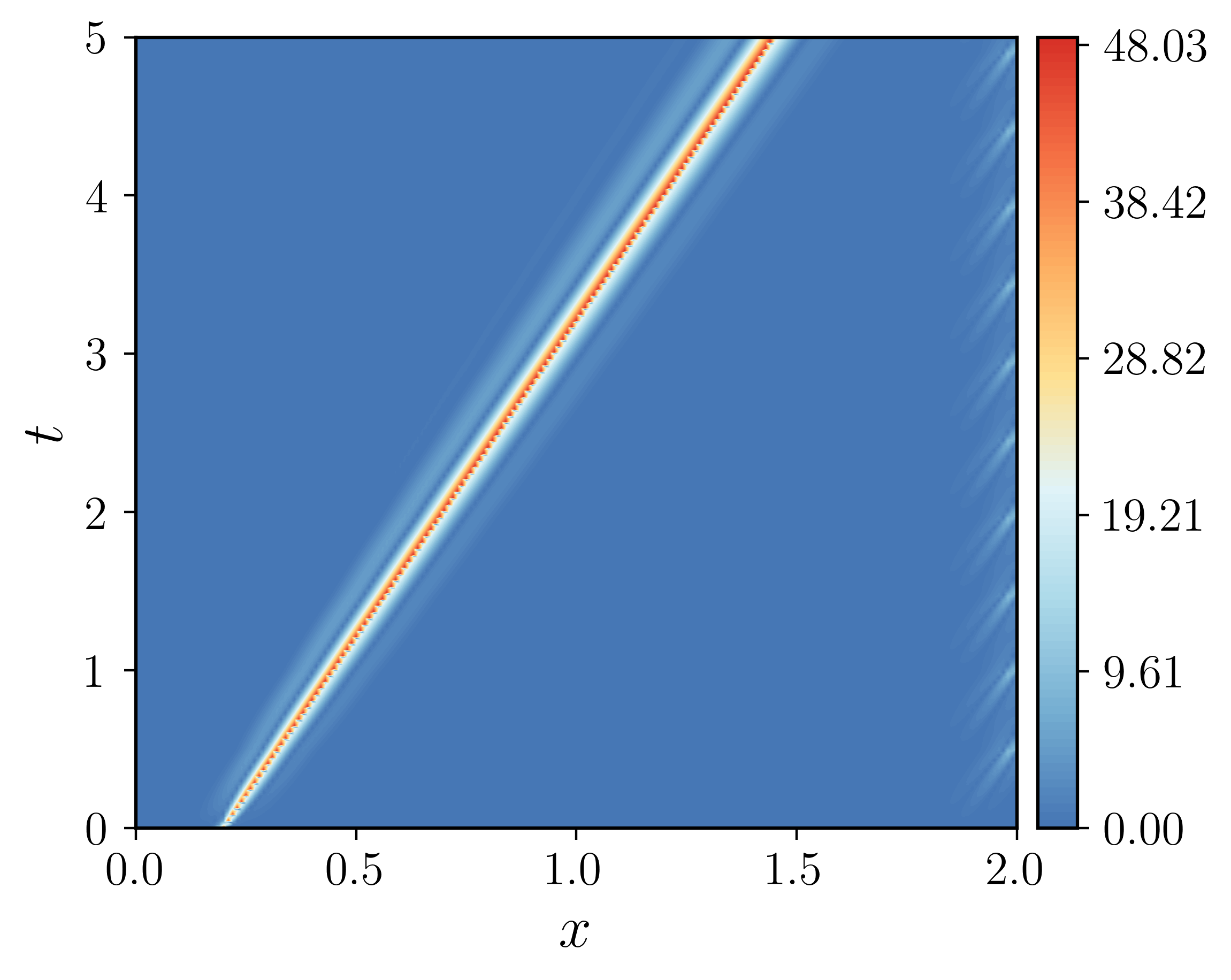}  
  \caption{$\%$ error in $u(x,t)$}
  \label{fig:sub-wave_p1_error}
\end{subfigure}
\caption{DtP mapping generated primal field $u$ for a discontinuous i.c. in wave equation. Mesh is $200\times55$ per stage, $T_e^{(s)}=0.5$, $T=5$ and $c=0.25$.}
\label{fig:wave_p1_primals}
\end{figure}

In order to simulate this problem for a longer period of time, we utilize the `time-slicing' strategy as alluded to in the preamble of Sec.~\ref{result}. We set a time per stage, $T^{(s)}$, and discard the results beyond  $T_e^{(s)}, 0 < T_e^{(s)} < T^{(s)}$. For any stage, $s$, we refer to the dual solution at $(x,t)$ as $\lambda^{(s)}(x,t), t \in [0,T^{(s)}]$, and the corresponding solution for $u$, generated through the DtP mapping, as $u^{(s)}(x,t)$. The global-in-time $u$ solution corresponding to $t \in [0,T^{(s)}]$ is then given by $u\left(x, t + \sum^{s-1}_{i = 1} T^{(i)}_e\right) = u^{(s)}(x,t)$ (with the sum taking the value $0$ for $s = 1$). For any stage, $s$, we set the initial condition, $u^{(s)}(x,0) = u^{(s-1)}\left(x,T_e^{(s-1)}\right) $, with $u^{(1)}(x,0) = u_0(x)$. We repeat this procedure to obtain the results until the simulation reaches a required time, $T$, such that
$$ T = \sum_s T_e^{(s)}.$$

We assess our scheme on this problem by setting $T^{(s)}=0.55$,  $T_e^{(s)}=0.5$ and $T=5$. 
Fig.~\ref{fig:sub-wave_p1_u} and \ref{fig:sub-wave_p1_error} illustrate the field $u(x,t), t \in [0,T]$ and its associated error, respectively.
These figures have been produced by `stitching' together the results for consecutive stages in time, as defined in the previous paragraph.
\begin{figure}
\begin{minipage}{0.5\textwidth}
%% left-hand side: a single subfigure
\begin{subfigure}{\textwidth}
\includegraphics[width=\linewidth]{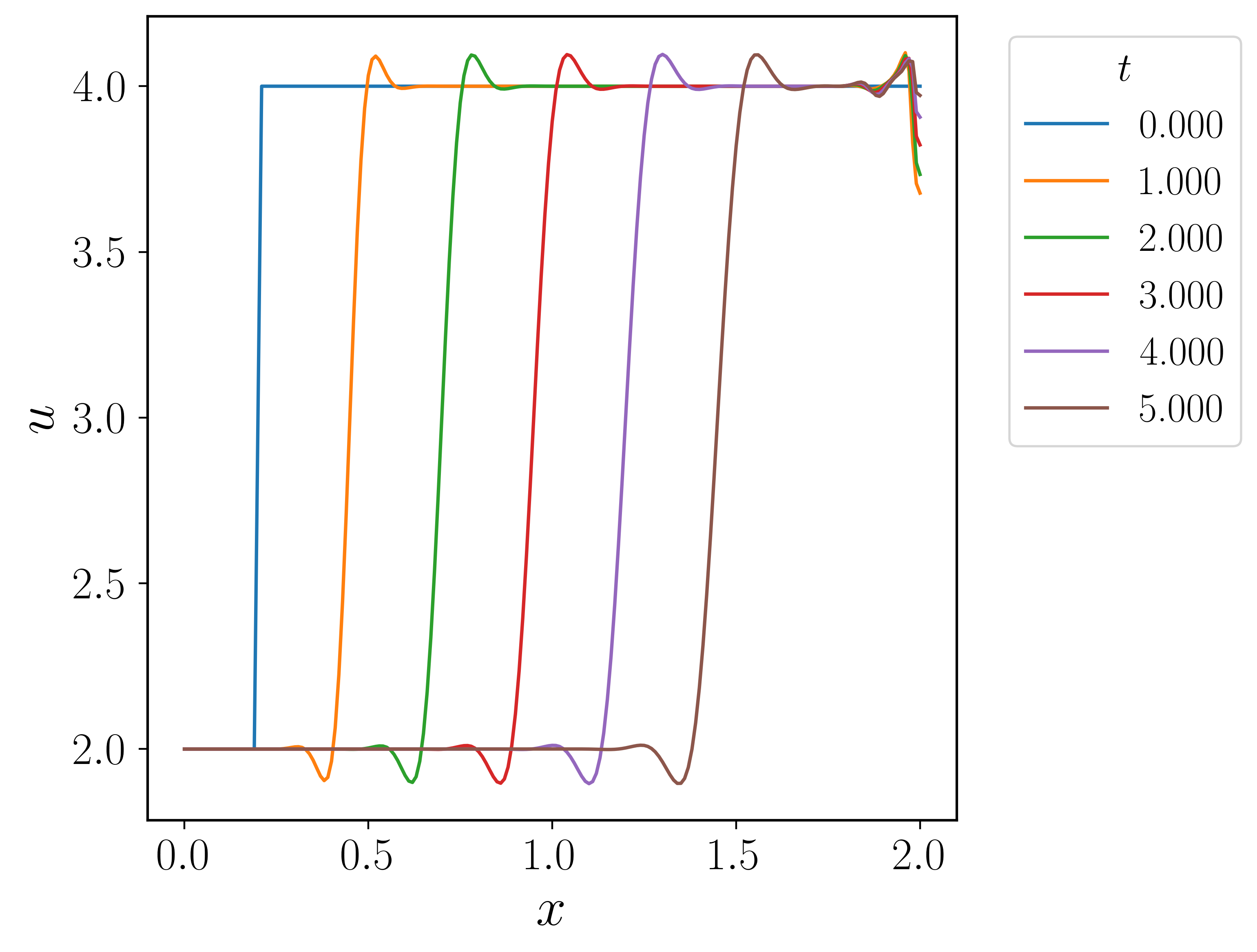}
\caption{Line plot for $u(x,t)$} \label{fig:sub-wave_p1_u_line}
\end{subfigure}
\end{minipage}
%% horizontal separation between the left and right hand sides
\hspace*{\fill}
%% right-hand side: a minipage that contains two more subfigures
\begin{minipage}{0.5\textwidth}
\begin{subfigure}{\linewidth}
\includegraphics[width=0.9\linewidth]{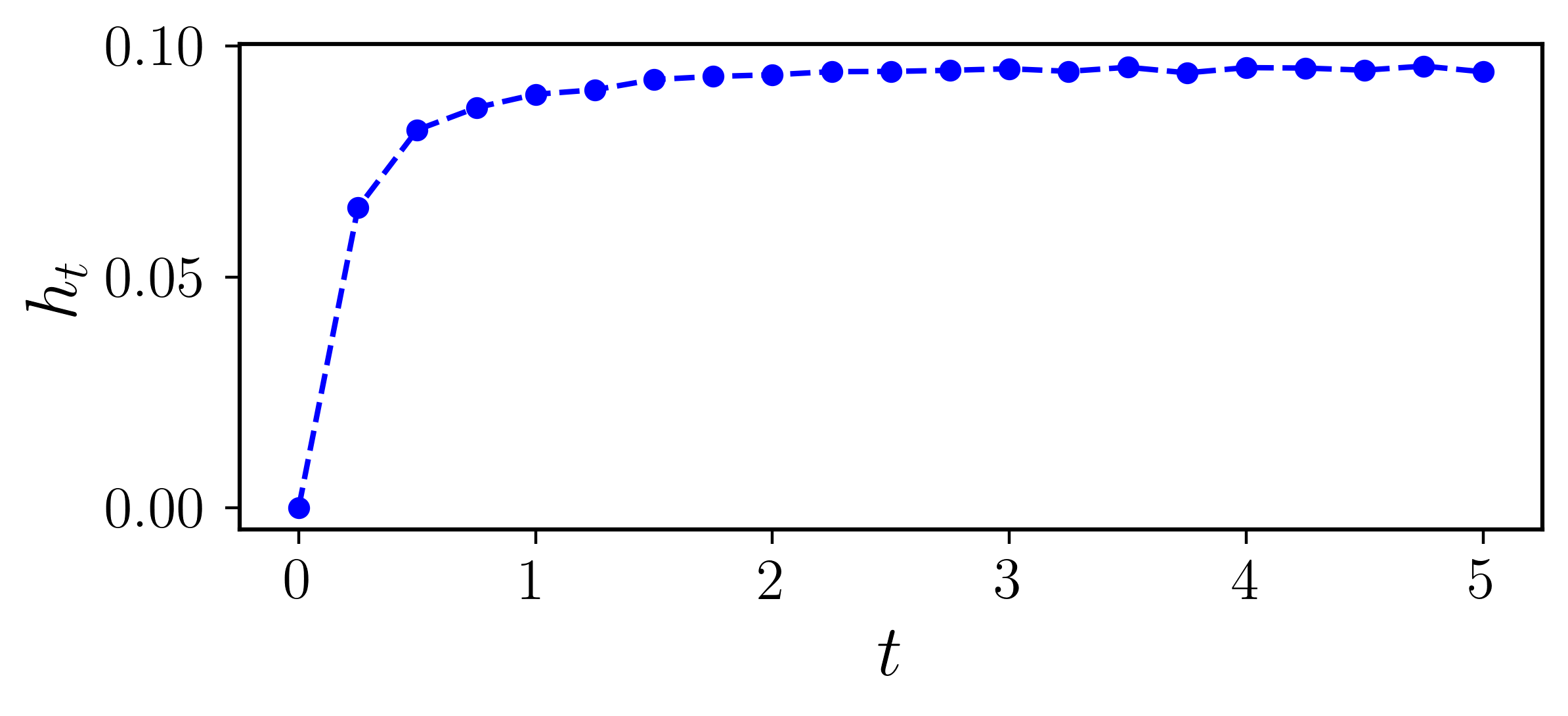}
\caption{Overshoot values} \label{fig:sub-wave_p1_jump_top}
\end{subfigure}

\vspace*{0.6cm}
\begin{subfigure}{\linewidth}
\includegraphics[width=0.9\linewidth]{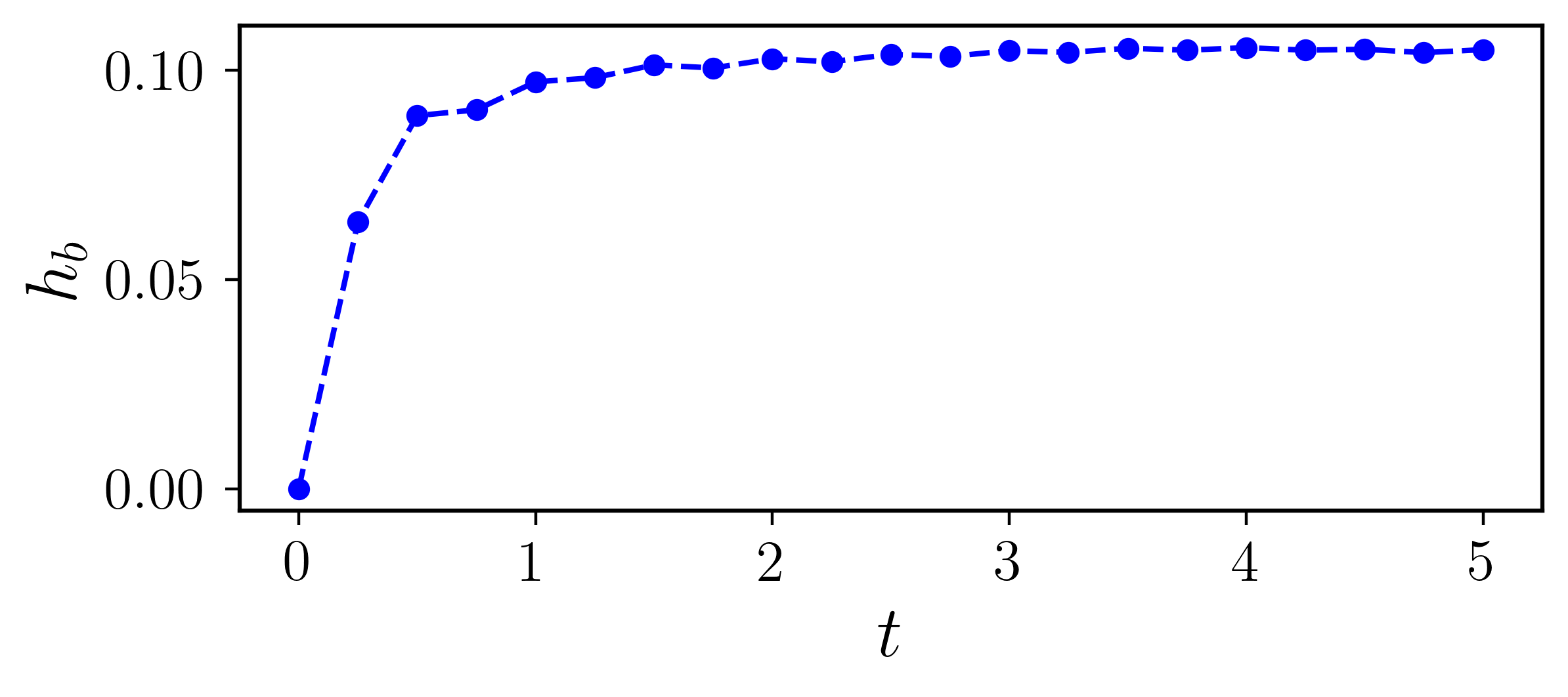}
\caption{Undershoot values} \label{fig:sub-wave_p1_jump_bottom}
\end{subfigure}
\end{minipage}

\caption{(a) illustrates DtP mapping generated primal field $u$ for a discontinuous i.c. in wave equation (Line plots). The error blips near $x=2$ are a consequence of the Dirichlet boundary conditions on $\lambda$. (b) and (c) illustrate the overshooting, $h_t$, and the undershooting, $h_b$, values at the top and the bottom of the jump, respectively.}
\end{figure}
%%%
\begin{figure}[htb]
\begin{subfigure}[t]{.49\textwidth}
  \centering
  % include first image
  \includegraphics[width=\linewidth]{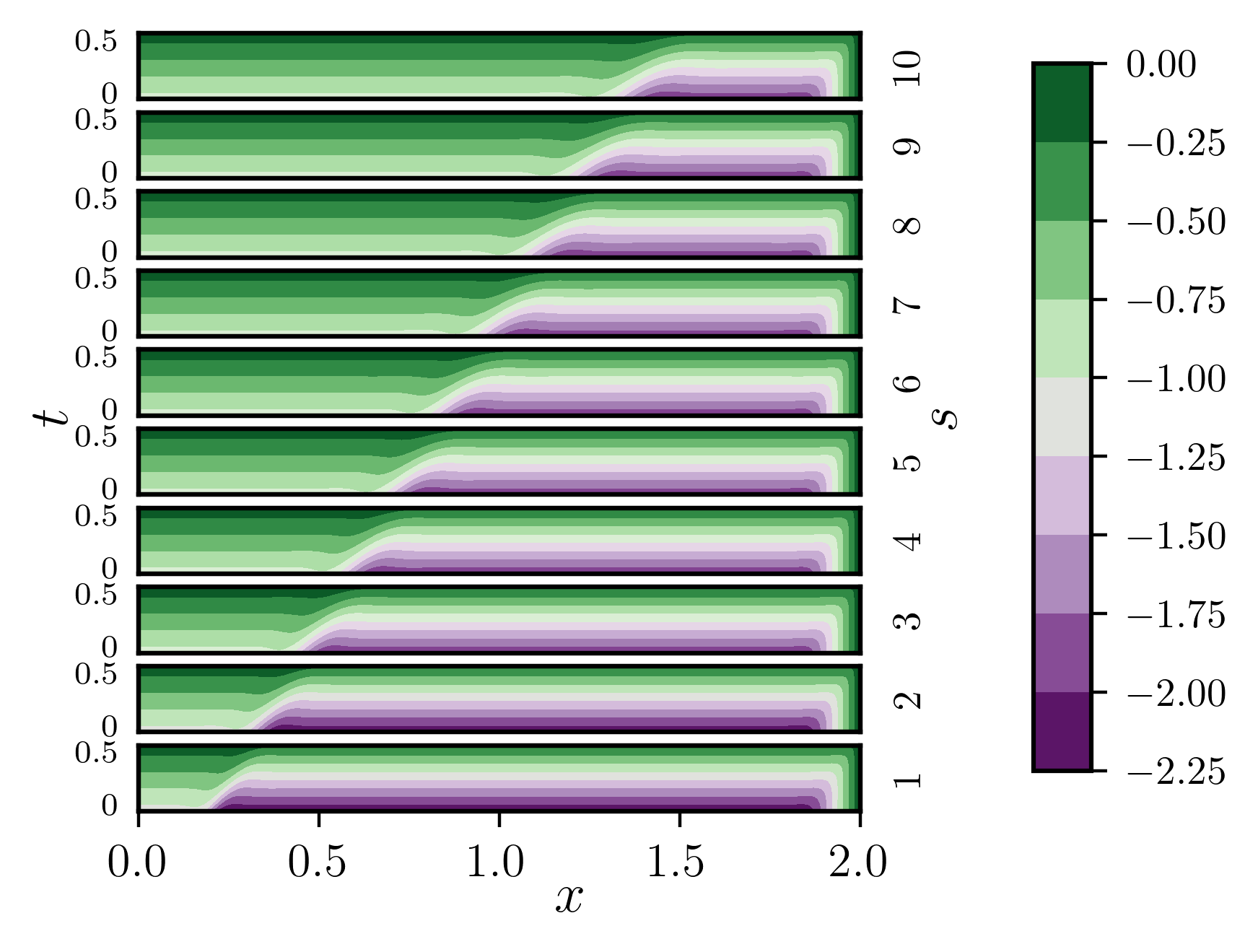}  
  \caption{$\lambda^{(s)}(x,t), t \in [0,T^{(s)}_e], s \in \{1,\ldots, N^{*}\}$, where $N^{*} = 10$ is the total number of computed stages for this problem. }
  \label{fig:sub-wave_p1_lambda}
\end{subfigure}
\begin{subfigure}[t]{.49\textwidth}
  \centering
  % include second image
  \includegraphics[width=\linewidth]{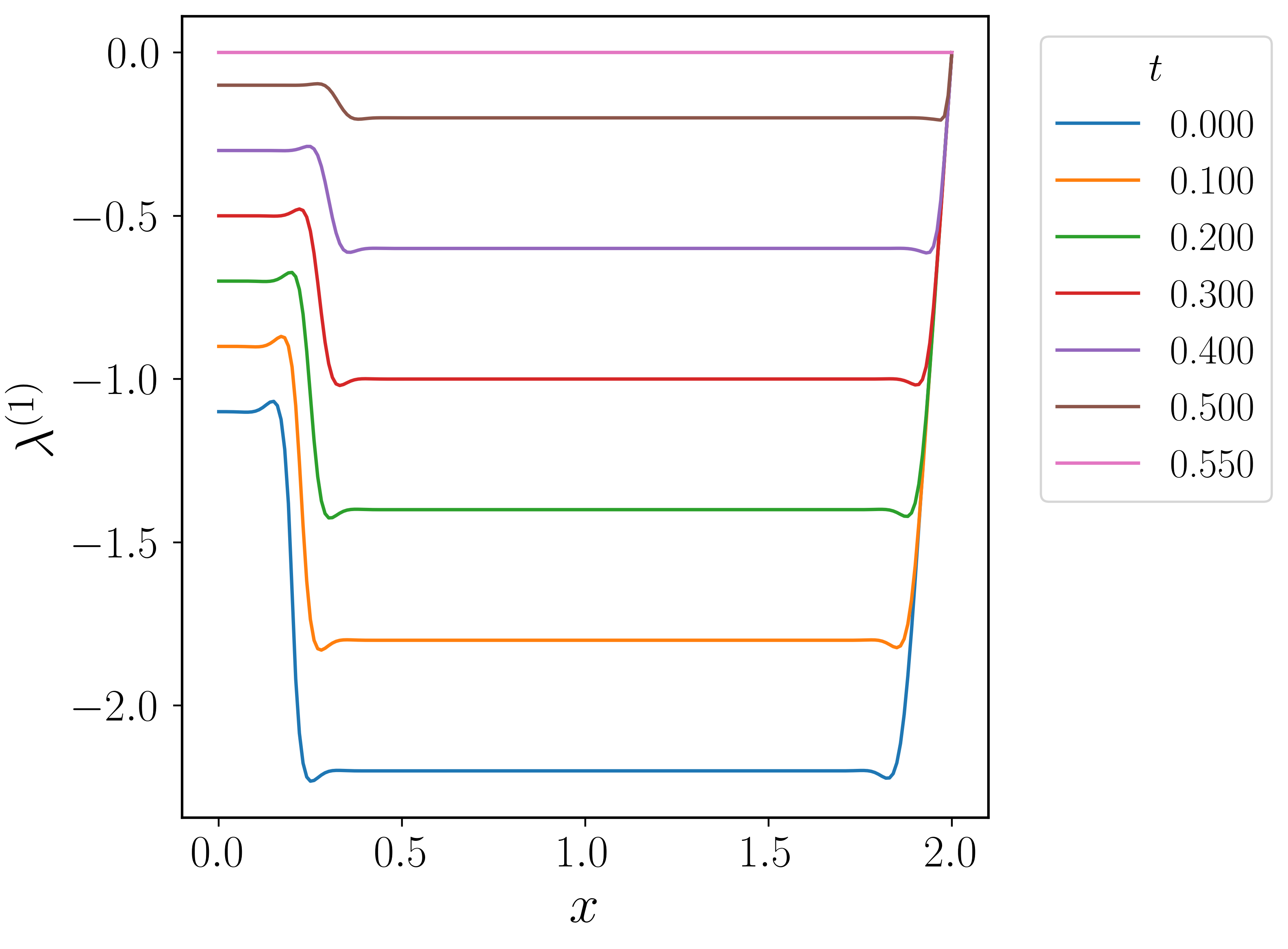}  
  \caption{Line plots for $\lambda^{(1)}(x,t)$}
  \label{fig:sub-wave_p1_lambda_line}
\end{subfigure}
\caption{Dual field $\lambda$ for a discontinuous i.c. in the wave equation. Mesh is $200\times55$ per stage, $T_e^{(s)}=0.5$, $T=5$ and $c=0.25$. For every stage, $\lambda_T^{(s)}(x)=\lambda_r^{(s)}(t)=0$. (b) Illustrates $\lambda^{(1)}$ at different times. }
\label{fig:wave_p1_duals}
\end{figure}

Fig.~\ref{fig:sub-wave_p1_u_line} shows localized errors near the right boundary throughout time. As explained in the preamble to this section, this is a mesh resolution related error arising from the imposition of the b.c. $\lambda_T(x)=\lambda_r(t)=0$ at each stage, resulting in high gradients of the dual $\lambda$ field, as shown in Fig.~\ref{fig:wave_p1_duals}.  The primal field $u$ is a function of the gradients in the $\lambda$ field, and errors in the latter transfer to the primal solution. In this specific problem an exact solution for the dual problem can be generated by the method of characteristics and b.cs for the dual field at the top and right boundaries can be determined to minimize such gradients, given the primal initial and inflow boundary conditions. However, this is not a procedure that generalizes to more involved, e.g. nonlinear, problems and hence we have not adopted such an approach. We have also checked that such errors decrease with mesh refinement.

%%%

The exact solution \eqref{eq:wave_p1_exact} involves a spatial discontinuity persisting in time. While the use of strictly piecewise smooth, globally $C^0$ finite element interpolation for the dual problem allows gradient discontinuities, such discontinuities are not necessarily aligned with the characteristic lines of the primal problem, which results in approximation errors.  
These errors manifest as overshooting and undershooting of the value of $u$ near the jump, as evident in Fig.~\ref{fig:sub-wave_p1_u_line}.
Fig.~\ref{fig:sub-wave_p1_jump_top} and Fig.~\ref{fig:sub-wave_p1_jump_bottom} show that these approximations remain stable as time progresses and, as shown in Fig.~\ref{fig:sub-wave_p1_error}, these errors are localized only in the vicinity of the discontinuity.

In any case, this approximated discontinuity at different times strengthens the claim that our dual scheme can deal with discontinuous functions given that appropriate computational methods are employed.

As a final remark on this problem, we note that while $u^{(s)}(x,0) = u^{(s-1)}\left(x,T_e^{(s-1)}\right)$ by design, $\lambda^{(s)}(x,0) \neq \lambda^{(s-1)}\left(x,T_e^{(s-1)}\right)$, in general. Nevertheless, in any stage $s$, the $u^{(s)}$ and $\lambda^{(s)}$ functions are consistently related by the DtP mapping, exemplifying the many-to-one nature of the DtP change of variables. The latter is a general feature of the duality scheme not particular to this problem (cf.~the distinct dual solutions to the steady heat equation in Sec.~\ref{sec:steady_heat}) which we exploit in our algorithm.

\subsection{Euler's Rigid Body system, with and without damping}\label{sec:euler_algorithm}
To solve the system of equations \eqref{eq:Euler}, we employ the Newton Raphson scheme. Correspondingly, we generate a residual $R[\lambda,\delta \lambda] $ given by
\begin{equation}\label{eq:euler_resi}
    R\,[\lambda;\delta \lambda] = 
        \int_0^T dt \, \sum_i (-I_i \, \hat{\omega}_i \, \dot{ \delta\lambda_i} + c_i \, \hat{\omega}_{i+1} \, \hat{\omega}_{i+2} \, \delta\lambda_i + \nu I_i \, \hat{\omega}_i \, \delta\lambda_i) - \sum_i I_i\,\omega_i^0\,\delta\lambda_i(0),
\end{equation}
where $\delta \lambda_i$ represents the test functions corresponding to each of the equations in \eqref{eq:Euler} and for each $i$ and $j$, $\hat{\omega}_i(\lambda_j,\dot{\lambda}_j) = \omega_i(t)$. Consequently,
the Jacobian matrix is generated by considering the variation of the residual \eqref{eq:euler_resi} in the direction $d\lambda$ and is given by
\begin{multline*}
J\,[d\lambda;\delta \lambda] =
    \int_0^T dt \, \sum_i  \left(\bigl(-I_i (\dot{ \delta\lambda_i}\bigl) + \nu \, I_i (\delta \lambda_i)\bigl ) \left(\parderiv{\hat{\omega}_i}{\lambda_j}d\lambda_j+\parderiv{\hat{\omega}_i}{\dot{\lambda_j}}d\dot{\lambda_j}\right)\right)  \\+ \int_0^T dt \, \sum_i  c_i (\delta \lambda_i)\, \hat{\omega}_{i+2}  \left(\parderiv{\hat{\omega}_{i+1}}{\lambda_j}d\lambda_j+\parderiv{\hat{\omega}_{i+1}}{\dot{\lambda_j}}d\dot{\lambda_j}\right) 
     \\
     +\int_0^T dt \, \sum_i c_i (\delta \lambda_i) \, \hat{\omega}_{i+1}  \left(\parderiv{\hat{\omega}_{i+2}}{\lambda_j}d\lambda_j+\parderiv{\hat{\omega}_{i+2}}{\dot{\lambda_j}}d\dot{\lambda_j}\right),
\end{multline*}
where $d\lambda_i$ represent the variations in each of the $\lambda_i$ respectively.

Next, we discretize the domain and approximate the dual fields as 
$$\lambda_i(t) = \lambda_i^A N^A(t).$$ The subscript index values are interpreted as modulo 3 below, and the index $A$ runs from $1$ to $N$, where $N$ represents the total number of nodes generated upon discretization of time. Using these approximated fields, the value of $\hat{\omega_i}(t)$ can be evaluated using \eqref{eq:euler_DtP}, which is being reproduced below:
\begin{equation*}
 \omega_i(t) - \tilde{\omega}_i(t) =  \sum_j \mathbb{K}^{-1}_{ij} \left(I_j \dot{\lambda}_j(t) - \nu I_j \lambda_j(t) \right); \qquad \hat{\omega}_i(\lambda_j,\dot{\lambda}_j) = \omega_i(t)
\end{equation*}
with
\begin{equation*}
\mathbb{K} = 
\begin{bmatrix}
a & c_3 \lambda_3 & c_2 \lambda_2 \\
c_3 \lambda_3 & a & c_1 \lambda_1 \\
c_2 \lambda_2 & c_1 \lambda_1 & a
\end{bmatrix}.
\end{equation*}
We utilize a base state $\tilde{\omega_i}$ (see Sec.~\ref{sec:euler_formulation}) in this nonlinear problem. The specific choice of $\tilde{\omega}_i$ is subsequently explained in the discussion of Algorithm \ref{algo:euler_algorithm}.

We utilize the residual \eqref{eq:euler_resi}, along with the approximated dual fields to generate a discrete residual $R_i^A$ 
at each node $A$ corresponding to the degree of freedom pair $(A,i)$ given by
\begin{equation*}\label{eq:euler_resi_discrete}
 R_i^A = \int_0^T dt \, (-I_i \, \hat{\omega}_i \, \dot{N^A}(t) + c_i \, \hat{\omega}_{i+1} \, \hat{\omega}_{i+2} \, N^A(t) + \nu I_i \, \hat{\omega}_i \, N^A(t)) - I_i\,\omega_i^0\,\delta_{i0}),
\end{equation*}
\textit{ no sum on $i$}. 
\begin{table}[!ht]
{\begin{algorithm}[H]
\caption*{\textbf{Algorithm}} 
\begin{algorithmic}
\State \textbf{Initialization}: Set $s=1$, $t_f^{(0)} = 0$ and $\bfomega_f^{(0)}=\boldsymbol{\omega_0}$. Choose a value for $N_e^{(s)}$ and $tol$. 
\\\hrulefill
\State \textbf{$\boldsymbol{s_{th}}$ stage}: 
\begin{enumerate}
    \item \label{start-stage} Set $\bflambda^{(s)(-1)}= \bf0$, $t_i^{(s)} = t_f^{(s-1)}$ and  $\bfomega_0^{(s)} = \bfomega_f^{(s-1)}$. Over the domain $\bfOmega^{(s)}$, set $\boldsymbol{\tilde{\omega}}^{(s)} = \bfomega_0^{(s)}$
    \item Newton's Method is used to evaluate the current-stage results. Set $i=0$
    \item[] \textbf{For} $i \geq 0$:
    \begin{enumerate}[i]
        \item \label{step-i} Set $\bflambda^{(s)(i)}=\bflambda^{(s)(i-1)}$.
        \item Evaluate $\bfR^{(s)(i)}$ and $\bfJ^{(s)(i)}$.
        \item Solve for $\boldsymbol{d\lambda}^{(i)(s)}$ using \eqref{eq:increment_vector} and set $\bflambda^{(s)(i+1)}=\bflambda^{(s)(i)} + \boldsymbol{d\lambda}^{(s)(i)}$.
        \item Evaluate $d^{(s)(i)} = \max_{(A,j)}\,|(\lambda^A_j)^{(s)(i+1)}-(\lambda^A_j)^{(s)(i)}|$. 
        \item[]\textbf{if} $d^{(s)(i)}<tol$ \textbf{then} go to step \ref{stepout}
        \item[]\textbf{else do} $i=i+1$ and go to step \ref{step-i}   
    \end{enumerate}
    \item \label{stepout} Evaluate the values of $\bfomega$ at the nodes using the projection \eqref{eq:L2_project}. Discard the last $N_c$ elements along with the results obtained on them. The retained nodes and their values establish the solution for the  current stage  $\bfomega^{(s)}$.
    \item Set $\bfomega_f^{(s)}$ and $t_f^{(s)}$ based on the retained values of $\bfomega^{(s)}$.
    \item Set $s=s+1$ and repeat steps 1-4 until $t_f^{(s)} \geq T$
\end{enumerate}
\end{algorithmic}
\end{algorithm}}
\caption{Algorithm to solve Euler's Rigid Body system }
\label{algo:euler_algorithm}
\end{table}

Additionally, $R_i^N=0$ is a consequence of the Dirichlet boundary condition applied as in \eqref{eq:weak_euler}. The discrete version of the Jacobian corresponding to the degree of freedom pair {$(A,i),(B,j)$} is given as follows:
\begin{multline*}
    J_{ij}^{AB} = \int_0^T dt \,  \bigl(-I_i \, \dot{ N^A}(t) + \nu \,I_i \,N^A(t) \bigl) \left(\parderiv{\hat{\omega}_i}{\lambda_j}N^B(t)+\parderiv{\hat{\omega}_i}{\dot{\lambda_j}}\dot{N^B}(t)\right)  \\ 
    + \int_0^T dt \,  c_i \,N^A(t) \,\hat{\omega}_{i+2} \left(\parderiv{\hat{\omega}_{i+1}}{\lambda_j}N^B(t)+\parderiv{\hat{\omega}_{i+1}}{\dot{\lambda_j}}\dot{N^B}(t)\right) \\
    + \int_0^T dt \, c_i \, N^A(t) \,\hat{\omega}_{i+1} \left(\parderiv{\hat{\omega}_{i+2}}{\lambda_j}N^B(t)+\parderiv{\hat{\omega}_{i+2}}{\dot{\lambda_j}}\dot{N^B}(t)\right),
\end{multline*}
\textit{no sum on $i$}, and the derivatives of $\hat{\omega_i}$ can be evaluated by differentiating the DtP mapping \eqref{eq:euler_DtP} with respect to $\lambda$ and $\dot{\lambda}$, respectively. The final expressions for these derivatives are shown below:
\begin{gather*}
    \parderiv{\hat{\omega_i}}{\lambda_k}=  \mathbb{K}^{-1}_{ij} [\bff_k]_j;
    \qquad \parderiv{\hat{\omega_i}}{\dot{\lambda_k}} =  \mathbb{K}^{-1}_{ij} [\bfg_k]_j,
\end{gather*}
where
\begin{gather*}
    \bff_1 = \begin{bmatrix}
        -\nu I_1 \\
        -c_1 \omega_3 \\
        -c_1 \omega_2
    \end{bmatrix}; \qquad
        \bff_2 = \begin{bmatrix}
        -c_2 \omega_3 \\
        -\nu I_2 \\
        -c_2 \omega_1
    \end{bmatrix} ;\qquad
        \bff_3 = \begin{bmatrix}
        -c_3 \omega_2 \\
        -c_3 \omega_1 \\
        -\nu I_3
    \end{bmatrix}; \\
    \bfg_1 = \begin{bmatrix}
        I_1 \\
        0 \\
        0
    \end{bmatrix}; \qquad
        \bfg_2 = \begin{bmatrix}
        0 \\
        I_2 \\
        0
    \end{bmatrix} ;\qquad
        \bfg_3 = \begin{bmatrix}
        0 \\
        0 \\
         I_3
    \end{bmatrix}. \\
% \mathbb{K} = 
% \begin{bmatrix}
% a & c_3 \lambda_3 & c_2 \lambda_2 \\
% c_3 \lambda_3 & a & c_1 \lambda_1 \\
% c_2 \lambda_2 & c_1 \lambda_1 & a
% \end{bmatrix}    
\end{gather*}
Let $(\cdot)^{(n)}$ denote the value of $(\cdot)$ at the $n^{th}$ iterate of a NR solve. Using this notation, at each NR iterate $n$, the corrections can be evaluated as:
\begin{equation*}
J^{AB\,(n)}_{ij} d\lambda^{B\,(n)}_j = -
R^{A\,(n)}_i.
\end{equation*}
The corrected $\lambda_i$ for the next iteration is obtained via: 
\begin{equation}\label{eq:increment_vector}
\bflambda^{(n+1)}=\bflambda^{(n)} + \boldsymbol{d\lambda}^{(n)},
\end{equation}
where $\boldsymbol{\lambda}$ and $\boldsymbol{d\lambda}$ represent collective terms for $\lambda_i^A$ and $d\lambda_i^A$, respectively.

To obtain a solution to the set of equations \eqref{eq:Euler} at a final time $T$, we proceed incrementally by solving the equations on smaller subdomains. Each of these subdomains is referred to as a stage. We select a time $T^{(s)} \leq T$ as the length of each stage. After each stage is completed, we accumulate its time and continue running the simulation until the cumulative time reaches the final time $T$.

We denote the domain for the $s^{th}$ stage as $\Omega^{(s)}$. We further discretize each of the stages and refer to the number of corresponding elements in stage $s$ as $N_e^{(s)}$. Within each stage, the Newton Raphson (NR) algorithm is employed to solve the set of equations \eqref{eq:Euler}. The algorithm employed to achieve this, along with the motivation for a choice of $\tilde\omega$, has been explained in the following paragraph.

An initial guess for the dual field $\bflambda$ is required by the NR algorithm within each stage. It is practical, and more feasible, to make an initial guess for a primal field $\bfomega$ with physical meaning, and infer a dual field $\bflambda$ consistent with it as a solution of the DtP mapping. By \eqref{eq:euler_DtP}, a primal initial guess of $\bfomega^{(s)}(t) := \tilde{\bfomega}^{(s)}(t)$ is consistent with the dual initial guess of $\bflambda^{(s)}(t) = \bfzero$, for $t$ belonging to the $s^{th}$ stage. We choose this initial guess for $\bfomega$ within each stage $s$ as $\bfomega^{(s)}(t)=\tilde{\bfomega}^{(s)}(t) = \bfomega_0^{(s)}$, where $\bfomega_0^{(s)}$ denotes the initial condition of $\bfomega$ in that stage. 
For the first stage, we set the specified initial condition for the primal problem, $\bfomega_0$, as this constant and for all subsequent stages, we set the constant equal to the value of $\bfomega$ obtained at the end of the previous stage. 
The algorithm capturing these ideas has been outlined in Table \ref{algo:euler_algorithm}.

Thus, our overall algorithm in effect designs different dual functionals for each stage, the set of functionals parametrized by the piecewise constant choice of a base state defining the potential, say $H^{(s)}$, for that stage, as well as the initial condition for the stage.

 The following notation is used for the description of the algorithm:
$s$ denotes the active stage. The collective vectors, $\bfomega$ and $\boldsymbol{\tilde{\omega}}$, consist of $\omega_i^A$ and the base state values $\tilde{\omega_i}^A$, respectively, where the index $i$ takes values from $1$ to $3$ and  the index $A$ ranges over all the nodes in the current stage. $\bfomega_0^{(s)}$ and $\bfomega_f^{(s)}$ denote the values of $\omega_i$ at the first and the last node under the current stage, respectively. For Newton's method, let $(\cdot)^{(s)(i)}$ denote the value of $(\cdot)$ at a stage $s$ and an iteration number $i$. It is important to note that, typically, each stage generates an error towards the end of the stage. To address this problem, we discard the results obtained for a specific number of elements at the end of each stage, and then start the next stage from the final node of the last retained element. We denote the number of discarded elements as $N_c$. Also, $t_i^{(s)}$ and $t_f^{(s)}$ denote the time at the start and end of any stage. $tol$ represents the tolerance for the iterative solve of Newton Raphson.

\subsubsection{Example 1: Free Rotation}
\begin{figure}[ht]
\centering
\begin{subfigure}{.32\textwidth}
  \centering
  % include first image
  \includegraphics[width=\linewidth]{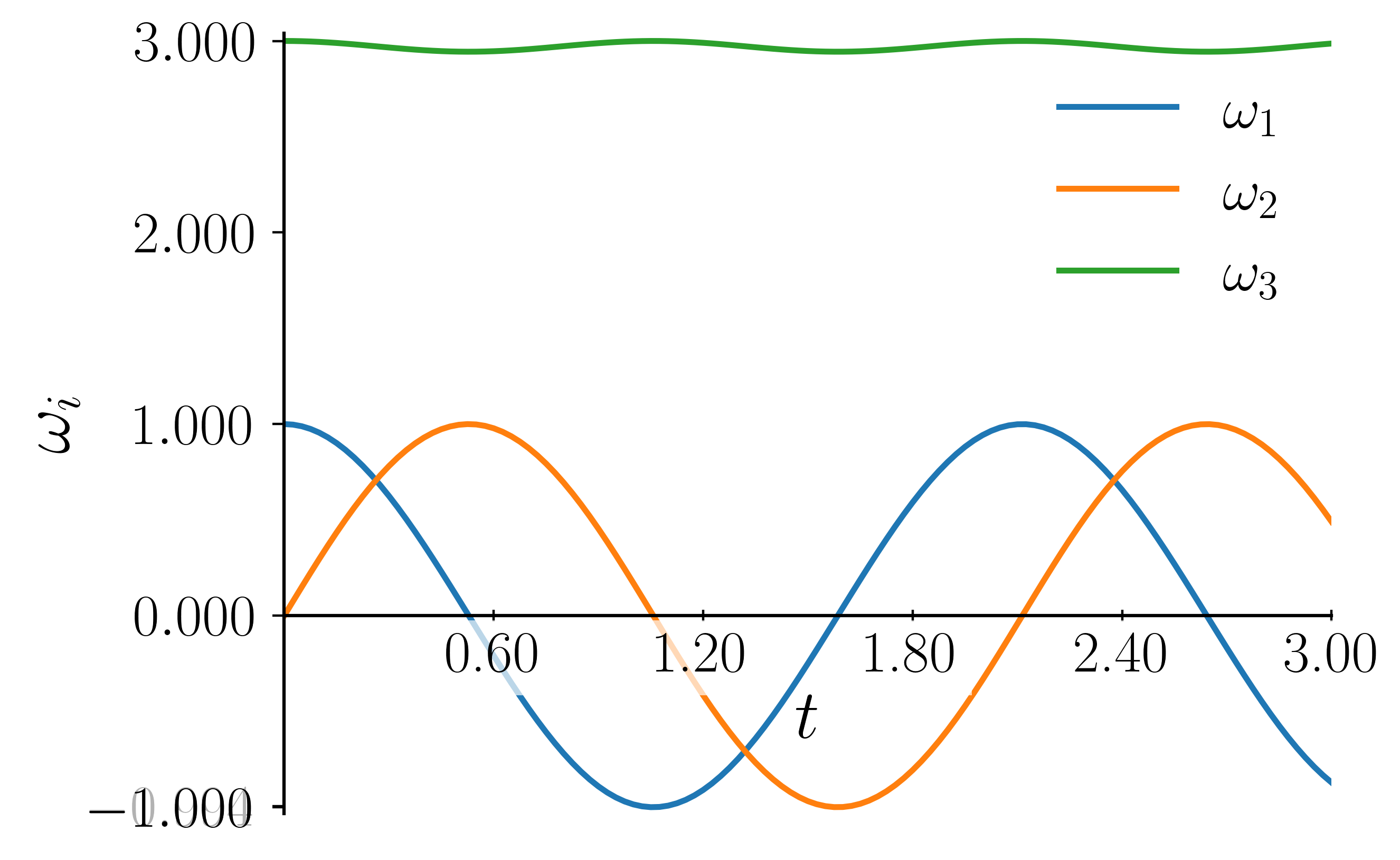}  
  \caption{$\omega_i(t)$}
  \label{fig:sub-euler_p1_omega}
\end{subfigure}
\begin{subfigure}{.32\textwidth}
  \centering
  % include second image
  \includegraphics[width=\linewidth]{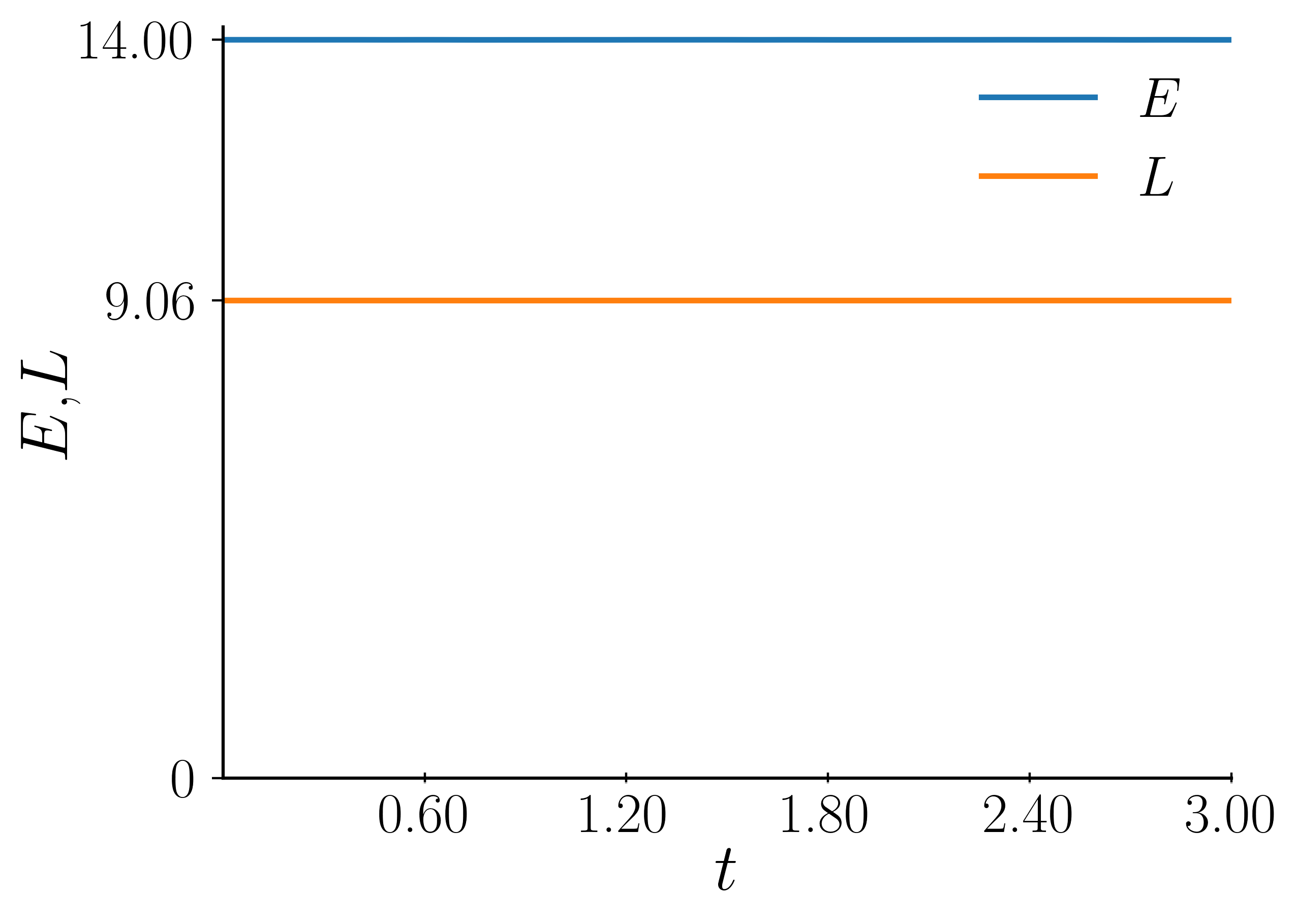}  
  \caption{$E(t)$ and $L(t)$}
  \label{fig:sub-euler_p1_EM}
\end{subfigure}
\begin{subfigure}{.32\textwidth}
  \centering
  % include second image
  \includegraphics[width=\linewidth]{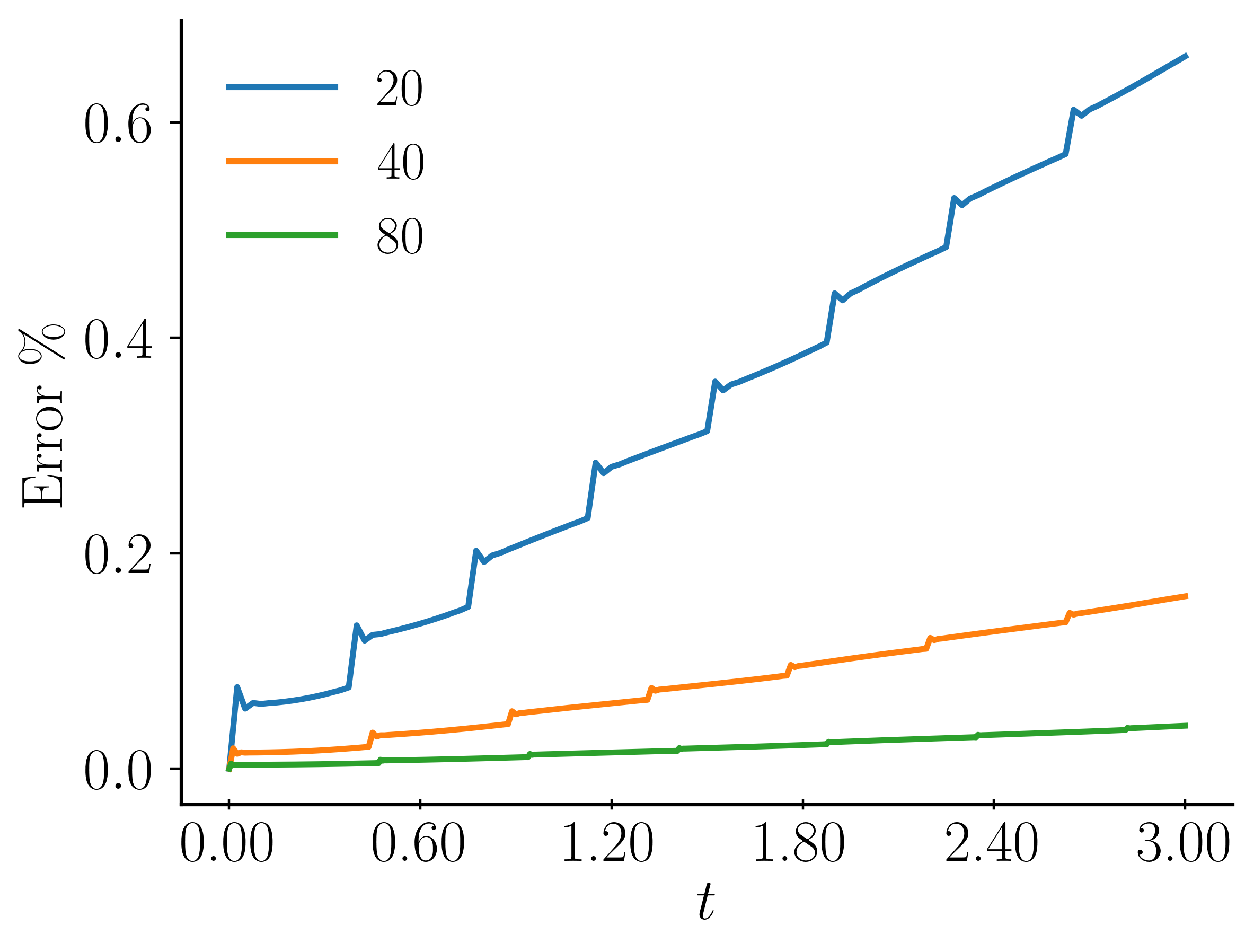}  
  \caption{$err(\omega)(t)$}
  \label{fig:sub-euler_p1_error_norm}
\end{subfigure}
\caption{DtP mapping generated results for a free rotation. Mesh is 20 elements per stage in (a) and (b) and $N_c=5$ for all stages. (c) illustrates error for $\omega_i$ obtained on additional meshes. The maximum error drops quadratically when compared against the mesh size.}
\label{fig:euler_p1_results}
\end{figure}

To demonstrate the algorithm, consider the set of equations \eqref{eq:Euler} with the following conditions:
\begin{equation*}
I = [1,2,3];  \qquad 
\omega_0 = 
    [1,0,3];
\qquad \nu = 0.
\end{equation*}
We run the problem with $T^{(s)}=0.5$ for every stage until $T=3$ is reached. The number of elements considered per stage is 20. The results obtained for $\omega_i$ are plotted in Fig.~\ref{fig:sub-euler_p1_omega}. 

Since this example represents a free rotation, the magnitude of angular momentum and kinetic energy are conserved. The obtained results approximate these features well, as shown in Fig.~\ref{fig:sub-euler_p1_EM}. 
We evaluate the error obtained in $\bfomega(t)$  by comparing it against a reference solution $\bfomega^\bfe(t)$.  
The percentage error for $\bfomega$ (represented by $err(\omega)$) is given by:
$$
err(\omega) = 100 \times\sqrt{\frac{\mathstrut (\omega_1 - \omega^e_1)\strut^2+(\omega_2 - \omega^e_2)\strut^2+(\omega_3 - \omega^e_3)\strut^2}{ (\omega^e_1)\strut^2+(\omega^e_2)\strut^2+(\omega^e_3)\strut^2}}.
$$
For the current example, we set the reference solution to be the exact solution obtained in Appendix \ref{app2}.
Similar to the cases of the heat and wave equation examples, error concentrations appear at the end of every stage.
 To overcome this shortcoming, we discard the results obtained for a few elements (denoted by $N_c$) near the end of the stage as explained at the start of the section \ref{sec:euler_algorithm}.

Additionally, small errors that are relatively larger than the domain errors are found near the beginning of each stage. Predictably, the refinement of the mesh per stage reduces these errors. The error obtained for the current case of 20 elements along with a few different types of meshes ($N_c=5$ for all the stages) are shown in Fig.\ref{fig:sub-euler_p1_error_norm}.

This example demonstrates that the dual methodology has the potential to solve nonlinear problems, preserving conservation properties when expected.

\subsubsection{Example 2: Dissipative System}
\begin{figure}[ht]
\centering
\begin{subfigure}{.49\textwidth}
  \centering
  % include first image
  \includegraphics[width=0.9\linewidth]{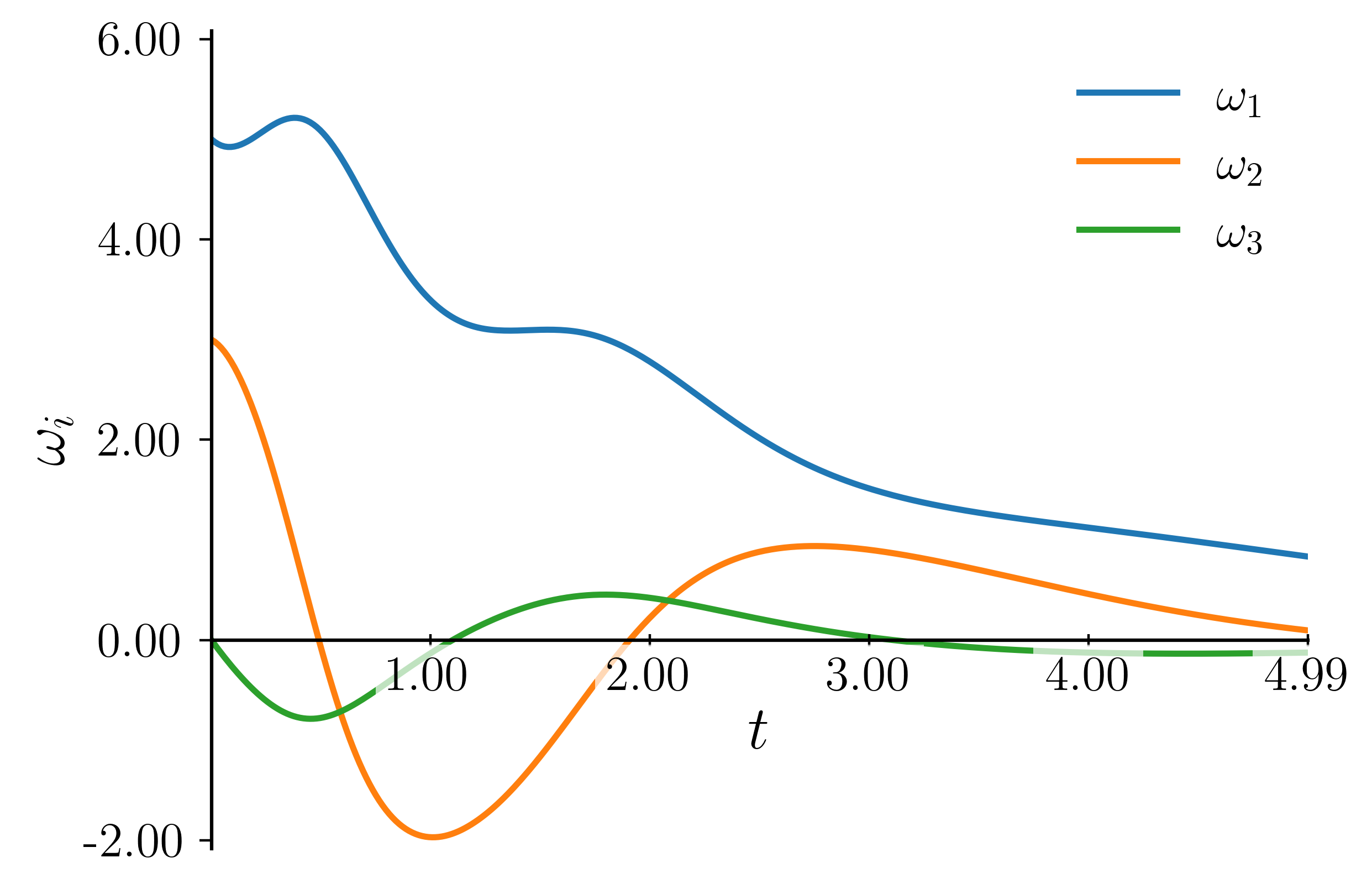}  
  \caption{$\omega_i(t)$}
  \label{fig:sub-euler_p2_omega}
\end{subfigure}
\begin{subfigure}{.49\textwidth}
  \centering
  % include second image
  \includegraphics[width=0.9\linewidth]{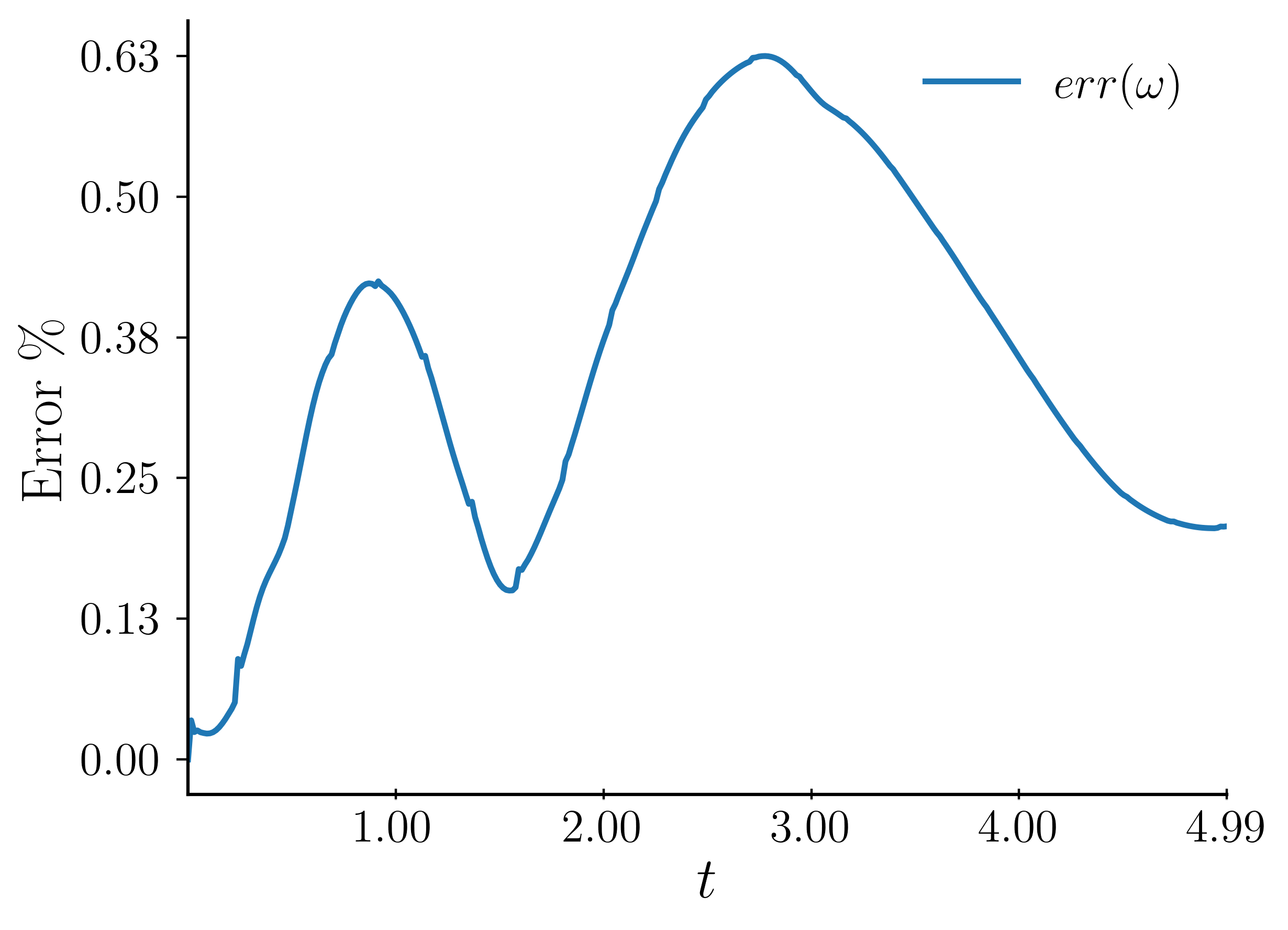}  
  \caption{$err(\omega)(t)$}
  \label{fig:sub-euler_p2_error}
\end{subfigure}
\\
\begin{subfigure}{.49\textwidth}
  \centering
  % include second image
  \includegraphics[width=0.9\linewidth]{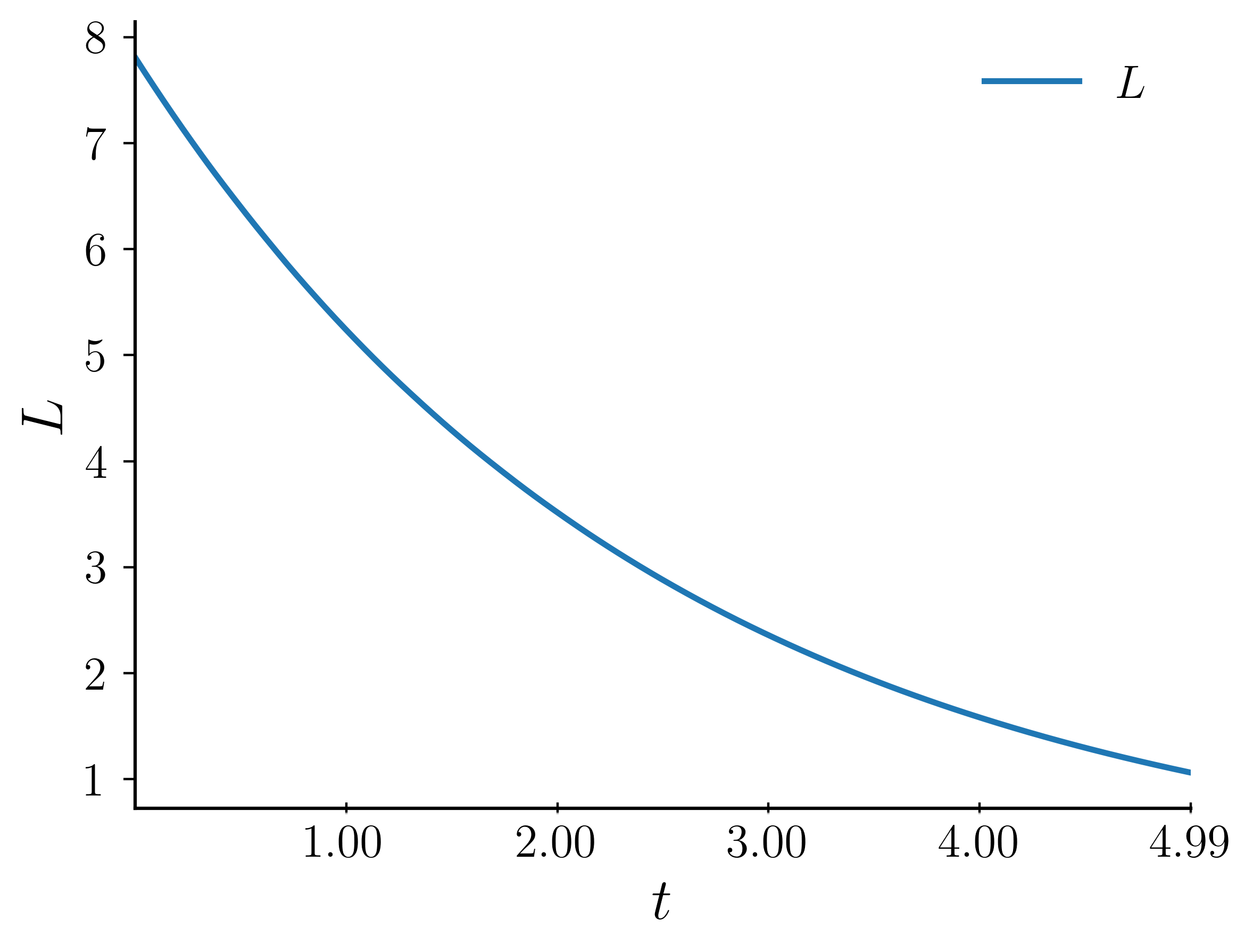}  
  \caption{$L(t)$}
  \label{fig:sub-euler_p2_L}
\end{subfigure}
\begin{subfigure}{.49\textwidth}
  \centering
  % include first image
  \includegraphics[width=0.9\linewidth]{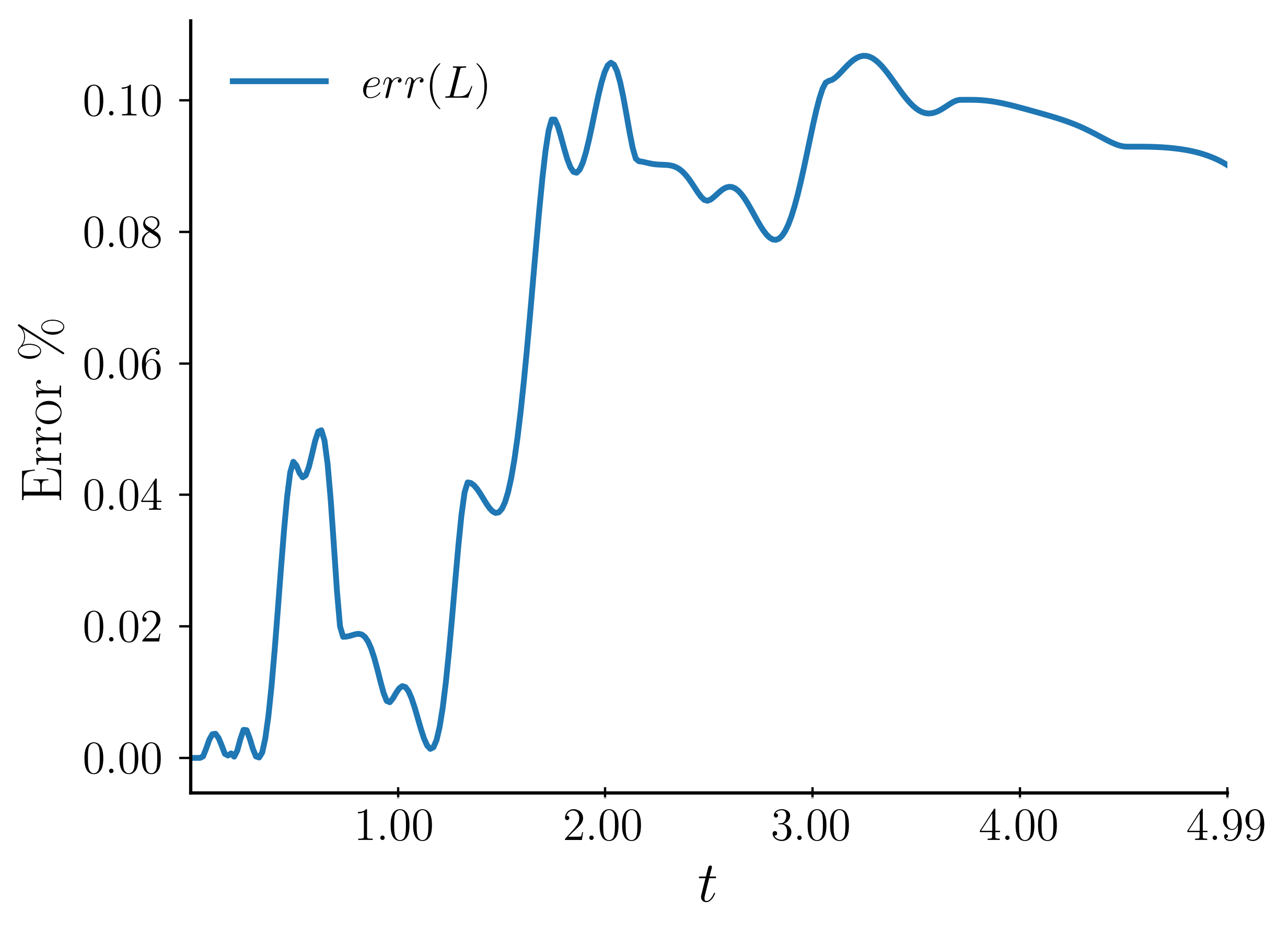}  
  \caption{$\%$ Error in $L(t)$}
  \label{fig:sub-euler_p2_L_error}
\end{subfigure}
\caption{DtP mapping generated results for a dissipative rigid body rotation. The mesh is 20 elements per stage. $N_c=5$, $T^{(s)}=0.3$ per stage and $T=5$.}
\label{fig:euler_p2_omega_error}
\end{figure}

The following example demonstrates the algorithm with an active viscous damping term. Consider the following:
\begin{equation*}
I = 
    [1,2,5];  \qquad 
\omega_0 = [5,3,0];
\qquad \nu = 0.4.
\end{equation*}
These parameter values have been chosen solely for the purpose of demonstrating an example. In fact, the algorithm remains applicable to any value other than the selected values. To evaluate the dual scheme on the current problem, we set the solution obtained from \textit{MATLAB's ode45 tool} with the default settings \cite{MATLAB:R2021a} as the benchmark solution. The plots for $\omega_i$ along with the error obtained for $\bfomega$ against the benchmark are shown in Fig.~\ref{fig:sub-euler_p2_omega} and \ref{fig:sub-euler_p2_error}, respectively. 
Additionally, plots for the magnitude of angular momentum $L$ along with its respective error $err(L)$ when compared against the exact expression obtained in Appendix \ref{app} are produced.
These have been illustrated in Fig.~\ref{fig:sub-euler_p2_L} and Fig.~\ref{fig:sub-euler_p2_L_error}, respectively. 

Surprisingly, these results are based on the Euler-Lagrange equations of the variational principle \eqref{eq:euler_dual_functional}, even though it is believed that, in general, this is not possible for dissipative (P)ODE (in this regard, it is perhaps worth noting that a postulated gradient flow is not an E-L equation of a functional defined over (space)-time fields).
\section{Concluding remarks}
We have provided proof-of-principle demonstrations of a duality-based scheme for generating solutions to differential equations. Non-standard dual formulations for linear parabolic and hyperbolic PDE, and a nonlinear ODE system have been approximated by the finite element method and shown to reproduce correct primal responses. Interestingly, our approach provides a somewhat unifying view of converting initial value problems to well-set degenerate elliptic boundary value problems in (space)-time. This works because the order of the equation in time derivatives (as well as space derivatives when present) increases, allowing for the imposition of an extra, final-time boundary condition on the dual variables which, however, does not get in the way of correct, `free' evolution of the primal field at the considered final time (i.e., constrained only by the primal equation and its initial and boundary conditions), as conjectured in \cite[Sec.~7]{action_2}. A curious observation is that in computationally solving problems with wave-propagation, no explicit account of `upwinding' or domains of dependence of the solution has been involved in the example considered, but much more study to precisely understand the details of such features is required.

Our approach falls within the umbrella of ideas of `Hidden Convexity' in PDE advanced by Yann Brenier \cite{brenier_book}, with some key differences of interpretation. In particular, our understanding and allowance for 
\begin{itemize}
    \item the use of a large class of potentials $H$ in formulating the dual problem, which includes the judicious use of `base states' in designing the dual variational principle for nonlinear problems;
    \item final-time dual boundary conditions as a necessary, routine matter for solving strict primal initial value problems without invoking the degenerate ellipticity of the dual E-L equations for associating a well-set dual variational principle with a primal initial value problem; and
    \item recovery of the primal equations as a critical point or Euler-Lagrange equations of an appropriate dual functional as the principal guiding requirement of the scheme
\end{itemize}
serve as the main conceptual points of departure from Brenier's approach.

\section*{Acknowledgment}
We thank Vladimir Sverak for helpful discussions. This work was supported by the grant NSF OIA-DMR \#2021019 and was also supported by a Simons Pivot Fellowship to AA. It was completed while he was on sabbatical leave at a) the Max Planck Institute for Mathematics in the Sciences in Leipzig, and b) the Hausdorff Institute for Mathematics at the University of Bonn funded by the Deutsche Forschungsgemeinschaft (DFG, German Research Foundation) under Germany's Excellence Strategy – EXC-2047/1 – 390685813, as part of the Trimester Program on Mathematics for Complex Materials. The support and hospitality of both institutions is acknowledged.

\appendix 
\renewcommand{\thesection}{\Alph{section}}

\section{Review of the general formalism}\label{app:rev_duality}
The following two subsections are excerpted from \cite{action_2,action_3} to make this paper self-contained.
\subsection{The essential idea: An optimization problem for an algebraic system of equations}\label{sec:fin_dim}
Consider a generally nonlinear system of algebraic equations in the variables $x \in \mathbb{R}^n$ given by
	\begin{equation}\label{eq:alg_sys}
	    A_\alpha (x) = 0,
	\end{equation}
	where $A: \mathbb{R}^n \to \mathbb{R}^N$ is a given function (a simple example would be $A_\alpha (x) = \bar{A}_{\alpha i} \, x^i - b_\alpha$, $\alpha = 1 \ldots N, i = 1 \ldots n$, where $\bar{A}$ is a constant matrix, \textit{not necessarily symmetric} (when $n = N$), and $b$ is a constant vector). We allow for all possibilities $0 < n \lesseqqgtr N > 0$.
	
	The goal is to construct an objective function whose critical points solve the system \eqref{eq:alg_sys} (when a solution exists) by defining an appropriate $x^* \in \mathbb{R}^n$ satisfying  $A_\alpha (x^*) = 0$.
	
	For this, consider first the auxiliary function
	\begin{equation*}
	    \widehat{S}_H(x,z) = z^\alpha A_\alpha (x) + H(x)
	\end{equation*}
	(where $H$ belongs to a class of scalar-valued function to be defined shortly) and define
	\begin{equation*}
	    S_H(z) = z^\alpha A_\alpha(x_H (z)) + H(x_H(z))
	\end{equation*}
	with the requirement that the system of equations
	\begin{equation}\label{eq:H_fin_dim}
	    z^\alpha \frac{\p A_\alpha}{\p x^i}(x) + \frac{\p H}{\p x^i}(x) = 0
	\end{equation}
	be solvable for the function $x = x_H(z)$ through the choice of $H$, and \textit{any} function $H$ that facilitates such a solution qualifies for the proposed scheme. 
	
	In other words, given a specific $H$, it should be possible to define a function $x_H(z)$ that satisfies 
	\begin{equation*} 
	z^\alpha \p_{x^i} A_\alpha (x_H(z)) + \p_{x^i} H (x_H(z)) = 0 \quad \forall z \in \mathbb{R}^N
	\end{equation*}
	(the domain of the function $x_H$ may accommodate more intricacies, but for now we stick to the simplest possibility). Note that \eqref{eq:H_fin_dim} is a set of $n$ equations in $n$ unknowns regardless of $N$ ($z$ for this argument is a parameter).
	
	Assuming this is possible, we have
	\begin{equation*}
	    \frac{\p S_H}{\p z^\beta} (z) = A_\beta(x_H(z)) +  \left( z^\alpha \frac{\p A_\alpha}{\p x^i}(x_H(z)) + \frac{\p H}{\p x^i}(x_H(z)) \right) \frac{\p x^i_H}{\p z^\beta}(z) = A_\beta(x_H(z)),
	\end{equation*}
	using \eqref{eq:H_fin_dim}. Thus,
	\begin{itemize}
	    \item if $z_0$ is a critical point of the objective function $S_H$ satisfying $\p_{z^\beta} S_H(z_0) = 0$, then the system $A_\alpha(x) = 0$ has a solution defined by $x_H(z_0)$; 
	    \item if the system $A_\alpha(x) = 0$ has a unique solution, say $y$, and if $z^H_0$ is any critical point of $S_H$, then $x_H\left(z^H_0 \right) = y$, for all admissible $H$.
	    \item If $A_\alpha(x) = 0$ has non-unique solutions, but $\p_{z^\beta} S(z) = 0$ ($N$ equations in $N$ unknowns) has a unique solution for a specific choice of the function $z \mapsto x_H(z)$ related to a choice of $H$, then such a choice of $H$ may be considered a selection criterion for imparting uniqueness to the problem $A_\alpha(x) = 0$.
        \item Finally, to see the difference of this approach with the Least-Squares (LS) Method, we note that the optimality condition for the objective $A_\alpha(x) A_\alpha(x)$ is $A_\alpha(x) \p_{x^i} A_\alpha(x) = 0 \centernot \implies A_\alpha(x) = 0$. 
        
        For a linear system $\bar{A} x = b$, the LS governing equations are given by
        \[
        \bar{A}^T \bar{A} z = \bar{A}^T b,
        \]
        with LS solution defined as $z$ even when the original problem $\bar{A} x = b$ does not have a solution (i.e., when $b$ is not in the column space of $\bar{A}$). The LS problem always has a solution, of course. In contrast, in the present duality-based approach with quadratic $H(x) = \frac{1}{2} x^T x$ the governing equation is
        \[
        - \bar{A}\bar{A}^T z = b
        \]
        with solution to $\bar{A} x = b$ given by $x = - \bar{A}^T z$, and the problem has a solution only when $\bar{A} x = b$ has a solution, since the column spaces of the matrices $\bar{A}$ and $\bar{A}\bar{A}^T$ are identical.
        \end{itemize}
 \subsection{The idea behind the general formalism}
The proposed scheme for generating variational principles for nonlinear PDE systems may be abstracted as follows: We first pose the given system of PDE as a \textit{first-order} system (introducing extra fields representing (higher-order) space and time derivatives of the fields of the given system); as before let us denote this collection of primal fields by $U$. `Multiplying' the primal equations by dual Lagrange multiplier fields, the collection denoted by $D$, adding a function $H(U)$, solely in the variables $U$ (the purpose of which, and associated requirements, will be clear shortly), and integrating by parts over the space-time domain, we form a `mixed' functional in the primal and dual fields given by
\begin{equation*}
    \widehat{S}_H [U,D] = \int_{[0,T]\times \Omega} dt d^3x \ \scl_H (\dee,U)
\end{equation*}
where $\dee$ is a collection of local objects in $D$ and at most its first order derivatives. We then require that the family of functions $H$ be such that it allows the definition of a function $U_H(\dee)$ such that
\begin{equation*}
    \frac{\p \scl_H}{\p U} (\dee, U_H(\dee)) = 0
\end{equation*}
so that the \emph{dual} functional, defined solely on the space of the dual fields $D$, given by
\begin{equation*}
    S_H[D] = \int_{[0,T]\times \Omega} dt d^3x \ \scl_H(\dee,U_H(\dee))
\end{equation*}
has the first variation
\begin{equation*}
    \delta S_H = \int_{[0,T]\times \Omega} dt d^3x \ \frac{\p \scl_H}{\p \dee} \delta \dee.
\end{equation*}
By the process of formation of the functional $\widehat{S}_H$, it can then be seen that the (formal) E-L equations arising from $\delta S_H$ have to be the original first-order primal system, with $U$ substituted by $U_H(\dee)$, regardless of the $H$ employed.

Thus, the proposed scheme may be summarized as follows: we wish to pursue the following (local-global) critical point problem
\begin{equation*}
   \begin{smallmatrix} \mbox{extremize}\\ D\end{smallmatrix} \int_{[0,T]\times \Omega} dt d^3x \ \begin{smallmatrix} \mbox{extremize}\\ U\end{smallmatrix} \  \scl_H (\dee(t,x),U),
\end{equation*}
where the pointwise extremization of $\scl_H$ over $U$, for fixed $\dee$, is made possible by the choice of $H$.

Furthermore, assume the Lagrangian $\scl_H$ can be expressed in the form
\begin{equation*}
    \scl_H(\dee, U) := - P(\dee)\cdot U + f(U,D) + H(U)
\end{equation*}
for some function $P$ defined by the structure of the primal first-order system ((linear terms in) first derivatives of $U$ after multiplication by the dual fields and integration by parts always produce such terms), and for some function $f$ which, when non-zero, does not contain any linear dependence in $U$. Our scheme requires the existence of a function $U_H$ defined from `solving $\frac{\p \scl}{\p U} (\dee, U) = 0$ for $U$,' i.e.~$\exists \  U_H(P(\dee),\dee)$ s.t. the equation
\begin{equation*}
    - P(\dee) + \frac{\p f}{\p U}(U_H(P(\dee),\dee), \dee) + \frac{\p H}{\p U}\left(U_H(P(\dee),\dee)\right) = 0
\end{equation*}
is satisfied. This requirement may be understood as follows: define
\begin{equation*}
    f(U, \dee) + H(U) =: M(U, \dee)
\end{equation*}
and assume that it is possible, through the choice of $H$, to make the function $\frac{\p M}{\p U}(U, \dee)$ \textit{monotone} in $U$ so that a function $U_H(P,\dee)$ can be defined that satisfies
\begin{equation*}
    \frac{\p M}{\p U}(U_H(P,\dee), \dee) = P, \quad \forall P.
\end{equation*}
Then the Lagrangian is
\begin{equation*}
    \scl(\dee, U_H(P(\dee),\dee)) = - P(\dee) \cdot U_H(P(\dee),\dee) + M(U_H(P(\dee),\dee), \dee) =: - M^*(P(\dee), \dee)
\end{equation*}
where $M^*(P,\dee)$ is the Legendre transform of the function $M$ w.r.t $U$, with $\dee$ considered as a parameter.

Thus, our scheme may also be interpreted as designing a concrete realization of abstract saddle point problems in optimization theory \cite{rockafellar1974conjugate}, where we exploit the fact that, in the context of `solving' PDE viewed as constraints implemented by Lagrange multipliers to generate an unconstrained problem, there is a good deal of freedom in choosing an objective function(al) to be minimized. We exploit this freedom in choosing the function $H$ to develop dual variational principles corresponding to general systems of PDE.

\section{Appendix: Exact solution of the heat equation with an almost discontinuous initial condition}\label{app:heat_discontinuous}
For given parameters $\veps$ and $\beta$, the solution of the heat equation \eqref{eq:primal_ht}, along with the boundary conditions $\theta_l(t)=\theta_r(t)=\beta$ and the initial condition 
\begin{equation*}
    \theta_0(x) = \begin{cases}
\beta+2x &  \mbox{for } 0 \leq x< 0.5-\veps  \\ 
kx+c    &  \mbox{for } 0.5 - \veps \leq x \leq 0.5+ \veps \\
\beta-2+2x &  \mbox{for }  0.5+\veps < x \leq 1,
\end{cases};    \qquad k = \frac{2\veps-1}{\veps}; \qquad c = \beta-\frac{2\veps-1}{2\veps},
\end{equation*}
is given by:
\begin{equation*}
    \theta(x,t) = \beta + \sum_{m=1}^\infty a_m \sin(2\pi mx)e^{-(2m)^2\pi^2 kt},
\end{equation*}
where 
$$a_m = a_{m_1} + a_{m_2} + a_{m_3};$$
\begin{equation*}
    a_{m_1} = \frac{1}{2\pi^2m^2}\sin(2\pi ml) - \frac{l}{\pi m}\cos(2 \pi ml) + \frac{\beta}{2\pi m}\bigl(1-\cos(2\pi ml)\bigl);
\end{equation*}
\begin{multline*}
    a_{m_2} = \frac{1}{4 \pi^2 m^2}\Bigl(2\pi c m\bigl(\cos(2\pi ml) - \cos(2\pi rm)\bigl) + k\bigl(-2\pi rm\cos(2\pi rm) \\ + \sin(2\pi rm) 
    + 2\pi lm \cos(2\pi lm) - \sin(2\pi lm)\bigl)\Bigl);
\end{multline*}
\begin{equation*}
    a_{m_3} = \frac{1}{2\pi^2m^2}\bigl(-\beta\pi m + (\beta-2+2r)\pi m \cos(2\pi mr) + \sin(2\pi m) - \sin(2\pi mr) \bigl);
\end{equation*}
\begin{equation*}
    l = 0.5 - \veps; \qquad r =0.5+\veps.
\end{equation*}

\section{Appendix: Derivation of Euler's system for rigid body motion with viscosity}\label{app}
Consider the motion of a rigid body with one of its points fixed. In the following part, $\dot{()}$ represents the time derivative $\deriv{()}{t}$. Let $\bfI$ represent the moment of Inertia tensor for the body with $\bfa_i$ and $I_i$ representing its principal directions and principal values, respectively. \textit{The indices presented in this section are modulo 3. Use of summation over repeated indices is followed in this section unless stated otherwise.}. The $\bfa_i$ vectors form an orthonormal triad. Let 
\begin{equation*}
    \bfI = \sum_i I_i \, \bfa_i \otimes \bfa_i.
\end{equation*}
Let $\bfomega$ denote the angular velocity of the body and 
\begin{equation*}
    \bfomega =  \omega_i \, \bfa_i.
\end{equation*}
While $\bfa_i$ and $\omega_i$ vary with time, $I_i$ remain constant due to the rotational invariance of eigenvalues. The angular momentum of the body can be represented by $\bfL = \bfI \bfomega$. Denoting $\bfN$ as the external torque acting on the body, the balance of angular momentum leads to
\begin{equation}\label{eq:angular_mom}
    \frac{d\bfL}{dt} = \frac{d(\bfI\bfomega)}{dt} = \bfN.
\end{equation}
Let $\bfa_i(t) = \bfR(t) \bfa_i(0)$ for some reference basis vectors being fixed at time zero and $\bfR(t)$ represent the rotation tensor of the triad at time $t$. The time derivative of $\bfL$ can be represented in terms of a convected derivative as follows:
\begin{equation}   \label{eq:angular_simplified}
       \frac{d\bm{L}}{dt}  =  \frac{d(L_i\bm{a}_i)}{dt} 
       = \dot{L_i}\bm{a}_i + L_i\dot{\bm{a}_i}
       =  \dot{L_i}\bm{a}_i
       + L_i \dot{\bm{R}} \bm{a}_i(0) 
       = \dot{L_i}\bm{a}_i+ L_i \dot{\bm{R}} \bm{R}^T \bm{a}_i(t).
\end{equation}
Let $\bfW$ be the unique skew tensor corresponding to the angular velocity $\bfomega$ satisfying $\bfomega \times \bfv = \bfW\bfv$ $\forall \bfv \in \mathbb{R}^3$ and $\bfW(t) = \dot{\bfR(t)} \bfR(t)^T$. Using \eqref{eq:angular_mom} and \eqref{eq:angular_simplified},
\begin{equation*}
    \frac{d\bm{L}}{dt} = \dot{L_i}\bm{a}_i 
    + \bm{\omega} \times (L_i \bm{a}_i) =
   \dot{L_i}\bm{a}_i + (\omega_i \bm{a}_i) \times (L_k \bm{a}_k) 
    = \bm{N}.
\end{equation*}
Writing out in components on the $\bfa_i(t)$ basis, and utilizing $\dot{L}_i=I_i\dot{\omega}_i$ (\textit{no sum on i}), Euler’s equations of motion for rigid-body rotation are obtained as
\begin{equation}\label{eq:euler_original}
\begin{aligned}
    I_1 \frac{d\omega_1}{dt} + (I_3 - I_2)\,\omega_2 \,\omega_3 = N_1, \\
    I_2 \frac{d\omega_2}{dt} + (I_1 - I_3)\,\omega_3 \,\omega_1 = N_2, \\
    I_3 \frac{d\omega_3}{dt} + (I_2 - I_1)\,\omega_1 \,\omega_2 = N_3, \\
\end{aligned}
\end{equation}
where $N_i$ are the components of torque on the $\bfa_i(t)$ basis. 

Free rotation of a rigid body about a fixed point is governed by the system \eqref{eq:euler_original} for $N_i = 0$. (Rayleigh) damped rotations of the body about a fixed point can be represented by considering viscous damping torques given by $\bfN= - \nu\bfL, \ \nu > 0$ in \eqref{eq:angular_mom}. In the rotating frame defined by the principal axes of inertia, the governing system reduces to 
\begin{equation}\label{eq:euler_free_equations}
    I_i \dot{\omega_i} + (I_{i+2} - I_{i+1})\,\omega_{i+2} \,\omega_{i+1} + \nu I_i \omega_i = 0,
\end{equation}
\textit{no sum on i},
with the system for free-rotation corresponding to $\nu = 0$.

An analytical solution for the free rotation problem can be found in \cite{Landau-Lifshitz}, and is reproduced in Appendix \ref{app2} for the sake of this paper being self-contained.
The dissipative system given by $\frac{d\bfL}{dt} = - \nu \bfL$ implies that the magnitude of the angular momentum $L = |\bfL|$ evolves as $L(t) = L(0)e^{-\nu t}$, where $L(0)$ is its initial value.

For $\nu = 0$, angular momentum is conserved and $L(t) = L(0)$.

\section{Appendix: Analytical solutions to Free Rotations of a Rigid Body \cite{Landau-Lifshitz}}\label{app2}
For no external torque $\bfN$, angular momentum is conserved (Appendix \ref{app}) and
\[
|\bfL(t)|^2 = L^2(0).
\]
 Similarly, 
\begin{equation*}
    \deriv{\bfL}{t} = 0 \Rightarrow \bfomega\cdot\deriv{(\bfI\bfomega)}{t} = 0\Rightarrow \frac{1}{2}\deriv{}{t}\left(\bfomega\cdot\bfI\bfomega\right) = 0\Rightarrow \frac{1}{2} \bfomega(t)\cdot\bfI(t)\bfomega(t) =: E(t) = E(0),
\end{equation*}
 which represents conservation of kinetic energy.
\textit{In the following, we refer to the constants $L(0), E(0)$ as simply $L, E$ for notational convenience.}

In the rotating frame defined by the principal axes of inertia, the above conditions reduce to:
\begin{equation*}
\begin{aligned}\label{eq:euler_conservation}
L_1^2 + L_2^2 &+ L_3^2 = L^2,\\
    I_1\omega_1^2 +I_2\omega_2^2&+I_3\omega_3^2=2E.
\end{aligned}
\end{equation*}
 Using \eqref{eq:euler_free_equations} with $\nu=0$, $\omega_1$ and $\omega_3$ can be expressed  in terms of $\omega_2$ as
\begin{equation}\label{eq:euler_substitution}
\begin{gathered}
    \omega_1^2 = \frac{(2EI_3-L^2)-I_2(I_3-I_2)\omega_2^2}{I_1(I_3-I_1)}; \qquad
    \omega_3^2 = \frac{(L^2-2EI_1)-I_2(I_2-I_1)\omega_2^2}{I_3(I_3-I_1)}.
\end{gathered}
\end{equation}
Also, when $i=2$ and $\nu=0$, the equation obtained from \eqref{eq:euler_free_equations} is: 
\begin{equation*}
    \left(\deriv{\omega_2}{t}\right)^2 = \left(\frac{(I_3-I_1)\omega_1\omega_3}{I_2} \right)^2.\nonumber   
\end{equation*}
Using the substitution \eqref{eq:euler_substitution} in the last equation results in:
\begin{equation*}
    \left(\deriv{\omega_2}{t}\right)^2   
    = \frac{[(2EI_3-L^2)-I_2(I_3-I_2)\omega_2^2]\times[(L^2-2EI_1)-I_2(I_2-I_1)\omega_2^2]}{I_2^2 I_1I_3}.  
\end{equation*}
Let
\begin{equation}\label{eq:euler_new_variables}
\begin{gathered}
\tau = t\sqrt{\frac{(I_3-I_2)(L^2-2EI_1)}{I_1I_2I_3}}; \qquad
s = \omega_2\sqrt{\frac{I_2(I_3-I_2)}{2EI_3-L^2}}; \qquad
k^2 = \frac{(I_2-I_1)(2EI_3-L^2)}{(I_3-I_2)(L^2-2EI_1)}.
\end{gathered}
\end{equation}
where $k^2<1$ forms a positive parameter for the elliptic function. This leads to the equation
\begin{equation}\label{eq:app2_elliptic_sn}
    \left(\deriv{s}{\tau}\right)^2 = (1-s^2)(1-k^2s^2).
\end{equation}
The value of $\varphi$ satisfying:
$$ \tau = \int_0^\varphi \frac{d\theta}{\sqrt{1-k^2\sin^2(\theta)}}$$ is known as the Jacobi amplitude  and denoted by am$(u,k^2)$. Correspondingly, the Jacobi Elliptic functions are defined as
\begin{equation*}
    \mbox{sn}(\tau,k^2) = \sin(\text{am}(\tau,k^2)); 
    \qquad 
    \mbox{cn}(\tau,k^2) = \cos(\text{am}(\tau,k^2));
    \qquad
    \mbox{dn}(\tau,k^2) = 
    \frac{d}{du}(\text{am}(\tau,k^2)).
\end{equation*}
The Jacobi Elliptic function $\mbox{sn}(\tau,k^2)$ satisfies the equation \eqref{eq:app2_elliptic_sn} exactly. Using \eqref{eq:euler_new_variables}, one can obtain $\omega_2$ explicitly in terms of the elliptic Jacobi functions. Also, $\omega_1$ and $\omega_3$ can algebraically obtained through equations \eqref{eq:euler_substitution}.  The final expressions for each of the $\omega_i$ are shown below:
\begin{equation*}
\begin{aligned}
    \omega_1 &= \sqrt{\frac{2EI_3-L^2}{I_1(I_3-I_1)}}\mbox{cn}\,(\tau,k^2); \quad
    \omega_2 &= \sqrt{\frac{2EI_3-L^2}{I_2(I_3-I_2)}}\mbox{sn}\,(\tau,k^2); \quad
    \omega_3 = \sqrt{\frac{L^2 - 2EI_1}{I_3(I_3 - I_1)}}\mbox{dn}\,(\tau,k^2). 
\end{aligned}
\end{equation*}
\bibliographystyle{IEEEtran}\bibliography{dual_heat_wave_ERB.bib}
\end{document}